\documentclass[fleqn,11pt,oneside]{article}

\usepackage{hyperref}
\usepackage[final]{graphicx}
\usepackage{bbm}  
\usepackage{bm}   
\usepackage{amsmath,amsthm,amssymb,amscd}
\usepackage{caption,subcaption}
\usepackage{booktabs,multirow}
\usepackage{setspace} 
\usepackage[top=1.1in, bottom=1.1in, left=1.6in, right=1.1in]{geometry}
\usepackage[inline,final]{showlabels}
\usepackage{dashbox}
\usepackage{microtype}
\usepackage[medium]{titlesec}
\usepackage{algpseudocode}
\usepackage[ruled]{algorithm}
\usepackage{xpatch}
\usepackage{accents}
\usepackage{bbding}
\usepackage[table]{xcolor}
\usepackage{srcltx}
\usepackage{mathtools}
\usepackage{mathrsfs}

\newtheorem{theorem}{Theorem}[section]
\newtheorem{lemma}[theorem]{Lemma}
\newtheorem{corollary}[theorem]{Corollary}
\newtheorem{definition1}{Definition}[section]
\newtheorem{observe}{Observation}[section]
\newtheorem{remark1}[observe]{Remark}
\newtheorem{example1}{Example}[section]
\newtheorem{aside1}[observe]{Aside}

\newenvironment{definition}[1][]{\begin{definition1}[#1] \rm}{\end{definition1}}

\newenvironment{remark}{\begin{remark1} \rm}{\end{remark1}}

\def\qed{\hfill$\blacksquare$\\} \renewenvironment{proof}{\noindent {\bf 
Proof.}}{\qed}

\algrenewcommand\algorithmicindent{3.0em}

\makeatletter
\xpatchcmd{\algorithmic}{\itemsep\z@}{\itemsep=2.3ex}{}{}
\makeatother

\newif\ifshowboxes \showboxestrue

\renewcommand{\hat}{\widehat}
\renewcommand{\tilde}{\widetilde}
\newcommand{\overbar}[1]{\mkern 1.5mu\overline{\mkern-1.5mu#1\mkern-1.5mu}\mkern%
  1.5mu}

\newcommand{\princip}{\mathrm{p.v.}}

\newcommand{\supp}{\mathop{\mathrm{supp}}}

\newcommand{\pp}[2]{\frac{\partial #1}{\partial #2}}

\newcommand{\erf}{\mathop{\mathrm{erf}}}

\newcommand{\norm}[1]{\ensuremath{ {\lVert #1 \rVert} }}

\def\N{\mathbb{N}}

\def\R{\mathbb{R}}
\def\C{\mathbb{C}}
\def\1{\mathbbm{1}}

\def\cF{\mathcal{F}}
\def\cI{\mathcal{I}}
\def\cL{\mathcal{L}}
\def\cO{\mathcal{O}}
\def\cS{\mathcal{S}}
\def\cT{\mathcal{T}}

\def\sE{\mathscr{E}}
\def\sS{\mathscr{S}}

\usepackage{xspace}
\newcommand{\etal}{\textit{et~al.}\xspace}

\begin{document}

\begin{center}
   \begin{minipage}[t]{5.8in}

We present an adaptive delaminating Levin method for evaluating
bivariate oscillatory integrals over rectangular domains.
Whereas previous analyses of Levin methods impose non-resonance conditions
that exclude stationary and resonance points,
we rigorously establish the existence of a slowly-varying, approximate solution
to the Levin PDE across all frequency regimes, even when the non-resonance
condition is violated.
This allows us to derive error estimates for the numerical solution of the Levin
PDE via the Chebyshev spectral collocation method, and for the evaluation of the
corresponding oscillatory integrals, showing that high accuracy can be achieved
regardless of whether or not stationary and resonance points are present.
We then present a Levin method 
incorporating adaptive subdivision in both two and one
dimensions, as well as delaminating Chebyshev spectral collocation, which
is effective in both the presence and absence of stationary and
resonance points. We demonstrate the effectiveness of
our algorithm with a number of numerical experiments.

\thispagestyle{empty}

  \vspace{ -100.0in}

  \end{minipage}
\end{center}

\vspace{ 2.60in}
\vspace{ 0.50in}

\begin{center}
  \begin{minipage}[t]{4.4in}
    \begin{center}

\textbf{An adaptive delaminating Levin method in two dimensions}\\

  \vspace{ 0.50in}

Shukui Chen$\mbox{}^{\dagger \, \diamond \, *}$,
Kirill Serkh$\mbox{}^{\ddagger \, \dagger \, \diamond}$,\\
James Bremer$\mbox{}^{\ddagger \, \circ}$,
Murdock Aubry$\mbox{}^{\dagger}$\\
June 3, 2025

    \end{center}
  \vspace{ -100.0in}
  \end{minipage}
\end{center}

\vspace{ 2.00in}

\vfill

\noindent 
$\mbox{}^{\diamond}$  This author's work was supported in part by the NSERC
Discovery Grants RGPIN-2020-06022 and DGECR-2020-00356.
\\

\noindent 
$\mbox{}^{\circ}$  This author's work was supported in part by the NSERC
Discovery Grant RGPIN-2021-02613.
\\


\vspace{2mm}

\noindent
$\mbox{}^{\dagger}$ Dept.~of Computer Science, University of Toronto,
Toronto, ON M5S 2E4 

\vspace{2mm}

\noindent
$\mbox{}^{\ddagger}$ Dept.~of Mathematics, University of Toronto,
Toronto, ON M5S 2E4 

\vspace{2mm}

\noindent
$\mbox{}^{*}$ 
Corresponding author. Email: shukui.chen@mail.utoronto.ca \\



\vfill
\eject

\section{Introduction}

The subject of this paper is the efficient evaluation of oscillatory integrals
of the form
\begin{align}
\label{intro:oscint}
I[f,g]
=
\int_{-1}^{1} \int_{-1}^{1}
f(x,y) \exp(i g(x,y))
\, \mathrm{d}x \, \mathrm{d}y
,
\end{align}
where $f \colon [-1,1]^2 \to \C$ and $g \colon [-1,1]^2 \to \R$ are smooth,
slowly-varying functions, and $\nabla g$ can be of large magnitude.
We refer to $f$ as the amplitude function, and $g$ as the phase function.

When evaluating the integral $I[f,g]$ using the most straightforward quadrature
methods, adaptive Gauss-Legendre quadrature for example, the incurred cost grows
at best linearly and at worst quadratically with the magnitude of $\nabla g$.
When the magnitude of $\nabla g$ is large, this exorbitant cost renders such
methods impractical, motivating the development of many approaches for
the numerical integration of highly oscillatory functions, including asymptotic
expansions, numerical steepest descent, Filon methods, and Levin methods.
For a review of these methods, we refer to the monograph \cite{Deano2018} and
references therein.
This paper focuses on Levin methods,
which date back to David Levin's seminal work \cite{Levin1982}.
These methods reduce the problem of evaluating an oscillatory
integral to that of solving a differential equation with slowly-varying
coefficients,
whose solution can be computed to a fixed desired accuracy in time independent
of the magnitude of $\nabla g$.

The classical Levin method, originally introduced in \cite{Levin1982}, works by
solving the second-order partial differential equation (PDE)
\begin{align}
\label{intro:levpde}
p_{xy} + i g_y p_x + i g_x p_y + (i g_{xy} - g_x g_y) p = f
\quad \text{in $(-1,1)^2$}
\end{align}
by collocation with the monomial basis, so that the slowly-varying solution $p$
satisfies
\begin{align}
\frac{\partial^2}{\partial x \partial y}
\left(
p(x,y) \exp(i g(x,y))
\right)
=
f(x,y) \exp(i g(x,y))
,
\end{align}
for all $(x,y) \in (-1,1)^2$. The value of the integral $I[f,g]$ is then readily
obtained as
\begin{align}
\begin{split}
I[f,g]
=
&-p(-1, 1) \exp(i g(-1, 1)) - p( 1,-1) \exp(i g( 1,-1))\\
&+p( 1, 1) \exp(i g( 1, 1)) + p(-1,-1) \exp(i g(-1,-1))
.
\end{split}
\end{align}
In order to guarantee the existence of a slowly-varying solution to
(\ref{intro:levpde}), which can then be approximated by a basis of
slowly-varying functions, Levin's method imposes a non-resonance condition
requiring that both $g_x$ and $g_y$ are nonvanishing and that their maximum
absolute values over the domain are large.
This non-resonance condition excludes stationary points of $g$, where both $g_x$
and $g_y$ vanish, as well as resonance points, which are boundary points where
$\nabla g$ aligns with the normal vector.

A generalization of the classical Levin method, applicable to non-polytopal
domains, was later proposed by Olver in \cite{Olver2006}.
Instead of the PDE (\ref{intro:levpde}), Olver's method considers the
vector-valued first-order PDE,
\begin{align}
\label{intro:divpde}
\cL[\bm{p}]
=
\nabla \cdot \bm{p} + i \, \nabla g \cdot \bm{p} = f
\quad\text{in $\Omega$},
\end{align}
which we refer to as the Levin PDE, where $\Omega \subset \R^n$ is a connected,
open, and bounded set with piecewise smooth boundary $\Gamma$.
If $\bm{p} = (p_1, \ldots, p_n)^T \colon \overbar{\Omega} \to \C^n$ is a
solution to (\ref{intro:divpde}), then
\begin{align}
\nabla \cdot (\bm{p}(\bm{x}) \exp(i g(\bm{x})))
=
f(\bm{x}) \exp(i g(\bm{x}))
,
\end{align}
for all $\bm{x} \in \Omega$, and, by the divergence theorem,
\begin{align}
\label{intro:divthm}
\int_{\Omega} f(\bm{x}) \exp(i g(\bm{x})) \, \mathrm{d}A
=
\int_{\Gamma}
\bm{p}(\bm{x}) \cdot \bm{n}(\bm{x}) \exp(i g(\bm{x}))
\, \mathrm{d}\ell
.
\end{align}
For $\Omega = (-1,1)^2$, the value of $I[f,g]$ is given by
\begin{align}
\label{intro:bdryint}
&\hspace*{-3em} \begin{split}
I[f,g]
=
&-\int_{-1}^{1} p_2( x,-1) \exp(i g( x,-1)) \, \mathrm{d}x
 +\int_{-1}^{1} p_1( 1, y) \exp(i g( 1, y)) \, \mathrm{d}y\\
&+\int_{-1}^{1} p_2( x, 1) \exp(i g( x, 1)) \, \mathrm{d}x
 -\int_{-1}^{1} p_1(-1, y) \exp(i g(-1, y)) \, \mathrm{d}y
.
\end{split}
\end{align}
In order for Olver's method to be applicable, it must be possible to both
accurately solve (\ref{intro:divpde}) over the domain, and to efficiently
evaluate the boundary integral in (\ref{intro:divthm}).
To ensure the solvability of (\ref{intro:divpde}) using collocation,
Olver's method eliminates the null space of the differential operator $\cL$
appearing in the left-hand side of the PDE,
which consists of vector fields of the form
\begin{align}
\label{intro:null}
\bm{p}(\bm{x}) = \bm{q}(\bm{x}) \exp(-i g(\bm{x}))
,
\end{align}
where $\nabla \cdot \bm{q}(\bm{x}) = 0$, by imposing a non-resonance condition
analogous to the condition found in the classical Levin method.
Olver's non-resonance condition requires that $\nabla g$ is both of large
magnitude and is nowhere orthogonal to the boundary of the domain.
This non-resonance condition again excludes both stationary points and resonance
points, noting that, in Olver's formulation, resonance points are precisely
stationary points of the phase function in the right-hand side of
(\ref{intro:divthm}).

When the stationary points have locations which are known beforehand,
certain methods, such as the modified Levin-GMRES method \cite{Olver2010} and
the augmented Levin method \cite{Wang2022}, can be used to efficiently evaluate
the corresponding integrals in the univariate case.
However, these methods have yet to be generalized to higher dimensions, and they
are not applicable when the number and locations of the stationary points of the
phase function are unknown \textit{a priori}.
These challenges have led to the long-standing belief that developing efficient
multivariate Levin methods is particularly difficult when stationary and
resonance points are present, and that their absence is essential for the
success of such methods.

It has been observed however that, in the univariate case, Levin methods are
much more effective in the presence of stationary points than was previously
believed.
For univariate oscillatory integrals of the form
\begin{align}
\label{intro:oscint1d}
\int_{-1}^{1} f(x) \exp(i g(x)) \, \mathrm{d}x,
\end{align}
a Chebyshev spectral collocation method was used by Li~\etal in
\cite{Li2008,Li2010} to discretize the Levin ordinary differential
equation (ODE)
\begin{align}
\label{intro:ode}
p' + i g' p = f \quad 
\text{in $(-1,1)$},
\end{align}
whose solution leads to the formula
\begin{align}
\begin{split}
\int_{-1}^{1} f(x) \exp(i g(x)) \, \mathrm{d}x
&=
\int_{-1}^{1}
\frac{\mathrm{d}}{\mathrm{d}x} ( p(x) \exp(i g(x)) )
\, \mathrm{d}x\\
&=
p(1) \exp(i g(1)) - p(-1) \exp(i g(-1))
.
\end{split}
\end{align}
Li~\etal present numerical evidence showing that, when the resulting linear
system is solved via a truncated singular value decomposition, high accuracy is
obtained even in the presence of stationary points.
In fact, a similar Levin collocation method was previously shown to be effective
in Moylan's thesis \cite{Moylan2008}, in which it was combined with an adaptive
subdivision strategy.
Recently, Chen~\etal \cite{Chen2024} proved the existence of a slowly-varying,
approximate solution to the ODE (\ref{intro:ode}) in the presence of stationary
points, and demonstrated that the univariate Levin method,
which discretizes (\ref{intro:ode}) via the Chebyshev spectral collocation
method and solves the resulting linear system using a truncated singular value
decomposition,
can be applied adaptively to accurately and efficiently evaluate the integral
(\ref{intro:oscint1d}).
We claim that, as a result, it is not actually necessary to exclude the presence
of resonance points when computing $I[f,g]$, since the univariate adaptive Levin
method can be used to evaluate the boundary integrals in (\ref{intro:divthm})
efficiently.

In the presence of stationary points, ill-conditioned linear systems are
encountered when solving the PDE (\ref{intro:divpde}) using collocation methods,
since the null space of the differential operator is fully resolved
by any collocation grid dense enough to resolve $f$ and $\nabla g$.
Ashton \cite{Ashton2025} observes that, with appropriate Cauchy data,
the collocation system can be transformed into a well-conditioned problem,
from which a solution $\bm{p}$ to (\ref{intro:divpde}) with
$
\| \nabla \cdot \bm{p} \|_{L^2(\Omega)}^{2}
=
\mathcal{O}(\| |\nabla g| \|_{L^\infty(\Omega)}^{-1/2})
$
can be recovered.
In particular, Ashton studies the PDE
\begin{align}
\label{intro:direction}
\bm{c} \cdot \nabla p + i (\bm{c} \cdot \nabla g) p = f
\quad \text{in $\Omega$},
\end{align}
where $\Omega$ is an $n$-simplex in $\R^n$ and $g$ has a stationary point at
a vertex of $\Omega$.
This PDE can be derived from (\ref{intro:divpde}) by considering solutions
of the form $\bm{p}=\bm{c}p$, where $\bm{c}$ is a nonzero constant and
$p \colon \overbar{\Omega} \to \C$.
Under a transversality condition on $\bm{c}$ and $g$, Ashton constructs initial
data that determines a solution $p$ satisfying
$
\| \bm{c} \cdot \nabla p \|_{L^2(\Omega)}^{2}
=
\mathcal{O}(\| |\nabla g| \|_{L^\infty(\Omega)}^{-1/2})
$,
and uses the initial data to obtain a well-conditioned collocation system.
We note that their solution may have higher order derivatives that grow with
$\| |\nabla g| \|_{L^\infty(\Omega)}$.

Instead of constructing an appropriate initial condition for
(\ref{intro:direction}), Li~\etal~\cite{Li2009} address the problem of
ill-conditioned linear systems using truncated singular value decompositions.
Specifically, they consider the PDE
\begin{align}
\label{intro:delam}
p_y + i g_y p = f
\quad \text{in $(-1,1)^2$},
\end{align}
which is a special case of (\ref{intro:direction}) with
$\bm{c} = (0,1)^T$ and $\Omega = (-1,1)^2$,
and present numerical evidence showing that, when (\ref{intro:delam}) is
discretized using the Chebyshev spectral collocation method and the resulting
linear system is solved using a truncated singular value decomposition,
the integral $I[f,g]$ can be computed much more accurately and efficiently than
the classical Levin method.
The delaminating quadrature method of Li~\etal exploits the block diagonal
structure of the resulting linear system, where each block corresponds to a
univariate fiber of the domain, in order to reduce computational complexity by
solving each block independently.
Furthermore, it is demonstrated numerically in~\cite{Li2009} that their
method can converge rapidly in the presence of stationary and resonance
points, provided that the domain is refined around them.

In this paper, we prove that the Levin PDE admits a slowly-varying, approximate
solution, both in the case where $\nabla g$ is nonvanishing, and in the case
where $\nabla g$ is of small magnitude and $g$ can have stationary points.
We further prove that, when the Levin PDE is discretized using the Chebyshev
spectral collocation method and the resulting linear system is solved
numerically using a truncated singular value decomposition, high accuracy can be
achieved regardless of whether or not $g$ has stationary and resonance points.
Finally, we present an adaptive
delaminating Levin method for evaluating integrals
of the form~(\ref{intro:oscint}), which solves the Levin PDE by collocation and
computes the boundary integrals in~(\ref{intro:bdryint}) using the univariate
adaptive Levin method presented in~\cite{Chen2024}.
Critically, our use of the univariate adaptive Levin method eliminates the
issue of resonances.
We then demonstrate the effectiveness of our algorithm through
a number of numerical experiments.
We also note in passing the great generality of integrals expressible in the
form~(\ref{intro:oscint}), as most oscillators of interest admit phase
function representations that can be efficiently constructed using the
algorithms of
\cite{Bremer2018,Bremer2023,Aubry2023}.

The remainder of this paper is structured as follows.
Section~\ref{section:prelim} reviews the necessary mathematical and numerical
preliminaries.
In Section~\ref{section:leveq}, we prove the existence of a slowly-varying,
approximate solution to the Levin PDE.
In Section~\ref{section:levnum}, we derive error estimates for the numerical
solution of the Levin PDE via the Chebyshev spectral collocation method,
and for the computed value of the integral.
We give a detailed description of the adaptive delaminating Levin method in
Section~\ref{section:algorithm}, and present the results of numerical
experiments demonstrating the properties of our algorithm in
Section~\ref{section:experiments}.
Finally, we provide a few comments regarding this work and future directions
for research in Section~\ref{section:conclusion}.

\section{Preliminaries}

\label{section:prelim}

\subsection{Notation and conventions}

An $n$-tuple of the form $\bm{\alpha} = (\alpha_1, \ldots, \alpha_n) \in \N_0^n$
is called a multi-index of order $|\bm{\alpha}|$, where
\begin{align}
|\bm{\alpha}| = \alpha_1 + \cdots + \alpha_n
.
\end{align}
Given a multi-index $\bm{\alpha} = (\alpha_1, \ldots, \alpha_n) \in \N_0^n$ and
$\bm{x} = (x_1, \ldots, x_n) \in \R^n$,
we define
\begin{align}
\bm{x}^{\bm{\alpha}} = x_1^{\alpha_1} x_2^{\alpha_2} \cdots x_n^{\alpha_n}
,
\end{align}
and
\begin{align}
D^{\bm{\alpha}}u(\bm{x})
=
\frac
{\partial^{|\bm{\alpha}|}u(\bm{x})}
{\partial x_1^{\alpha_1} \cdots \partial x_n^{\alpha_n}}
.
\end{align}
For $u \colon U \subset \R^n \to \C$, the notation
\begin{align}
\pp{u}{x}, \quad \frac{\partial^2 u}{\partial x \partial y}
\end{align}
is used interchangeably with $u_x$, $u_{xy}$ to denote partial derivatives.
For $\bm{u} \colon U \subset \R^n \to \R^m$, we employ the notation $D\bm{u}$
for the Jacobian matrix of $\bm{u}$.

If $X$ is a topological space, we denote by $C(X)$ the space of continuous
complex-valued functions on $X$.
For a locally compact Hausdorff space $X$, we denote the space of complex Radon
measures on $X$ by $M(X)$, and use $|\mu|$ for the total variation of
$\mu \in M(X)$.

For an open subset $U$ of $\R^n$, we denote by $C^k(U)$ the space of $k$-times
continuously differentiable functions, and by $C^k(\overbar{U})$ the subspace of
functions $u$ in $C^k(U)$ such that $D^{\bm{\alpha}}u$ is uniformly continuous
on bounded subsets of $U$, and hence extends continuously to $\overbar{U}$,
for all multi-indices $\bm{\alpha}$ with $|\bm{\alpha}| \leq k$.
We denote by $C_b^k(U)$ the subspace of functions $u$ in $C^k(U)$ such that
$D^{\bm{\alpha}}u$ is bounded for all multi-indices $\bm{\alpha}$ with
$|\bm{\alpha}| \leq k$, and by $C_c^k(U)$ the subspace of functions in $C^k(U)$
with compact support.
We use the notations
$C^\infty(U) = \bigcap_{k=0}^{\infty} C^k(U)$
and
$C^\infty(\overbar{U}) = \bigcap_{k=0}^{\infty} C^k(\overbar{U})$ for the spaces
of infinitely differentiable functions.
We use $C(U)$, $C_b(U)$, and $C_c(U)$ to denote $C^0(U)$, $C_b^0(U)$, and
$C_c^0(U)$, respectively.  For a Lebesgue measurable subset $U$ of $\R^2$,
we denote by $L^\infty(U)$ the Banach space of measurable functions which
are essentially bounded on $U$.

We let $\sS(\R^n)$ denote the space of infinitely differentiable functions
$\varphi$ such that
\begin{align}
\left\| \varphi \right\|_{\bm{\alpha},\bm{\beta}}
=
\sup_{\bm{x} \in \R^n}
\left| \bm{x}^{\bm{\alpha}} D^{\bm{\beta}}\varphi(\bm{x}) \right|
<
\infty
,
\end{align}
for all multi-indices $\bm{\alpha}$, $\bm{\beta}$.
Its topological dual, which is the space of tempered distributions, is denoted
by $\sS'(\R^n)$.
We write $\langle \cT, \varphi \rangle$ for the action of the tempered
distribution $\cT$ on the test function $\varphi \in \sS(\R^n)$.
The support of $\cT \in \sS'(\R^n)$ is defined as the complement of the union of
all open subsets $U$ of $\R^n$ such that $\langle \cT, \varphi \rangle = 0$
whenever $\varphi \in \sS(\R^n)$ is supported in $U$.
We use $\sE'(\R^n)$ to denote the subspace of $\sS'(\R^n)$ consisting of
compactly supported tempered distributions.
We say that a tempered distribution $\cT \in \sS'(\R^n)$ is of order at most $m$
if, for each compact set $K \subset \R^n$, there exists a constant $M_K$
such that
\begin{align}
\label{prelim:order}
|\langle \cT, \varphi \rangle|
\leq
M_K \sup_{|\bm{\alpha}| \leq m}
\sup_{\bm{x} \in K} \left| D^{\bm{\alpha}}\varphi(\bm{x}) \right|
\end{align}
for every function $\varphi \in \sS(\R^n)$ supported in $K$.
The order of a tempered distribution $\cT$ is the least nonnegative integer $m$
with the property above, noting that all tempered distributions are of finite
order.

We use the convention
\begin{align}
\hat{\varphi}(\xi_1,\xi_2)
=
\frac{1}{4 \pi^2} \int_{-\infty}^{\infty} \int_{-\infty}^{\infty}
\varphi(x,y) \exp(-i (\xi_1 x + \xi_2 y))
\, \mathrm{d}x \, \mathrm{d}y
\end{align}
for the Fourier transform of $\varphi \in \sS(\R^2)$, and, accordingly, we have
\begin{align}
\varphi(x,y)
=
\int_{-\infty}^{\infty} \int_{-\infty}^{\infty}
\hat{\varphi}(\xi_1,\xi_2) \exp(i (\xi_1 x + \xi_2 y))
\, \mathrm{d}\xi_1 \, \mathrm{d}\xi_2
.
\end{align}
The Fourier transform is extended to a topological automorphism in $\sS'(\R^2)$
via the formula
$\langle \hat{\cT}, \varphi \rangle = \langle \cT, \hat{\varphi} \rangle$.
We use the operators
\begin{align}
\cF_y[\varphi](x,\xi_2)
=
\frac{1}{2 \pi} \int_{-\infty}^{\infty} \varphi(x,y) \exp(-i \xi_2 y)
\, \mathrm{d}y
,
\end{align}
and
\begin{align}
\cF_{y}^{-1}[\hat{\varphi}](\xi_1,y)
=
\int_{-\infty}^{\infty} \hat{\varphi}(\xi_1,\xi_2) \exp(i \xi_2 y)
\, \mathrm{d}\xi_2
\end{align}
to denote the partial Fourier transform with respect to $y$ and its inverse.
If the Fourier transform of $\varphi$ is supported in $[-c,c]^2$, we say that
$\varphi$ has bandlimit $c$. Note that we do not require $c$ to be the smallest
positive real number with this property.

We use bold capital letters for matrices, and bold lowercase letters for vectors.
For $1 \leq p \leq \infty$,
we denote by $\|\bm{\mathsf{x}}\|_p$ the $p$-norm of a vector
$\bm{\mathsf{x}} \in \C^n$, and by $\|\bm{\mathsf{A}}\|_p$ the matrix norm of
a matrix $\bm{\mathsf{A}} \in \C^{m \times n}$ induced by the vector $p$-norms.

The notation $x \lesssim y$ is used to indicate that there is some constant
$C>0$, not depending on either $x$ or $y$, such that $x \leq Cy$.
For $\bm{\mathsf{x}} \in \C^n$, the notations $\|\bm{\mathsf{x}}\|_p \lesssim y$
and $\|\bm{\mathsf{x}}\|_p \gtrsim y$ are used to indicate that there is some
nonvanishing polynomial $C(n)$ such that $\|\bm{\mathsf{x}}\|_p \leq C(n) y$ and
$C(n) \|\bm{\mathsf{x}}\|_p \geq y$, respectively.
For $\bm{\mathsf{A}} \in \C^{m \times n}$, the notations
$\|\bm{\mathsf{A}}\|_p \lesssim y$ and $\|\bm{\mathsf{A}}\|_p \gtrsim y$ are
used to indicate that there is some nonvanishing polynomial $C(m,n)$ such that
$\|\bm{\mathsf{A}}\|_p \leq C(m,n) y$ and $C(m,n) \|\bm{\mathsf{A}}\|_p \geq y$,
respectively.

\subsection{Chebyshev interpolation}

\label{cheb}

We denote by $T_n$ the Chebyshev polynomial of degree $n$.
For $k > 1$, the extrema of $T_{k-1}$ on the interval $[-1,1]$ are termed the
Chebyshev nodes, which form the $k$-point Chebyshev grid, and are denoted by
\begin{align}
-1 = t_{1,k} < t_{2,k} < \cdots < t_{k,k} = 1
,
\end{align}
where
\begin{align}
t_{j,k} = \cos\left( \frac{k - j}{k - 1} \pi \right),
\quad j = 1, \ldots, k
.
\end{align}
The $k \times k$ tensor product Chebyshev grid on $[-1,1]^2$ is given by
\begin{align}
\label{cheb:tpnodes}
\left\{
\bm{x}_{i,k}^{(j)}
\right\}_{i,j=1}^{k},
\quad
\text{where}
\quad
\bm{x}_{i,k}^{(j)}
=
( t_{i,k}, t_{j,k} ),
\end{align}
so that $i$ indexes nodes in the $x$-direction and $j$ in the $y$-direction.
We will find it useful to represent the Chebyshev grid as a single vector of
length $k^2$, for which we use the notation
\begin{align}
\bm{x}_{i + k (j-1),k^2}
=
\bm{x}_{i,k}^{(j)}
,
\quad
i, j = 1, \ldots, k
.
\end{align}
For $f \in C^1([-1,1]^2)$, we denote by $\mathsf{P}_k[f]$ the unique Chebyshev
expansion of the form
\begin{align}
\sum_{m,n=0}^{k-1} a_{mn} T_m(x) T_n(y), \quad a_{mn} \in \C
,
\end{align}
which agrees with $f$ at the nodes (\ref{cheb:tpnodes}).
It is well-known that $\mathsf{P}_k[f]$ converges to $f$ uniformly as
$k \to \infty$.
For bivariate polynomials, we use the term ``degree'' to mean the maximum of the
degrees in $x$ and the degrees in $y$, whichever is larger, of all terms in the
polynomial, rather than the total degree of the polynomial.
Evidently, if $f$ is a bivariate polynomial of degree less than $k$, then
\begin{align}
\mathsf{P}_k[f](x,y) = f(x,y)
,
\end{align}
and
\begin{align}
\mathsf{P}_k\biggl[ \pp{f}{x} \biggr](x,y)
=
\pp{f}{x}(x,y)
=
\pp{}{x}\mathsf{P}_k[f](x,y)
.
\end{align}

We denote by $\bm{\mathsf{D}}_k$ the $k \times k$ Chebyshev differentiation
matrix which takes the vector
\begin{align}
\begin{pmatrix}
f(t_{1,k}) & f(t_{2,k}) & \cdots & f(t_{k,k})
\end{pmatrix}
^T
\end{align}
of the values of a univariate polynomial $f(x)$ of degree less than $k$ at the
Chebyshev nodes to the vector
\begin{align}
\begin{pmatrix}
f'(t_{1,k}) & f'(t_{2,k}) & \cdots & f'(t_{k,k})
\end{pmatrix}
^T
\end{align}
of the values of its derivative at the same nodes.
Since the nodes (\ref{cheb:tpnodes}) are indexed in the usual lexicographical
order on $\R^2$, the matrix
\begin{align}
(\bm{\mathsf{D}}_{x})_k = \bm{\mathsf{I}}_k \otimes \bm{\mathsf{D}}_k
\end{align}
computes the values of the partial derivative of a bivariate polynomial
$f(x,y)$ of degree less than $k$ with respect to $x$ on the $k \times k$ tensor
product Chebyshev grid, where the symbol $\otimes$ denotes the Kronecker product
of matrices and $\bm{\mathsf{I}}_k$ is the $k \times k$ identity matrix.
We denote by $\bm{\mathsf{P}}_{(k,\ell)}$ the $k^2 \times \ell^2$ interpolation
matrix which takes the vector
\begin{align}
\begin{pmatrix}
f\bigl(\bm{x}_{1,\ell^2}\bigr) &
f\bigl(\bm{x}_{2,\ell^2}\bigr) &
\cdots &
f\bigl(\bm{x}_{\ell^2,\ell^2}\bigr)
\end{pmatrix}
^T
\end{align}
of the values of a bivariate polynomial $f(x,y)$ of degree less than $\ell$ on
the $\ell \times \ell$ tensor product Chebyshev grid to the vector
\begin{align}
\begin{pmatrix}
f\bigl(\bm{x}_{1,k^2}\bigr) &
f\bigl(\bm{x}_{2,k^2}\bigr) &
\cdots &
f\bigl(\bm{x}_{k^2,k^2}\bigr)
\end{pmatrix}
^T
\end{align}
of the values of the polynomial on the $k \times k$ tensor product Chebyshev
grid.

It is well known from classic Sturm-Liouville theory (see, for instance,
Section~4 of \cite{Bernardi1989}) that,
if $f \in C^\infty([-1,1]^2)$, then it admits a uniformly convergent Chebyshev
series
\begin{align}
f(x,y) = \sum_{m,n=0}^{\infty} b_{mn} T_m(x) T_n(y)
,
\end{align}
such that $|b_{mn}|$ decays super-algebraically as $m,n \to \infty$.
It follows from the aliasing formula for Chebyshev coefficients that
\begin{align}
\label{cheb:ftail}
\left\| \mathsf{P}_k[f] - f \right\|_{L^\infty([-1,1]^2)}
\leq
2 \sum_{\max\{m,n\} \geq k} |b_{mn}|
,
\end{align}
and
\begin{align}
\label{cheb:fptail}
\left\| \pp{}{x}\mathsf{P}_k[f] - \pp{f}{x} \right\|_{L^\infty([-1,1]^2)}
\leq
2 \sum_{\max\{m,n\} \geq k} m^2 |b_{mn}|
.
\end{align}

\begin{lemma}
\label{cheb:lebesgue}
Suppose that $f$ is a bivariate polynomial of degree less than $k$, and that
\begin{align}
\bm{\mathsf{f}}
=
\begin{pmatrix}
f\bigl(\bm{x}_{1,k^2}\bigr) &
f\bigl(\bm{x}_{2,k^2}\bigr) &
\cdots &
f\bigl(\bm{x}_{k^2,k^2}\bigr)
\end{pmatrix}
^T
.
\end{align}
Then,
\begin{align}
\left\| f \right\|_{L^\infty([-1,1]^2)}
\leq
\Lambda_k^2 \, \|\bm{\mathsf{f}}\|_\infty
,
\end{align}
where $\Lambda_k$ is the Lebesgue constant for the $k$-point Chebyshev grid on
$[-1,1]$.
\end{lemma}

\begin{proof}
Let $(x,y) \in [-1,1]^2$ with
\begin{align}
|f(x,y)|
=
\left|
\sum_{i=1}^{k} \sum_{j=1}^{k} f( t_{i,k}, t_{j,k} ) \ell_i(x) \ell_j(y)
\right|
=
\left\| f \right\|_{L^\infty([-1,1]^2)}
,
\end{align}
where $\ell_1, \ldots, \ell_k$ are the Lagrange polynomials associated with the
$k$-point Chebyshev grid.
The conclusion of the lemma follows immediately from this and the definition
\begin{align}
\Lambda_k
=
\max_{-1 \leq x \leq 1} \sum_{i=1}^{k} |\ell_i(x)|
.
\end{align}
\end{proof}

\subsection{Tempered distributions}

In this subsection, we collect some standard results from the theory of
distributions that we will use in the analysis of the Levin PDE.

\begin{definition}
For $\cT \in \sS'(\R^n)$ and $\varphi \in \sS(\R^n)$, the convolution
$\cT \ast \varphi$ is defined as
$
\cT \ast \varphi(\bm{x})
=
\langle
\cT, \tau_{\bm{x}}\tilde{\varphi}
\rangle
$
for all $\bm{x} \in \R^n$, where
$\tau_{\bm{x}} \varphi(\bm{y}) = \varphi(\bm{y} - \bm{x})$ denotes translation
by $\bm{x}$ and $\tilde{\varphi}(\bm{x}) = \varphi(-\bm{x})$ denotes reflection
through the origin.
\end{definition}

The following lemma is reproduced from Proposition~9.10 of \cite{Folland1999}.

\begin{lemma}
\label{dist:conv}
If $\cT \in \sS'(\R^n)$ and $\varphi \in \sS(\R^n)$, then $\cT \ast \varphi$
is a slowly increasing $C^\infty$ function which satisfies
\begin{align}
\label{dist:convp}
D^{\bm{\alpha}}(\cT \ast \varphi)
=
(D^{\bm{\alpha}}\cT) \ast \varphi
=
\cT \ast (D^{\bm{\alpha}}\varphi)
\end{align}
for all multi-indices $\bm{\alpha}$.
Furthermore,
\begin{align}
\supp (\cT \ast \varphi) \subset \supp \cT + \supp \varphi
.
\end{align}
\end{lemma}

The following provides the standard definition for the tensor product of two
tempered distributions.

\begin{definition}
For $\cS \in \sS'(\R^m)$ and $\cT \in \sS'(\R^n)$,
the tensor product $\cS \otimes \cT$ is a tempered distribution on
$\R^m \times \R^n$ defined via the formula
\begin{align}
\langle
\cS(\bm{x}) \otimes \cT(\bm{y}), \varphi(\bm{x},\bm{y})
\rangle
=
\langle
\cS(\bm{x}),
\langle
\cT(\bm{y}), \varphi(\bm{x},\bm{y})
\rangle
\rangle
\end{align}
for all $\varphi(\bm{x},\bm{y}) \in \sS(\R^m \times \R^n)$.
\end{definition}

The notation we have used has the obvious meaning that,
for each fixed $\bm{x} \in \R^m$,
$\langle \cT(\bm{y}), \varphi(\bm{x},\bm{y}) \rangle$ is the action of the
distribution $\cT$ on the function $\varphi(\bm{x},\cdot) \in \sS(\R^n)$.
A proof of the fact that $\langle \cT(\bm{y}), \varphi(\bm{x},\bm{y}) \rangle$
defines a function in $\sS(\R^m)$ can be found in Chapter~7 of
\cite{Richards1995}.

A proof of the following lemma, which states that the tensor product commutes,
can be found in the same place.

\begin{lemma}
\label{dist:tpcommute}
The tensor product is commutative, that is,
$\cS(\bm{x}) \otimes \cT(\bm{y}) = \cT(\bm{y}) \otimes \cS(\bm{x})$.
\end{lemma}

In Theorem~\ref{dist:structthm}, we prove a structure theorem for tempered
distributions that are compactly supported in a bounded convex set with nonempty
interior.

A key ingredient in our proof is the continuity of the following extension
operator~$E$, which is constructed in Chapter~VI of \cite{Stein2016}.

\begin{lemma}
\label{dist:extn}
Let $U \subset \R^n$ be an open, bounded, convex set. Then there exists a linear
extension operator $E \colon C_b(U) \to C_b(\R^n)$ such that $E$ maps
$C_b^m(U)$ continuously into $C_b^m(\R^n)$ for all nonnegative integers $m$, and
$E[f]$, $f \in C_b(U)$, is supported in a fixed compact set that depends only on
the domain $U$.
\end{lemma}

Note that, if $\cT \in \sE'(\R^n)$ is of order at most $m$, then formula
(\ref{prelim:order}), together with the fact that $C_c^\infty(\R^n)$ is dense in
$C^m(\R^n)$, means that $\cT$ can be continuously extended to a bounded 
linear functional on $C^m(\R^n)$.
The following lemma can be proven by constructing a sequence of bells
$(\rho_n)_{n=1}^{\infty}$ in $C_c^\infty(\R^n)$ such that
$\{ \bm{x} : \rho_n(\bm{x}) = 1 \}^\circ \supset \supp \cT$
and $\bigcap_{n=1}^{\infty} \supp \rho_n = \supp \cT$,
and showing that
$\langle \cT, \varphi \rangle = \langle \cT, \rho_n \varphi \rangle \to 0$.

\begin{lemma}
\label{dist:vanish}
Let $\cT \in \sE'(\R^n)$ be of order at most $m$, and $\varphi \in C^m(\R^n)$
with $D^{\bm{\alpha}}\varphi = 0$ on $\supp \cT$ for all $|\bm{\alpha}| \leq m$.
Then $\langle \cT, \varphi \rangle = 0$.
\end{lemma}

The following theorem provides a structure theorem for tempered distributions
that are compactly supported in a bounded convex set with nonempty interior.
The proof closely follows Appendix~2 of \cite{Richards1995}.

\begin{theorem}
\label{dist:structthm}
Let $U \subset \R^n$ be an open, bounded, convex set, and suppose that
$\cT \in \sE'(\R^n)$ is of order at most $m$ with
$\supp \cT \subset \overbar{U}$.
Then there exist complex Radon measures $\mu_{\bm{\alpha}} \in M(\overbar{U})$,
corresponding to the set of multi-indices $\bm{\alpha}$ with
$|\bm{\alpha}| \leq m$, such that
\begin{align}
\langle \cT, \varphi \rangle
=
\sum_{|\bm{\alpha}| \leq m}
\int_{\overbar{U}}
D^{\bm{\alpha}}\varphi(\bm{x})
\, \mathrm{d}\mu_{\bm{\alpha}}(\bm{x}),
\end{align}
for all $\varphi \in C^m(\R^n)$.
\end{theorem}

\begin{proof}
Using the extension operator $E$ in Lemma~\ref{dist:extn},
we know that there exists
a compact set $K \supset \overbar{U}$, depending only on $U$, and 
a constant $A$, depending only on $U$ and $m$,
such that, for each $\varphi \in C^m(\R^n)$,
the function
$\tilde{\varphi} = E[\varphi\vert_{U}]$ is supported in $K$, satisfies
$D^{\bm{\alpha}}\tilde{\varphi} = D^{\bm{\alpha}}\varphi$ on $\overbar{U}$
for all $|\bm{\alpha}| \leq m$, and satisfies the inequality
\begin{align}
\label{dist:phibound}
\sup_{|\bm{\alpha}| \leq m} \sup_{\bm{x} \in K}
\left| D^{\bm{\alpha}}\tilde{\varphi}(\bm{x}) \right|
\leq
A
\sup_{|\bm{\alpha}| \leq m} \sup_{\bm{x} \in \overbar{U}}
\left| D^{\bm{\alpha}}\varphi(\bm{x}) \right|.
\end{align}
It then follows that
\begin{align}
\label{dist:tbound}
\begin{split}
\left|
\left\langle \cT, \varphi \right\rangle
\right|
=
\left|
\left\langle \cT, \tilde{\varphi} \right\rangle
\right|
&\leq
M_{K}
\sup_{|\bm{\alpha}| \leq m} \sup_{\bm{x} \in K}
\left| D^{\bm{\alpha}}\tilde{\varphi}(\bm{x}) \right|\\
&\leq
A M_{K}
\sup_{|\bm{\alpha}| \leq m} \sup_{\bm{x} \in \overbar{U}}
\left| D^{\bm{\alpha}}\varphi(\bm{x}) \right|,
\end{split}
\end{align}
where the equality comes from Lemma~\ref{dist:vanish}, the first
inequality comes from the fact that $\cT$ extends to a bounded linear
functional on $C^m(\R^n)$, and the second inequality follows from
inequality~(\ref{dist:phibound}).

Now let $N = |\{\bm{\alpha} \colon |\bm{\alpha}| \leq m\}|$, and $S$ be the
union of
$N$ disjoint copies of $\overbar{U}$ endowed with the disjoint union topology.
The space $C(S)$ is the Banach space of all continuous complex-valued functions
$\psi = (\psi_1, \ldots, \psi_N)$ on $S$, where $\psi_k \in C(\overbar{U})$
for all $1 \leq k \leq N$.
Let $Y \subset C(S)$ be the subspace consisting of all $m$-jets
$(D^{\bm{\alpha}}\varphi)_{|\bm{\alpha}| \leq m}$ for $\varphi \in C^m(\R^n)$.

From (\ref{dist:tbound}), it is clear that $\cT$ induces a bounded linear
operator on $Y$. By the Hahn-Banach theorem, $\cT$ has an extension to a
bounded linear operator on $C(S)$.
Since $S$ is compact, it then follows from the Riesz representation theorem that
there exists a complex Radon measure $\mu \in M(S)$ such that,
for all $\varphi \in C^m(\R^n)$,
\begin{align}
\langle \cT, \varphi \rangle
=
\sum_{|\bm{\alpha}| \leq m}
\int_{\overbar{U}}
D^{\bm{\alpha}}\varphi(\bm{x})
\, \mathrm{d}\mu_{\bm{\alpha}}(\bm{x})
,
\end{align}
where $\mu_{\bm{\alpha}} \in M(\overbar{U})$ are the restrictions of $\mu$ to
the disjoint copies of $\overbar{U}$.
\end{proof}

\subsection{Bandlimited approximation of bivariate functions}

In this subsection, we provide an explicit construction of certain bandlimited
approximations of a function $f$ defined on a closed subset of $[-1,1]^2$, under
the mild regularity condition that $f$ extends smoothly to a neighborhood of its
domain.

\begin{theorem}
\label{approx:thm}
Let $\Omega \subset [-1,1]^2$ be closed.
Suppose that $f \colon \Omega \to \C$ admits an infinitely differentiable
extension to an open neighborhood of $\Omega$.
Then, for each positive integer $m$ and each real number $c>1$, there exist a
constant $C(m,f)$, depending on $m$ and $f$, but not $c$, and a function
$f_b \in \sS(\R^2)$ such that
\begin{enumerate}
\item
$\hat{f}_b$ is supported in $[-c-4,c+4]^2$,

\item
$
\begin{aligned}
\left\| f_b - f \right\|_{L^\infty(\Omega)}
\leq
\frac{C(m,f)}{c^{2m}}
\end{aligned}
,
$

\item
$
\begin{aligned}
\left\| f_b \right\|_{L^\infty(\R^2)}
\leq
2 \left\| f \right\|_{L^\infty(\Omega)}
+
\frac{C(m,f)}{c^{2m}}
\end{aligned}
,
$

\item
$
\begin{aligned}
\left\| D^{\bm{\alpha}}\hat{f}_b \right\|_{L^\infty(\R^2)}
\leq
10 \left\| f \right\|_{L^\infty(\Omega)}
\quad
\text{for all $|\bm{\alpha}| \leq 2$}.
\end{aligned}
$
\end{enumerate}
\end{theorem}

\begin{proof}
Let $\tilde{f} \colon U \to \C$ be an infinitely differentiable extension of $f$
to an open neighborhood $U$ of $\Omega$.
Let $M = \left\| f \right\|_{L^\infty(\Omega)}$ and choose an open neighborhood
$V$ of $\Omega$ such that $V \subset U \cap [-2,2]^2$ and $\tilde{f}$ is bounded
in magnitude by $2M$ on $V$.

By the $C^\infty$ Urysohn's lemma (see Chapter~8 of \cite{Folland1999}),
there exists an infinitely differentiable function $g_1 \colon \R^2 \to \R$
such that $0 \leq g_1(x,y) \leq 1$ for all $(x,y) \in \R^2$, $g_1(x,y) = 1$ on
$\Omega$ and $\supp g_1 \subset V$.
Then $f_1(x,y) = f(x,y) g_1(x,y)$ is an element of $\sS(\R^2)$, and so is its
Fourier transform $\hat{f}_1$. Consequently, there exists a constant $k_1$
(depending on $m$ and $f$) such that
\begin{align}
\label{approx:mdecay}
\left| \hat{f}_1(\xi_1,\xi_2) \right|
\leq
\frac{k_1}{\left( 1 + \xi_1^2 + \xi_2^2 \right)^{m+1}}
,
\end{align}
for all $(\xi_1,\xi_2) \in \R^2$.
Since $g_1$ is bounded in magnitude by $1$, we have
\begin{align}
\label{approx:f1bound}
\left| f_1(x,y) \right| \leq 2M
,
\end{align}
for all $(x,y) \in \R^2$.
Since $\supp f_1 \subset V \subset [-2,2]^2$, it follows that,
for all multi-indices $\bm{\alpha}$,
\begin{align}
\label{approx:f1fpbound}
\begin{split}
\left| D^{\bm{\alpha}}\hat{f}_1(\xi_1,\xi_2) \right|
&=
\left|
\frac{(-i)^{|\bm{\alpha}|}}{4 \pi^2}
\int_{V} \bm{x}^{\bm{\alpha}} f_1(x,y) \exp(-i (\xi_1 x + \xi_2 y))
\, \mathrm{d}x \, \mathrm{d}y
\right|\\
&\leq
2^{|\bm{\alpha}|} \frac{8}{\pi^2} M
\end{split}
\end{align}
holds for all $(\xi_1,\xi_2) \in \R^2$, where $\bm{x} = (x,y) \in \R^2$.

Using the infinitely differentiable ramp function
\begin{align}
H(x)
=
\begin{cases}
0, & x \leq -1,\\
\frac{1}{2} \left( 1 + \mbox{erf}\left( \frac{x}{\sqrt{1-x^2}} \right) \right),
& x \in (-1,1),\\
1, & x \geq 1,
\end{cases}
\end{align}
suggested in \cite{Boyd2003}, one can construct an infinitely differentiable
window function $g_2 \colon \R \to \R$ such that $0 \leq g_2(x) \leq 1$,
$|g_2'(x)| \leq 1$ and $|g_2''(x)| \leq 1$ for all $x \in \R$,
$g_2(x) = 1$ for all $|x| \leq c$ and $g_2(x) = 0$ for all $|x| \geq c+4$.
We then define $f_b$ by its Fourier transform
\begin{align}
\label{approx:deffbf}
\hat{f}_b(\xi_1,\xi_2)
=
\hat{f}_1(\xi_1,\xi_2) g_2(\| (\xi_1,\xi_2) \|)
,
\end{align}
so that the first condition listed above is satisfied.
From (\ref{approx:mdecay}) and the definition of $g_2$, it is clear that, for
all $(x,y) \in \R^2$,
\begin{align}
\begin{split}
&\left| f_1(x,y) - f_b(x,y) \right|\\
&=
\left|
\int_{\R^2}
\hat{f_1}(\xi_1,\xi_2)
\left( 1 - g_2(\| (\xi_1,\xi_2) \|) \right)
\exp(i (\xi_1 x + \xi_2 y))
\, \mathrm{d}\xi_1 \, \mathrm{d}\xi_2
\right|\\
&\leq
\int_{0}^{2\pi} \int_{c}^{\infty}
\frac{k_1 r}{\left( 1+r^2 \right)^{m+1}}
\, \mathrm{d}r \, \mathrm{d}\theta\\
&=
\frac{k_1 \pi}{m \left( 1+c^2 \right)^m}
.
\end{split}
\end{align}
This inequality leads directly to the second condition listed above, and we
obtain the third condition by combining it with (\ref{approx:f1bound}).
Using (\ref{approx:f1fpbound}) and the properties of $g_2$,
it is straightforward to verify that the remaining conditions listed above are
satisfied by repeated differentiation of (\ref{approx:deffbf}).
\end{proof}

\begin{definition}
\label{approx:def}
Let $\Omega \subset [-1,1]^2$ be closed.
Suppose that $f \colon \Omega \to \C$ admits an infinitely differentiable
extension to an open neighborhood of $\Omega$.
For $0 < \varepsilon < 1$, we say that a function $f_b \in \sS(\R^2)$ is an
$\varepsilon$-bandlimited approximation of $f$ if $f_b$ is a bandlimited
function with the following properties:
\begin{enumerate}
\item
$
\begin{aligned}
\left\| f_b - f \right\|_{L^\infty(\Omega)}
\leq
\varepsilon \left\| f \right\|_{L^\infty(\Omega)}
\end{aligned}
,
$

\item
$
\begin{aligned}
\left\| f_b \right\|_{L^\infty(\R^2)}
\leq
3 \left\| f \right\|_{L^\infty(\Omega)}
\end{aligned}
,
$

\item
$
\begin{aligned}
\left\| D^{\bm{\alpha}}\hat{f}_b \right\|_{L^\infty(\R^2)}
\leq
10 \left\| f \right\|_{L^\infty(\Omega)}
\quad
\text{for all $|\bm{\alpha}| \leq 2$}.
\end{aligned}
$
\end{enumerate}
Furthermore, we denote by $c_f(\varepsilon)$ the infimum of the set of positive
real numbers $c$ such that there exists an $\varepsilon$-bandlimited
approximation of $f$ with bandlimit $c$.
Note that Theorem~\ref{approx:thm} shows that this set is nonempty, hence
$c_f(\varepsilon) < \infty$.
\end{definition}

Throughout this paper, we will assume that $c_f(\varepsilon)$ is positive
because it is zero only in the trivial case where $f$ is the zero function.

\subsection{Legendre expansions of bandlimited functions}

Here we characterize the super-exponential decay of the coefficients of the
Legendre expansion of a bandlimited function in terms of its bandlimit.
Throughout this subsection, $P_n$ denotes the Legendre polynomial of degree $n$,
and $j_\nu$ denotes the spherical Bessel function of order $\nu$.

\begin{lemma}
\label{lege:lemma}
Let $\varphi \in \sS'(\R^2)$. If the Fourier transform of $\varphi$ is a
tempered distribution of order at most $p$ supported in $[-c,c]^2$, then
$\varphi$ coincides with an entire function of the form
\begin{align}
\varphi(x,y)
=
\sum_{|\bm{\alpha}| \leq p}
\bm{x}^{\bm{\alpha}}
\int_{[-c,c]^2}
\exp(i (\xi_1 x + \xi_2 y))
\, \mathrm{d}\mu_{\bm{\alpha}}(\xi_1,\xi_2)
,
\end{align}
where $\mu_{\bm{\alpha}} \in M([-c,c]^2)$ for all multi-indices $\bm{\alpha}$
with $|\bm{\alpha}| \leq p$.
\end{lemma}

\begin{proof}
In view of the formula
\begin{align}
\varphi(x,y)
=
\langle \hat{\varphi}(\xi_1,\xi_2), \exp(i (\xi_1 x + \xi_2 y)) \rangle
,
\end{align}
the conclusion of the lemma is a direct consequence of
Theorem~\ref{dist:structthm}.
\end{proof}

\begin{theorem}
Let $\varphi \in \sS'(\R^2)$ such that the Fourier transform of $\varphi$ is a
tempered distribution of order at most $p$ supported in $[-c,c]^2$ for
$c \geq 1$. Then $\varphi$ admits a uniformly convergent Legendre expansion
\begin{align}
\label{lege:expn}
\varphi(x,y) = \sum_{m,n=0}^{\infty} a_{mn} P_m(x) P_n(y)
,
\end{align}
and there exists a constant $C(\varphi,p)$, depending only on $\varphi$ and $p$,
but not $m$ and $n$, such that
\begin{align}
|a_{mn}|
\leq
C(\varphi,p) \frac{\left( c / 2 \right)^{m+n+p}}{\Gamma(m-p+1) \Gamma(n-p+1)}
\end{align}
for all $\min\{m,n\} \geq p$, and
\begin{align}
|a_{mn}|
\leq
C(\varphi,p) \frac{\left( c / 2 \right)^{p'+p}}{\Gamma(p'-p+1)}
\end{align}
for all $p' = \max\{m,n\} \geq p$.
\end{theorem}

\begin{proof}
By Lemma~\ref{lege:lemma},
\begin{align}
\label{lege:amnbound}
\begin{split}
&
\left|
\int_{-1}^{1} \int_{-1}^{1}
\varphi(x,y) P_m(x) P_n(y)
\, \mathrm{d}x \, \mathrm{d}y
\right|\\
&\leq
\max_{|\bm{\alpha}| \leq p} | \mu_{\bm{\alpha}} |([-c,c]^2)
\sum_{0 \leq k + \ell \leq p}
\max_{\xi \in [-c,c]} |I_{km}(\xi)|
\max_{\xi \in [-c,c]} |I_{\ell n}(\xi)|
,
\end{split}
\end{align}
where
\begin{align}
I_{k \ell}(\xi)
=
\int_{-1}^{1}
x^k P_\ell(x) \exp(i \xi x)
\, \mathrm{d}x
.
\end{align}
For $k \leq \ell$, we have the expansion
\begin{align}
x^k P_\ell(x)
=
\sum_{\nu = \ell-k}^{\ell+k} b_{k \ell \nu} P_\nu(x)
,
\end{align}
where
\begin{align}
b_{k \ell \nu}
=
\frac{2\nu + 1}{2} \int_{-1}^{1} x^k P_\ell(x) P_\nu(x) \, \mathrm{d}x
.
\end{align}
Using the three-term recurrence relation for Legendre polynomials,
it can be shown that $|b_{k \ell \nu}| \leq 1$.
It then follows from Formula 10.54.2 and 10.14.4 in \cite{NIST2010} that
\begin{align}
\label{lege:iklbound}
\begin{split}
|I_{k \ell}(\xi)|
&\leq
\sum_{\nu = \ell-k}^{\ell+k}
\left|
\int_{-1}^{1}
P_\nu(x) \exp(i \xi x)
\, \mathrm{d}x
\right|\\
&=
2
\sum_{\nu = \ell-k}^{\ell+k}
\left|
j_\nu(|\xi|)
\right|\\
&\leq
2
\sum_{\nu = \ell-k}^{\ell+k}
\frac{\left| \xi / 2 \right|^{\nu}}{\Gamma(\nu + 1)}\\
&\leq
(4k + 2)
\frac{\left( c / 2 \right)^{\ell + k}}{\Gamma(\ell-k+1)}
\end{split}
\end{align}
holds for all $|\xi| \leq c$.
Finally, inserting the bound
\begin{align}
\label{lege:iklbound2}
|I_{k \ell}(\xi)| \leq 2
\end{align}
and (\ref{lege:iklbound}) into (\ref{lege:amnbound})
yields the conclusion of the theorem.
\end{proof}

\begin{remark}
Using the asymptotic expansion
\begin{align}
j_\nu(\xi)
\sim
\frac{\sqrt{e / 2}}{2\nu + 1}
\left(\frac{e \xi}{2\nu + 1}\right)^{\nu}
\quad
\text{as $\nu \to \infty$},
\end{align}
which is an application of Formula 10.19.1 in \cite{NIST2010}, in place of the
inequality used in (\ref{lege:iklbound}), we obtain the approximation
\begin{align}
\label{lege:iklapprox}
\begin{split}
|I_{k \ell}(\xi)|
&\leq
2
\sum_{\nu = \ell-k}^{\ell+k}
\left|
j_\nu(|\xi|)
\right|\\
&\sim
2
\sum_{\nu = \ell-k}^{\ell+k}
\frac{\sqrt{e / 2}}{2\nu + 1}
\left(\frac{e |\xi|}{2\nu + 1}\right)^{\nu}\\
&\leq
\sqrt{2e}
\sum_{\nu = \ell-k}^{\ell+k}
\frac{1}{2\nu + 1}
\left(\frac{e c}{2\nu + 1}\right)^{\nu}
,
\end{split}
\end{align}
as $\ell \to \infty$,
which holds for all $|\xi| \leq c$. Applying the bounds (\ref{lege:iklbound2})
and (\ref{lege:iklapprox}) to (\ref{lege:amnbound}), together with the condition
\begin{align}
\label{lege:condition}
\max\{m,n\} > \frac{e}{2}c + p - \frac{1}{2}
,
\end{align}
we see that
\begin{align}
\frac{ec}{2\nu + 1} < 1
\end{align}
holds for all $\nu$ involved, so we expect the coefficients in the Legendre
expansion (\ref{lege:expn}) of $\varphi$ to decay super-exponentially as soon as
(\ref{lege:condition}) is satisfied.
\end{remark}

\subsection{Truncated singular value decompositions}

If $\bm{\mathsf{A}}$ is a complex-valued $m \times n$ matrix with $m \geq n$,
then any decomposition of the form
\begin{align}
\label{tsvd:svd}
\bm{\mathsf{A}}
=
\begin{pmatrix}
\bm{\mathsf{u}}_1 & \bm{\mathsf{u}}_2 & \cdots & \bm{\mathsf{u}}_m
\end{pmatrix}
\begin{pmatrix}
\sigma_1 &          &        &          \\
         & \sigma_2 &        &          \\
         &          & \ddots &          \\
         &          &        & \sigma_n \\
0 & 0 & \cdots & 0 \\
\vdots & \vdots & \ddots & \vdots \\
0 & 0 & \cdots & 0
\end{pmatrix}
\begin{pmatrix}
\bm{\mathsf{v}}_1 & \bm{\mathsf{v}}_2 & \cdots & \bm{\mathsf{v}}_n
\end{pmatrix}
^{*}
,
\end{align}
where $\sigma_1, \ldots, \sigma_n \geq 0$,
$\{\bm{\mathsf{u}}_1, \ldots, \bm{\mathsf{u}}_m\}$ is an orthonormal basis of
$\C^m$, and
$\{\bm{\mathsf{v}}_1, \ldots, \bm{\mathsf{v}}_n\}$ is an orthonormal basis of
$\C^n$,
is known as a singular value decomposition of $\bm{\mathsf{A}}$.
The scalars $\sigma_1, \ldots, \sigma_n$ are uniquely determined up to ordering,
and they are known as the singular values of $\bm{\mathsf{A}}$.
It is conventional to arrange them in descending order, and we will assume that
this is the case with all singular value decompositions that we consider.

A $k$-truncated singular value decomposition of $\bm{\mathsf{A}}$ is an
approximation of the form
\begin{align}
\label{tsvd:tsvd}
\bm{\mathsf{A}}
\approx
\bm{\mathsf{A}}_k
=
\begin{pmatrix}
\bm{\mathsf{u}}_1 & \bm{\mathsf{u}}_2 & \cdots & \bm{\mathsf{u}}_k
\end{pmatrix}
\begin{pmatrix}
\sigma_1 &          &        &          \\
         & \sigma_2 &        &          \\
         &          & \ddots &          \\
         &          &        & \sigma_k
\end{pmatrix}
\begin{pmatrix}
\bm{\mathsf{v}}_1 & \bm{\mathsf{v}}_2 & \cdots & \bm{\mathsf{v}}_k
\end{pmatrix}
^{*}
,
\end{align}
where (\ref{tsvd:svd}) is a singular value decomposition of $\bm{\mathsf{A}}$
with $\sigma_1 \geq \sigma_2 \geq \cdots \geq \sigma_n$ and $1 \leq k \leq n$.
Typically, a desired precision $\varepsilon > 0$ is specified and $k$ is taken
to be the least integer between $1$ and $n-1$ such that
$\sigma_{k+1} < \varepsilon$, if such an integer exists, or $k=n$ otherwise.
In this case, we say that (\ref{tsvd:tsvd}) is a singular value decomposition
which has been truncated at precision $\varepsilon$.
We call the vector $\bm{\mathsf{x}} = \bm{\mathsf{A}}_k^\dagger \bm{\mathsf{y}}$
given by
\begin{align}
\bm{\mathsf{x}}
=
\begin{pmatrix}
\bm{\mathsf{v}}_1 & \bm{\mathsf{v}}_2 & \cdots & \bm{\mathsf{v}}_k
\end{pmatrix}
\begin{pmatrix}
\frac{1}{\sigma_1} &                    &        &                    \\
                   & \frac{1}{\sigma_2} &        &                    \\
                   &                    & \ddots &                    \\
                   &                    &        & \frac{1}{\sigma_k}
\end{pmatrix}
\begin{pmatrix}
\bm{\mathsf{u}}_1 & \bm{\mathsf{u}}_2 & \cdots & \bm{\mathsf{u}}_k
\end{pmatrix}
^{*}
\bm{\mathsf{y}}
\end{align}
the solution of the linear system
$\bm{\mathsf{A}}\bm{\mathsf{x}}=\bm{\mathsf{y}}$ obtained from the truncated
singular value decomposition (\ref{tsvd:tsvd}).

The following theorem implies that, when a linear system admits an approximate
solution with a modest norm, and it is solved numerically using a truncated
singular value decomposition, the computed solution will have both a small
residual and a modest norm.
We refer the reader to \cite{Zhao2025} for a proof of the theorem.

\begin{theorem}
\label{tsvd:thm}
Suppose that $\bm{\mathsf{A}} \in \C^{m \times n}$, where $m \geq n$, and let
$\sigma_1 \geq \sigma_2 \geq \cdots \geq \sigma_n$ denote the singular values of
$\bm{\mathsf{A}}$. Let $\bm{\mathsf{b}} \in \C^m$, and suppose that
$\bm{\mathsf{x}} \in \C^n$ satisfies
\begin{align}
\bm{\mathsf{A}} \bm{\mathsf{x}} = \bm{\mathsf{b}}
.
\end{align}
Let $\varepsilon > 0$, and suppose further that $\bm{\mathsf{E}} \in \C^{m
\times n}$ and $\bm{\mathsf{e}} \in \C^m$, with $\|\bm{\mathsf{E}}\|_2 <
\varepsilon/2$, and let $\hat{\sigma}_1 \ge \hat{\sigma}_2 \ge \cdots \ge
\hat{\sigma}_n$
denote the singular values of $\bm{\mathsf{A}} + \bm{\mathsf{E}}$,
defining $\hat{\sigma}_{n+1} = 0$. Suppose that
\begin{align}
\hat{\bm{\mathsf{x}}}_k
=
(\bm{\mathsf{A}} + \bm{\mathsf{E}})_k^\dagger
(\bm{\mathsf{b}} + \bm{\mathsf{e}}),
\end{align}
where $(\bm{\mathsf{A}} + \bm{\mathsf{E}})_k^\dagger$ is the pseudoinverse of
the $k$-truncated singular value decomposition of
$\bm{\mathsf{A}} + \bm{\mathsf{E}}$, so that
\begin{align}
\hat{\sigma}_k \geq \varepsilon \ge \hat{\sigma}_{k+1}.
\end{align}
%
Then
\begin{align}
\| \hat{\bm{\mathsf{x}}}_k \|_2
\leq
\frac{1}{\hat{\sigma}_k}
(2 \varepsilon \|\bm{\mathsf{x}}\|_2 + \|\bm{\mathsf{e}}\|_2 )
+
\|\bm{\mathsf{x}}\|_2
,
\end{align}
and
\begin{align}
\| \bm{\mathsf{A}} \hat{\bm{\mathsf{x}}}_k - \bm{\mathsf{b}} \|_2
\leq
5 \varepsilon \|\bm{\mathsf{x}}\|_2
+
\frac{3}{2} \|\bm{\mathsf{e}}\|_2
.
\end{align}
\end{theorem}

\section{Analysis of the Levin PDE}

\label{section:leveq}

In this section, we prove that, under mild conditions on $g_x$,
the Levin PDE
\begin{align}
\nabla \cdot \bm{p} + i \, \nabla g \cdot \bm{p} = f
\quad \text{in $(-1,1)^2$}
\end{align}
admits a slowly-varying, approximate solution of the form
\begin{align}
\label{leveq:soln}
\bm{p}(x,y)
=
\begin{pmatrix}
p_b(x,y) \\ 0
\end{pmatrix}
,
\end{align}
where $p_b \colon [-1,1]^2 \to \C$ is an infinitely differentiable function.

The main results of this section are Theorem~\ref{leveq:thmhigh} and
Theorem~\ref{leveq:thmlow}.
Theorem~\ref{leveq:thmhigh} shows that, when $g_x$ is nonzero over the domain
and the ratio of its maximum absolute value to its minimum absolute value is
small, the Levin PDE admits a slowly-varying, approximate solution of the form
(\ref{leveq:soln}), whose complexity is independent of the magnitude of $g_x$.
Theorem~\ref{leveq:thmlow} shows that, when $g_x$ is of small magnitude over the
domain, the Levin PDE also admits a slowly-varying, approximate solution of the
form (\ref{leveq:soln}), whose complexity decreases as the magnitude of $g_x$
decreases.

We begin our analysis with the following technical lemma, which applies when the
domain is an arbitrary closed subset of $[-1,1]^2$ and $g_x$ is equal to a
nonzero constant $W$.

\begin{lemma}
\label{leveq:lemma1}
Let $\Omega \subset [-1,1]^2$ be closed.
Suppose that $f \colon \Omega \to \C$ admits an infinitely differentiable
extension to an open neighborhood of $\Omega$ and $W \neq 0$.
Suppose further that $0 < \varepsilon < 1$, and let $W_0 = c_f(\varepsilon)$.
Then there exist a right-continuous function $C: (0, \infty) \to \R$ and, for
each $\delta > 0$, an entire function $p_{b,\delta} \in \sS'(\R^2)$ such that
\begin{enumerate}
\item
$\hat{p}_{b,\delta}$ is a tempered distribution of order $1$ supported in
$[-W_0 - \delta, W_0 + \delta]^2$,

\item
$
\begin{aligned}
\left|
\pp{p_{b,\delta}}{x}(x,y) + i W p_{b,\delta}(x,y) - f(x,y)
\right|
\leq
\varepsilon \left\| f \right\|_{L^\infty(\Omega)}
\end{aligned}
$
for all $(x,y) \in \Omega$,

\item
$
\begin{aligned}
\left\| p_{b,\delta} \right\|_{L^\infty(\Omega)}
\leq
C(W_0 + \delta)
\min\left\{ 1, \frac{1}{|W|} \right\}
\left\| f \right\|_{L^\infty(\Omega)}
\end{aligned}
,
$

\item
$
\begin{aligned}
\left\| \pp{p_{b,\delta}}{x} \right\|_{L^\infty(\Omega)}
\leq
C(W_0 + \delta)
\min\left\{ 1, \frac{1}{|W|} \right\}
\left\| f \right\|_{L^\infty(\Omega)}
\end{aligned}
.
$
\end{enumerate}
\end{lemma}

\begin{proof}
We let $f_{b,\delta} \in \sS(\R^2)$ be an $\varepsilon$-bandlimited
approximation of $f$ with bandlimit $W_0 + \delta$, and define the function
$p_{b,\delta}$ via the formula
\begin{align}
p_{b,\delta}(x,y)
=
\int_{-\infty}^{\infty}
\left[
\princip \int_{-\infty}^{\infty}
\frac{\hat{f}_{b,\delta}(\xi_1,\xi_2)}{i (\xi_1 + W)}
\exp(i (\xi_1 x + \xi_2 y))
\, \mathrm{d}\xi_1
\right]
\mathrm{d}\xi_2
.
\end{align}
Let $\cS$ denote the tempered distribution induced by the constant function
$1$, and $\cT_W$ denote the tempered distribution defined via the formula
\begin{align}
\label{leveq:deftw}
\langle
\cT_W, \varphi
\rangle
=
\princip \int_{-\infty}^{\infty}
\frac{\varphi(\xi)}{i (\xi + W)}
\, \mathrm{d}\xi
.
\end{align}
Since $\hat{f}_{b,\delta}$ is compactly supported, $p_{b,\delta}$ is exactly
the inverse Fourier transform of the product of $\cT_W \otimes \cS$ and the
infinitely differentiable function $\hat{f}_{b,\delta}$.
Since $\cT_W$ is of order $1$ and $\cS$ is of order $0$, and
$f_{b,\delta}$ has bandlimit $W_0 + \delta$, it is clear that
$\hat{p}_{b,\delta}$ is a tempered distribution of order $1$ which is supported
in $[-W_0 - \delta, W_0 + \delta]^2$, so the first condition listed above is
satisfied. In particular, the Schwartz-Paley-Wiener theorem asserts that
$p_{b,\delta}$ is an entire function.

From Lemma~(\ref{dist:tpcommute}), we arrive at the alternative formula
\begin{align}
\begin{split}
p_{b,\delta}(x,y)
&=
\left\langle
\cT_W(\xi_1),
\left\langle
\cS(\xi_2),
\hat{f}_{b,\delta}(\xi_1,\xi_2) \exp(i (\xi_1 x + \xi_2 y))
\right\rangle
\right\rangle\\
&=
\left\langle
\cT_W(\xi_1),
\cF_{y}^{-1}\Bigl[ \hat{f}_{b,\delta} \Bigr](\xi_1,y) \exp(i \xi_1 x)
\right\rangle\\
&=
\princip \int_{-\infty}^{\infty}
\frac{\cF_{y}^{-1}\Bigl[ \hat{f}_{b,\delta} \Bigr](\xi_1,y)}{i (\xi_1 + W)}
\exp(i \xi_1 x)
\, \mathrm{d}\xi_1
.
\end{split}
\end{align}
Then we see that, for all $(x,y) \in \Omega$,
\begin{align}
\begin{split}
&\pp{p_{b,\delta}}{x}(x,y) + i W p_{b,\delta}(x,y)\\
&=
\princip \int_{-\infty}^{\infty}
\frac
{i (\xi_1 + W) \cF_{y}^{-1}\Bigl[ \hat{f}_{b,\delta} \Bigr](\xi_1,y)}
{i (\xi_1 + W)}
\exp(i \xi_1 x)
\, \mathrm{d}\xi_1\\
&=f_{b,\delta}(x,y)
,
\end{split}
\end{align}
so the second condition listed above follows from this and property~(1) in
Definition~\ref{approx:def}.

To establish the remaining conditions listed above, we will make use of the
following bounds on the partial inverse Fourier transform of
$\hat{f}_{b,\delta}$ and its derivative. From property~(3) in
Definition~\ref{approx:def}, we have
\begin{align}
\left|
\cF_{y}^{-1}\Bigl[ \hat{f}_{b,\delta} \Bigr](\xi_1,y)
\right|
\leq
20 (W_0 + \delta) \left\| f \right\|_{L^\infty(\Omega)}
,
\end{align}
and
\begin{align}
\left|
\pp{}{\xi_1} \cF_{y}^{-1}\Bigl[ \hat{f}_{b,\delta} \Bigr](\xi_1,y)
\right|
\leq
20 (W_0 + \delta) \left\| f \right\|_{L^\infty(\Omega)}
.
\end{align}
Then an analogous proof to that of Lemma~4 in \cite{Chen2024} shows
that
\begin{align}
\left\| p_{b,\delta} \right\|_{L^\infty(\Omega)}
\leq
C(W_0 + \delta)
\min\left\{ 1, \frac{1}{|W|} \right\}
\left\| f \right\|_{L^\infty(\Omega)}
,
\end{align}
and
\begin{align}
\left\| \pp{p_{b,\delta}}{x} \right\|_{L^\infty(\Omega)}
\leq
C(W_0 + \delta)
\min\left\{ 1, \frac{1}{|W|} \right\}
\left\| f \right\|_{L^\infty(\Omega)}
,
\end{align}
where $C$ is a right-continuous function.
\end{proof}

The construction of $p_{b,\delta}$ outlined above fails in the case of
$\delta = 0$,
since the definition of $c_f(\varepsilon)$ does not guarantee the existence of
an $\varepsilon$-bandlimited approximation of $f$ with bandlimit
$c_f(\varepsilon)$.
Nonetheless, the topological structure of $\sS'(\R^2)$ ensures the existence of
a slowly-varying, approximate solution $p_b$ to the Levin PDE with bandlimit
$c_f(\varepsilon)$.
In the following corollary, we present the necessary topological arguments to
construct~$p_b$.

\begin{corollary}
\label{leveq:cor1}
Under the assumptions of Lemma~\ref{leveq:lemma1},
there exist an entire function $p_b \in \sS'(\R^2)$ and a constant $C(W_0)$
such that $p_b$ satisfies the properties of $p_{b,\delta}$ in
Lemma~\ref{leveq:lemma1} for $\delta = 0$,
with the exception that $\hat{p}_b$ is of order at most $1$,
rather than of order exactly $1$.
\end{corollary}

\begin{proof}
For each $n \in \N$, we invoke Lemma~\ref{leveq:lemma1} to obtain an
$\varepsilon$-bandlimited approximation $f_{b,1/n} \in \sS(\R^2)$ of $f$ with
bandlimit $W_0 + \frac{1}{n}$ and an entire function $p_{b,1/n} \in \sS'(\R^2)$
such that
\begin{enumerate}
\item
$\hat{p}_{b,1/n}$ is a tempered distribution of order $1$ supported in
$\left[ -W_0 - \frac{1}{n}, W_0 + \frac{1}{n} \right]^2$,

\item
$
\begin{aligned}
\left|
\pp{p_{b,1/n}}{x}(x,y) + i W p_{b,1/n}(x,y) - f(x,y)
\right|
\leq
\varepsilon \left\| f \right\|_{L^\infty(\Omega)}
\end{aligned}
$
for all $(x,y) \in \Omega$,

\item
$
\begin{aligned}
\left\| p_{b,1/n} \right\|_{L^\infty(\Omega)}
\leq
C\left( W_0 + \frac{1}{n} \right)
\min\left\{ 1, \frac{1}{|W|} \right\}
\left\| f \right\|_{L^\infty(\Omega)}
\end{aligned}
,
$

\item 
$
\begin{aligned}
\left\| \pp{p_{b,1/n}}{x} \right\|_{L^\infty(\Omega)}
\leq
C\left( W_0 + \frac{1}{n} \right)
\min\left\{ 1, \frac{1}{|W|} \right\}
\left\| f \right\|_{L^\infty(\Omega)}
\end{aligned}
.
$
\end{enumerate}

Recall that the weak dual topology on $\sS'(\R^2)$, which is the topology of
pointwise convergence in $\sS(\R^2)$, is defined by the family of seminorms
$\{ q_\varphi \colon \varphi \in \sS(\R^2) \}$, where
$q_\varphi(\cT) = |\langle \cT, \varphi \rangle|$ for $\cT \in \sS'(\R^2)$.

To see that the sequence $(\hat{p}_{b,1/n})_{n=1}^{\infty}$ is weakly bounded,
we observe that
\begin{align}
\label{leveq:qphibound}
\begin{split}
q_\varphi(\hat{p}_{b,1/n})
&=
\left|
\int_{-W_0 - \frac{1}{n}}^{W_0 + \frac{1}{n}}
\left\langle
\cT_W(\xi_1), \hat{f}_{b,1/n}(\xi_1,\xi_2) \varphi(\xi_1,\xi_2)
\right\rangle
\mathrm{d}\xi_2
\right|\\
&\leq
4 \left( W_0 + \frac{1}{n} \right)
\left(
\left\| \hat{f}_{b,1/n} \varphi \right\|_{\bm{e}_1,\bm{0}}
+
\left\| \hat{f}_{b,1/n} \varphi \right\|_{\bm{0},\bm{e}_1}
\right)\\
&\leq
120 \left( W_0 + 1 \right)^2
\left\| f \right\|_{L^\infty(\Omega)}
\sup_{|\bm{\beta}| \leq 1}
\left\| \varphi \right\|_{\bm{0},\bm{\beta}}
,
\end{split}
\end{align}
where $\bm{e}_1 = (1,0)$, and where the first inequality follows from the fact
that
\begin{align}
\left|
\left\langle
\princip \, \frac{1}{x}, \varphi
\right\rangle
\right|
\leq
2
\left(
\left\| \varphi \right\|_{1,0} + \left\| \varphi \right\|_{0,1}
\right)
\end{align}
and second inequality follows from
property~(3) in Definition~\ref{approx:def}.
Since $\sS(\R^2)$ is barreled, it follows that, since
$(\hat{p}_{b,1/n})_{n=1}^{\infty}$ is weakly bounded, it is also relatively
compact in the weak dual topology (see Theorem~33.2 of \cite{Treves2006}),
and hence contains a convergent subnet $(\hat{p}_{b,1/F(\nu)})_{\nu \in \cI}$
for some monotone final function $F \colon \cI \to \N$ and directed system $\cI$.
We denote by $\hat{p}_b$ the limit of the subnet,
and define $p_b$ as the inverse Fourier transform of $\hat{p}_b$.

From (\ref{leveq:qphibound}) and the pointwise convergence of
$(\hat{p}_{b,1/F(\nu)})_{\nu\in\cI}$, we observe that
\begin{align}
\left|
\langle
\hat{p}_b, \varphi
\rangle
\right|
\leq
120 \left( W_0 + 1 \right)^2
\left\| f \right\|_{L^\infty(\Omega)}
\sup_{|\bm{\beta}| \leq 1}
\left\| \varphi \right\|_{\bm{0},\bm{\beta}}
\end{align}
for all $\varphi \in \sS(\R^2)$, which implies that
$\hat{p}_b$ is of order at most $1$.
For each test function $\varphi$ whose support is disjoint from $[-W_0,W_0]^2$,
there exists some $n \in \N$ such that $\hat{p}_{b,1/n}$ and $\varphi$ have
disjoint support, so we can find some $\nu' \in \cI$ such that
$q_\varphi(\hat{p}_{b,1/F(\nu)}) = 0$ for all $\nu \succ \nu'$.
This implies that the limit $q_\varphi(\hat{p}_b)$ must be zero, and thus the
support of $\hat{p}_b$ is contained in $[-W_0,W_0]^2$.
It is then a direct consequence of the Schwartz-Paley-Wiener theorem that $p_b$
is an entire function.

Owing to the formula
\begin{align}
p_{b,1/F(\nu)}(x,y)
=
\left\langle
\hat{p}_{b,1/F(\nu)}, \exp(i (\xi_1 x + \xi_2 y))
\right\rangle
,
\end{align}
we conclude that $(D^{\bm{\alpha}}p_{b,1/F(\nu)})_{\nu\in\cI}$ converges
pointwise to $D^{\bm{\alpha}}p_b$ for all multi-indices $\bm{\alpha}$.
By combining the pointwise convergence of
$D^{\bm{\alpha}}p_{b,1/F(\nu)}$ with the
right-continuity of $C$, the remaining conditions follow immediately.
\end{proof}

We now consider the case where the domain is $[-1,1]^2$, and $g_x$ is
nonconstant with no zeros over the domain.
We suppose that $f \colon [-1,1]^2 \to \C$ and $g \colon [-1,1]^2 \to \R$ admit
infinitely differentiable extensions to an open neighborhood of $[-1,1]^2$, and
that the extension of $g_x$ is nonzero in the neighborhood.
We let
\begin{align}
G_0
=
\min_{(x,y) \in [-1,1]^2}
\left| \pp{g}{x}(x,y) \right|
,
\quad
G_1
=
\max_{(x,y) \in [-1,1]^2}
\left| \pp{g}{x}(x,y) \right|
,
\end{align}
and
\begin{align}
W
=
\max_{y \in [-1,1]}
\left| \int_{-1}^{1} \pp{g}{x}(x,y) \, \mathrm{d}x \right|
.
\end{align}
Furthermore, we define $\bm{u} \colon [-1,1]^2 \to \R^2$ via the formula
\begin{align}
\label{leveq:defu}
\bm{u}(x,y)
=
\left(
\frac{1}{W} \int_{-1}^{x} \pp{g}{x}(t,y) \, \mathrm{d}t
,\,
y
\right)
.
\end{align}
Since $g_x$ is nonvanishing on $[-1,1]^2$,
$\bm{u}$ is injective and has invertible Jacobian on $[-1,1]^2$.
By the generalized inverse function theorem
(see Chapter~1 of \cite{Guillemin2010}),
$\bm{u}$ extends to a diffeomorphism in an open neighborhood of $[-1,1]^2$.
We denote the image of $(-1,1)^2$ under $\bm{u}$ by $\Omega_{\bm{u}}$,
noting that $\overbar{\Omega}_{\bm{u}} \subset [-1,1]^2$,
and define $h \colon \overbar{\Omega}_{\bm{u}} \to \R^2$ via the formula
\begin{align}
\label{leveq:defh}
h(x',y')
=
f(\bm{u}^{-1}(x',y')) \det(D\bm{u}^{-1})(x',y')
.
\end{align}
Notice that both $\bm{u}$ and $h$ are independent of the magnitude of $g_x$,
in the sense that they remain unchanged under any rescaling of $g_x$ by a
nonzero constant.

We are now ready to prove the first main result of this section.

\begin{theorem}
\label{leveq:thmhigh}
Suppose that $0 < \varepsilon < 1$, and let $W_0 = c_h(\varepsilon)$.
Then there exists an infinitely differentiable function
$p_b \colon [-1,1]^2 \to \C$ and a constant $C(W_0)$, depending only on $W_0$,
such that
\begin{enumerate}
\item
the Fourier transform of $p_b(\bm{u}^{-1}(x',y'))$ is a tempered distribution of
order at most $1$ supported in $[-W_0,W_0]^2$,

\item
$
\begin{aligned}
\left|
\pp{p_b}{x}(x,y) + i \pp{g}{x}(x,y) p_b(x,y) - f(x,y)
\right|
\leq
\varepsilon \frac{G_1}{G_0} \left\| f \right\|_{L^\infty([-1,1]^2)}
\end{aligned}
$
for all $(x,y) \in [-1,1]^2$,

\item
$
\begin{aligned}
\left\| p_{b} \right\|_{L^\infty([-1,1]^2)}
\leq
C(W_0)
\frac{W}{G_0}
\min\left\{ 1, \frac{1}{W} \right\}
\left\| f \right\|_{L^\infty([-1,1]^2)}
\end{aligned}
,
$

\item
$
\begin{aligned}
\left\| \pp{p_{b}}{x} \right\|_{L^\infty([-1,1]^2)}
\leq
C(W_0)
\frac{G_1}{G_0}
\min\left\{ 1, \frac{1}{W} \right\}
\left\| f \right\|_{L^\infty([-1,1]^2)}
\end{aligned}
.
$
\end{enumerate}
\end{theorem}

\begin{proof}
Using the change of variables $(x',y') = \bm{u}(x,y)$, the PDE
\begin{align}
\pp{p}{x}(x,y) + i \pp{g}{x}(x,y) p(x,y) = f(x,y) \quad \text{in $(-1,1)^2$}
\end{align}
can be reduced to the simpler form
\begin{align}
\label{leveq:reducedpde}
\pp{p}{x'}(x',y') + i W p(x',y') = h(x',y') \quad \text{in $\Omega_{\bm{u}}$}.
\end{align}
Since $h$ admits an infinitely differentiable extension to an open neighborhood
of $\overbar{\Omega}_{\bm{u}}$, an application of Corollary~\ref{leveq:cor1} to
(\ref{leveq:reducedpde}) shows that there exist an entire function $p_1$ and a
constant $C(W_0)$ such that
\begin{enumerate}
\item
$\hat{p}_1$ is a tempered distribution of order at most $1$ supported in
$[-W_0,W_0]^2$,

\item
$
\begin{aligned}
\left|
\pp{p_1}{x'}(x',y') + i W p_1(x',y') - h(x',y')
\right|
\leq
\varepsilon \left\| h \right\|_{L^\infty(\overbar{\Omega}_{\bm{u}})}
\end{aligned}
$
for all $(x',y') \in \overbar{\Omega}_{\bm{u}}$,

\item
$
\begin{aligned}
\left\| p_{1} \right\|_{L^\infty(\overbar{\Omega}_{\bm{u}})}
\leq
C(W_0)
\min\left\{ 1, \frac{1}{W} \right\}
\left\| h \right\|_{L^\infty(\overbar{\Omega}_{\bm{u}})}
\end{aligned}
,
$

\item
$
\begin{aligned}
\left\| \pp{p_{1}}{x'} \right\|_{L^\infty(\overbar{\Omega}_{\bm{u}})}
\leq
C(W_0)
\min\left\{ 1, \frac{1}{W} \right\}
\left\| h \right\|_{L^\infty(\overbar{\Omega}_{\bm{u}})}
\end{aligned}
.
$
\end{enumerate}

We define $p_b$ as the composition of $p_1$ and $\bm{u}$.
It is clear that the first condition on $p_b$ is satisfied.
From (\ref{leveq:defu}) and (\ref{leveq:defh}), we see that
\begin{align}
\left\| h \right\|_{L^\infty(\overbar{\Omega}_{\bm{u}})}
\leq
\frac{W}{G_0}
\left\| f \right\|_{L^\infty([-1,1]^2)}
,
\end{align}
and therefore
\begin{align}
\left\| p_{b} \right\|_{L^\infty([-1,1]^2)}
\leq
C(W_0)
\frac{W}{G_0}
\min\left\{ 1, \frac{1}{W} \right\}
\left\| f \right\|_{L^\infty([-1,1]^2)}
,
\end{align}
and
\begin{align}
\left\| \pp{p_{b}}{x} \right\|_{L^\infty([-1,1]^2)}
\leq
C(W_0)
\frac{G_1}{G_0}
\min\left\{ 1, \frac{1}{W} \right\}
\left\| f \right\|_{L^\infty([-1,1]^2)}
.
\end{align}
Moreover, we have
\begin{align}
\left| \pp{p_b}{x}(x,y) + i \pp{g}{x}(x,y) p_b(x,y) - f(x,y) \right|
\leq
\varepsilon \frac{G_1}{G_0} \left\| f \right\|_{L^\infty([-1,1]^2)}
\end{align}
for all $(x,y) \in [-1,1]^2$.
\end{proof}

We note that Theorem~\ref{leveq:thmhigh} does not assert that the function $p_b$
is bandlimited.
However, we conclude that its smoothness is characterized independently of the
magnitude of $g_x$.
This follows from the observation that $p_1$ has a bandlimit independent of the
magnitude of $g_x$ and that the function $\bm{u}$ is also independent of the
magnitude of $g_x$.
It is clear that $p_b$ can be approximated to a fixed relative accuracy via a
Legendre expansion at a cost that does not depend on the magnitude of $g_x$, but
which necessarily depends on the complexity of $\bm{u}$ and $h$.

We close this section with the second main result of this section, which proves
the existence of a slowly-varying, approximate solution to the Levin PDE if $g_x$
is of small magnitude.

\begin{theorem}
\label{leveq:thmlow}
Suppose that $f \colon [-1,1]^2 \to \C$ and $g \colon [-1,1]^2 \to \R$ admit
infinitely differentiable extensions to an open neighborhood of $[-1,1]^2$, and
that $G_1 < 1/2$.
Then, for $0 < \varepsilon < 1$ and integer $n$ defined via the formula
\begin{align}
\label{leveq:defn}
n = \left\lfloor \frac{\log(\varepsilon)}{\log(2 G_1)} \right\rfloor,
\end{align}
there exists an entire function $p_b \in \sS'(\R^2)$ such that
\begin{enumerate}
\item
$\hat{p}_b$ is a tempered distribution of order at most $1$ supported in
\begin{align*}
\left[
-c_f(\varepsilon) - n c_{g_x}(\varepsilon),
 c_f(\varepsilon) + n c_{g_x}(\varepsilon)
\right]^2
,
\end{align*}

\item
$
\begin{aligned}
\left|
\pp{p_{b}}{x}(x,y) + i \pp{g}{x}(x,y) p_{b}(x,y) - f(x,y)
\right|
\leq
3 \varepsilon \left( 1 + \frac{G_1}{1 - 2 G_1} \right)
\left\| f \right\|_{L^\infty([-1,1]^2)}
\end{aligned}
$
for all $(x,y) \in [-1,1]^2$,

\item
$
\begin{aligned}
\left\| p_{b} \right\|_{L^\infty([-1,1]^2)}
\leq
\frac{2}{1 - 2 G_1}
\left\| f \right\|_{L^\infty([-1,1]^2)}
\end{aligned}
,
$

\item
$
\begin{aligned}
\left\| \pp{p_{b}}{x} \right\|_{L^\infty([-1,1]^2)}
\leq
4 \left( 1 + \frac{G_1}{1 - 2 G_1} \right)
\left\| f \right\|_{L^\infty([-1,1]^2)}
\end{aligned}
.
$
\end{enumerate}
\end{theorem}

\begin{proof}
We will first construct an entire function $p_{b,\delta} \in \sS'(\R^2)$ for
$\delta > 0$ such that $p_{b,\delta}$ satisfies conditions (2)-(4) listed above,
and $\hat{p}_{b,\delta}$ is a tempered distribution of order at most $1$
supported in
\begin{align}
\left[
-c_f(\varepsilon) - n c_{g_x}(\varepsilon) - (n+1) \delta,
 c_f(\varepsilon) + n c_{g_x}(\varepsilon) + (n+1) \delta
\right]^2
.
\end{align}
We let $f_{b,\delta}$ and $(g_{x})_{b,\delta}$ be
$\varepsilon$-bandlimited approximations of $f$ and $g_x$ with
bandlimit $c_{f}(\varepsilon) + \delta$ and $c_{g_{x}}(\varepsilon) + \delta$,
respectively.
We then define $A_{\delta} \colon L^\infty([-1,1]^2) \to L^\infty([-1,1]^2)$ via
the formula
\begin{align}
A_\delta[\varphi](x,y)
=
-i
\int_{0}^{x}
\left(\pp{g}{x}\right)_{b,\delta}(t,y) \, \varphi(t,y)
\, \mathrm{d}t
,
\end{align}
and define
\begin{align}
h_\delta(x,y) = \int_{0}^{x} f_{b,\delta}(t,y) \, \mathrm{d}t
.
\end{align}
The function $p_{b,\delta}$ is then defined as
\begin{align}
\label{leveq:defpiter}
p_{b,\delta}(x,y) = \sum_{k=0}^{n} A_{\delta}^{k}[h_\delta](x,y)
,
\end{align}
where $A_\delta^k$ denotes the repeated application of the operator $A_\delta$
and $n$ is defined by (\ref{leveq:defn}). By the same derivation presented in
the proof of Theorem~5 in \cite{Chen2024}, we see that
\begin{align}
\begin{split}
&
\left|
\pp{p_{b,\delta}}{x}(x,y) + i \pp{g}{x}(x,y) p_{b,\delta}(x,y) - f(x,y)
\right|\\
&\leq
3 \varepsilon \left( 1 + \frac{G_1}{1 - 2 G_1} \right)
\left\| f \right\|_{L^\infty([-1,1]^2)}
\end{split}
\end{align}
for all $(x,y) \in [-1,1]^2$,
\begin{align}
\left\| p_{b,\delta} \right\|_{L^\infty([-1,1]^2)}
&\leq
\frac{2}{1 - 2 G_1}
\left\| f \right\|_{L^\infty([-1,1]^2)}
,
\end{align}
and
\begin{align}
\left\| \pp{p_{b,\delta}}{x} \right\|_{L^\infty([-1,1]^2)}
&\leq
4 \left( 1 + \frac{G_1}{1 - 2 G_1} \right)
\left\| f \right\|_{L^\infty([-1,1]^2)}
.
\end{align}

It remains to verify the bandlimit of $p_{b,\delta}$ and the order of its
Fourier transform.
We first observe that
\begin{align}
\label{leveq:hdeltaf}
\hat{h}_\delta
=
\hat{f}_{b,\delta} (\cT_0 \otimes \cS)
-
\delta(\xi_1)
\otimes
\left(
\left\langle \cT_0(\xi), \hat{f}_{b,\delta}(\xi,\xi_2) \right\rangle
\cS(\xi_2)
\right)
,
\end{align}
where $\cT_0$ is defined by (\ref{leveq:deftw}), and $\cS$ denotes integration.
Since $f_{b,\delta}$ has bandlimit $c_f(\varepsilon) + \delta$, it is clear that
$h_\delta$ has the same bandlimit.
Analogously, we have
\begin{align}
\label{leveq:adeltaf}
\hspace*{-2em}
\hat{A_{\delta}[\varphi]}
=
-i
\left(
\hat{I_{\delta}[\varphi]} (\cT_0 \otimes \cS)
-
\delta(\xi_1)
\otimes
\left(
\left\langle \cT_0(\xi), \hat{I_{\delta}[\varphi]}(\xi,\xi_2) \right\rangle
\cS(\xi_2)
\right)
\right)
,
\end{align}
where $I_{\delta}[\varphi]$ denotes the integrand of $A_{\delta}[\varphi]$.
Suppose that $\varphi$ is a slowly increasing $C^\infty$ function with
bandlimit $c$. Then, as a consequence of Lemma~\ref{dist:conv} and the
fact that ${(g_x)_{b,\delta}} \in \sS(\R^2)$ has bandlimit
$c_{g_{x}}(\varepsilon) + \delta$,
\begin{align}
\hat{I_{\delta}[\varphi]}(\bm{x})
=
\hat{\varphi} \ast \hat{(g_x)}_{b,\delta}(\bm{x})
\end{align}
defines a $C^\infty$ function supported in
$[-c - c_{g_{x}}(\varepsilon) - \delta, c + c_{g_{x}}(\varepsilon) + \delta]^2$.
Hence it is clear that $A_{\delta}[\varphi]$ has bandlimit
$c + c_{g_{x}}(\varepsilon) + \delta$.
It then follows by induction that $A_\delta^k[h_\delta]$ has bandlimit
$c_f(\varepsilon) + k c_{g_{x}}(\varepsilon) + (k+1) \delta$.
Moreover, $A_\delta^k[h_\delta]^\wedge$ is of order at most $1$ since $\cT_0$ is
of order $1$.
We conclude that $p_{b,\delta}$ has bandlimit
$c_f(\varepsilon) + n c_{g_{x}}(\varepsilon) + (n+1) \delta$, and the order of 
its Fourier transform is at most $1$.

Recall that the weak dual topology on $\sS'(\R^2)$ is defined by the family of
seminorms $\{ q_\varphi \colon \varphi \in \sS(\R^2) \}$, where
$q_\varphi(\cT) = |\langle \cT, \varphi \rangle|$ for $\cT \in \sS'(\R^2)$.
To see that the sequence $(\hat{h}_{1/m})_{m=1}^{\infty}$ is weakly bounded,
we first observe that
\begin{align}
\begin{split}
&\left\langle
\delta(\xi_1)
\otimes
\left(
\left\langle \cT_0(\xi), \hat{f}_{b,1/m}(\xi,\xi_2) \right\rangle
\cS(\xi_2)
\right)
,
\varphi(\xi_1,\xi_2)
\right\rangle\\
&=
\left\langle
\hat{f}_{b,1/m}(\xi,\xi_2) (\cT_0(\xi) \otimes \cS(\xi_2)), \varphi(0,\xi_2)
\right\rangle
,
\end{split}
\end{align}
and then apply the same reasoning involved in (\ref{leveq:qphibound}) to
(\ref{leveq:hdeltaf}) to arrive at the bound
\begin{align}
q_\varphi(\hat{h}_{1/m})
\leq
240 \left( c_f(\varepsilon) + 1 \right)^2
\left\| f \right\|_{L^\infty([-1,1]^2)}
\sup_{|\bm{\beta}| \leq 1}
\left\| \varphi \right\|_{\bm{0},\bm{\beta}}
.
\end{align}
Suppose now that $(\hat{\psi}_m)_{m=1}^{\infty}$ is a sequence in
$\sE'(\R^2)$ with bandlimit $c$ satisfying
\begin{align}
\label{leveq:weakbdd}
q_\varphi(\hat{\psi}_m)
\leq
C
\sup_{|\bm{\beta}| \leq 1}
\left\| \varphi \right\|_{\bm{0},\bm{\beta}}
\end{align}
for some constant $C$ independent of $m$ and $\varphi$.
From (\ref{dist:convp}), (\ref{leveq:weakbdd}) and property~(3) of
$((g_x)_{b,1/m})_{m=1}^{\infty}$ in Definition~\ref{approx:def},
we see that $(I_{1/m}[\psi_m]^{\wedge})_{m=1}^{\infty}$ and their first
derivatives are all uniformly bounded, and are supported in
$[-c-c_{g_{x}}(\varepsilon)-1,c+c_{g_{x}}(\varepsilon)+1]^2$. It follows from
the same arguments again that $(A_{1/m}[\psi_m]^{\wedge})_{m=1}^{\infty}$ is
weakly bounded and satisfies a bound of the form (\ref{leveq:weakbdd}).

By induction, we see that $(A_{1/m}^{k}[h_{1/m}]^{\wedge})_{m=1}^{\infty}$
is weakly bounded for all positive $k$.
Consequently, $(\hat{p}_{b,1/m})_{m=1}^{\infty}$ is weakly bounded due to
subadditivity of seminorms. The existence of an entire function $p_b$ satisfying
the conditions listed above can be established by the same arguments used in the
proof of Corollary~\ref{leveq:cor1}.
\end{proof}

\section{Numerical aspects of the Levin method}

\label{section:levnum}

In this section, we derive error estimates for the Levin method that discretizes
the PDE
\begin{align}
\label{levnum:pde}
p_x + i g_x p = f
\quad \text{in $(-1,1)^2$}
\end{align}
using collocation on a tensor product Chebyshev grid and then solves the
resulting linear system numerically using a truncated singular value
decomposition.

In this section, we use the symbol $\hat{x}$ to denote the computed
approximation to $x$, rather than the Fourier transform of $x$.

Throughout this section, we will assume that $f \colon [-1,1]^2 \to \C$ and
$g \colon [-1,1]^2 \to \R$ admit infinitely differentiable extensions to an open
neighborhood of $[-1,1]^2$, and that $0 < \varepsilon < 1/2$.
As in the preceding section, we let
\begin{align}
G_0
=
\min_{(x,y) \in [-1,1]^2}
\left| \pp{g}{x}(x,y) \right|
,
\quad
G_1
=
\max_{(x,y) \in [-1,1]^2}
\left| \pp{g}{x}(x,y) \right|
,
\end{align}
and
\begin{align}
W
=
\max_{y \in [-1,1]}
\left| \int_{-1}^{1} \pp{g}{x}(x,y) \, \mathrm{d}x \right|
.
\end{align}

\subsection[Case I: G0 > 0]{Case I: $G_0 > 0$}

We let $W_0 = c_h(\varepsilon)$, where $h$ is defined via the formula
(\ref{leveq:defh}) and $c_h(\varepsilon)$ is defined in
Definition~\ref{approx:def}.
According to Theorem~\ref{leveq:thmhigh}, there exist an infinitely
differentiable function $p_b \colon [-1,1]^2 \to \C$ and a constant
$C(W_0)$, which we assume to satisfy $C(W_0) \geq 1$ in this subsection,
such that
\begin{align}
\label{levnum:pbhigherr}
\left|
\pp{p_b}{x}(x,y) + i \pp{g}{x}(x,y) p_b(x,y) - f(x,y)
\right|
\leq
\varepsilon \, \frac{G_1}{G_0} \left\| f \right\|_{L^\infty([-1,1]^2)}
\end{align}
for all $(x,y) \in [-1,1]^2$, and both
\begin{align}
\label{levnum:pbhighbound}
\left\| p_{b} \right\|_{L^\infty([-1,1]^2)}
&\leq
C(W_0) \,
\frac{W}{G_0}
\min\left\{ 1, \frac{1}{W} \right\}
\left\| f \right\|_{L^\infty([-1,1]^2)}
\end{align}
and
\begin{align}
\label{levnum:pbphighbound}
\left\| \pp{p_{b}}{x} \right\|_{L^\infty([-1,1]^2)}
&\leq
C(W_0) \,
\frac{G_1}{G_0}
\min\left\{ 1, \frac{1}{W} \right\}
\left\| f \right\|_{L^\infty([-1,1]^2)}
.
\end{align}
From (\ref{cheb:ftail}), (\ref{cheb:fptail}), and the discussion following
Theorem~\ref{leveq:thmhigh}, we observe that an integer $k$ can be chosen
independently of the magnitude of $g_x$ such that
\begin{align}
\label{levnum:pbcheb}
\left\| \mathsf{P}_k[p_b] - p_b \right\|_{L^\infty([-1,1]^2)}
\leq
\varepsilon
\left\| p_b \right\|_{L^\infty([-1,1]^2)}
,
\end{align}
and
\begin{align}
\label{levnum:pbpcheb}
\left\| \pp{}{x}\mathsf{P}_k[p_b] - \pp{p_b}{x} \right\|_{L^\infty([-1,1]^2)}
\leq
\varepsilon
\left\| \pp{p_b}{x} \right\|_{L^\infty([-1,1]^2)}
.
\end{align}
Moreover,
from the uniform convergence of Chebyshev interpolants, the value of $k$ can be
chosen such that
\begin{align}
\label{levnum:gxcheb}
\left\|
\mathsf{P}_k\biggl[ \pp{g}{x} \biggr] - \pp{g}{x}
\right\|_{L^\infty([-1,1]^2)}
\leq
\varepsilon \left\| \pp{g}{x} \right\|_{L^\infty([-1,1]^2)}
,
\end{align}
and
\begin{align}
\label{levnum:fcheb}
\left\|
\mathsf{P}_{2k-1}[f] - f
\right\|_{L^\infty([-1,1]^2)}
\leq
\varepsilon \left\| f \right\|_{L^\infty([-1,1]^2)}
\end{align}
are also satisfied.
With this choice of $k$ and the assumption that $0 < \varepsilon < 1/2$, we have
\begin{align}
\label{levnum:prodcheb}
\begin{split}
&
\left\|
\mathsf{P}_k\biggl[ \pp{g}{x} \biggr] \mathsf{P}_k[p_b] - \pp{g}{x} p_b
\right\|_{L^\infty([-1,1]^2)}\\
&\leq
\left\|
\mathsf{P}_k[p_b]
\left(
\mathsf{P}_k\biggl[ \pp{g}{x} \biggr] - \pp{g}{x}
\right)
+
\pp{g}{x}
\left(
\mathsf{P}_k[p_b] - p_b
\right)
\right\|_{L^\infty([-1,1]^2)}\\
&\leq
3 \varepsilon
\left\| \pp{g}{x} \right\|_{L^\infty([-1,1]^2)}
\left\| p_b \right\|_{L^\infty([-1,1]^2)}
.
\end{split}
\end{align}
From (\ref{levnum:pbhighbound}), (\ref{levnum:pbphighbound}),
(\ref{levnum:pbpcheb}), (\ref{levnum:prodcheb}),
and the fact that
\begin{align}
\max\{W, 1\} \min\left\{ 1, \frac{1}{W} \right\} = 1
,
\end{align}
we see that
\begin{align}
\label{levnum:pbhighcheb}
\begin{split}
&
\left\|
\pp{}{x}\mathsf{P}_k[p_b]
+
i \, \mathsf{P}_k\biggl[ \pp{g}{x} \biggr] \mathsf{P}_k[p_b]
-
\left(
\pp{p_b}{x} + i \pp{g}{x} p_b
\right)
\right\|_{L^\infty([-1,1]^2)}\\
&\leq
3 \varepsilon
\left\| \pp{g}{x} \right\|_{L^\infty([-1,1]^2)}
\left\| p_b \right\|_{L^\infty([-1,1]^2)}
+
\varepsilon \left\| \pp{p_b}{x} \right\|_{L^\infty([-1,1]^2)}\\
&\leq
4 \, \varepsilon \, C(W_0) \, \frac{G_1}{G_0}
\left\| f \right\|_{L^\infty([-1,1]^2)}
\end{split}
\end{align}
holds regardless of the magnitude of $g_x$.
Combining (\ref{levnum:pbhigherr}) with (\ref{levnum:fcheb}) and
(\ref{levnum:pbhighcheb}) yields
\begin{align}
\label{levnum:highcheberr}
\begin{split}
&
\left|
\pp{}{x}\mathsf{P}_k[p_b](x,y)
+
i \, \mathsf{P}_k\biggl[ \pp{g}{x} \biggr](x,y) \, \mathsf{P}_k[p_b](x,y)
-
\mathsf{P}_{2k-1}[f](x,y)
\right|\\
&\leq
\varepsilon \left( 1 + \left( 1 + 4 C(W_0) \right) \frac{G_1}{G_0} \right)
\left\| f \right\|_{L^\infty([-1,1]^2)}\\
&\leq
6 \, \varepsilon \, C(W_0) \, \frac{G_1}{G_0}
\left\| f \right\|_{L^\infty([-1,1]^2)}
,
\end{split}
\end{align}
for all $(x,y) \in [-1,1]^2$, since $C(W_0) \geq 1$.
If we let
\begin{align}
\bm{\mathsf{p_b}}
=
\begin{pmatrix}
p_b\bigl(\bm{x}_{1,k^2}\bigr) \\
p_b\bigl(\bm{x}_{2,k^2}\bigr) \\
\vdots \\
p_b\bigl(\bm{x}_{k^2,k^2}\bigr)
\end{pmatrix}
,
\quad
\bm{\mathsf{f}}
=
\begin{pmatrix}
f\bigl(\bm{x}_{1,(2k-1)^2}\bigr) \\
f\bigl(\bm{x}_{2,(2k-1)^2}\bigr) \\
\vdots \\
f\bigl(\bm{x}_{(2k-1)^2,(2k-1)^2}\bigr)
\end{pmatrix}
,
\end{align}
and define the $(2k-1)^2 \times (2k-1)^2$ diagonal matrix
$\bm{\mathsf{G}}_{2k-1}$ with entries given by
\begin{align}
(\bm{\mathsf{G}}_{2k-1})_{i,i}
=
g_x\bigl(\bm{x}_{i,(2k-1)^2}\bigr)
,
\end{align}
for all $1 \leq i \leq (2k-1)^2$,
then it follows from (\ref{levnum:highcheberr}) that
\begin{align}
\label{levnum:linsys}
\left(
\bm{\mathsf{P}}_{(2k-1,k)} (\bm{\mathsf{D}}_{x})_k
+
i \bm{\mathsf{G}}_{2k-1} \bm{\mathsf{P}}_{(2k-1,k)}
\right)
\bm{\mathsf{p_b}}
=
\bm{\mathsf{f}} + \bm{\mathsf{r}}
\end{align}
with
\begin{align}
\label{levnum:pbnorm}
\| \bm{\mathsf{p_b}} \|_2
\lesssim
C(W_0) \, \frac{W}{G_0}
\min\left\{ 1, \frac{1}{W} \right\}
\left\| f \right\|_{L^\infty([-1,1]^2)}
\end{align}
and
\begin{align}
\label{levnum:dfnorm}
\| \bm{\mathsf{r}} \|_2
\lesssim
\varepsilon \, C(W_0) \, \frac{G_1}{G_0}
\left\| f \right\|_{L^\infty([-1,1]^2)}
,
\end{align}
where $(\bm{\mathsf{D}}_{x})_k$ and $\bm{\mathsf{P}}_{(2k-1,k)}$ are the
Chebyshev differentiation and interpolation matrices defined in
Section~\ref{cheb}.

We introduce the notation
\begin{align}
\bm{\mathsf{A}}
=
\bm{\mathsf{P}}_{(2k-1,k)} (\bm{\mathsf{D}}_{x})_k
+
i \bm{\mathsf{G}}_{2k-1} \bm{\mathsf{P}}_{(2k-1,k)}
\end{align}
in order to simplify the following discussion.
We observe that
\begin{align}
\label{levnum:anormupper}
\|\bm{\mathsf{A}}\|_2
\lesssim
\max\{G_1, 1\}
,
\end{align}
and, by combining Lemma~\ref{cheb:lebesgue} with (\ref{levnum:gxcheb}),
we see that
\begin{align}
\label{levnum:anormlower}
\|\bm{\mathsf{A}}\|_2
\geq
\frac{1}{\sqrt{2} k}
\left(
1 + \frac{(1 - \varepsilon) G_1}{\Lambda_k^2}
\right)
\geq
\frac{1 + G_1}{2 \sqrt{2} k \Lambda_k^2}
\geq
\frac{\max\{G_1, 1\}}{2 \sqrt{2} k \Lambda_k^2}
,
\end{align}
and thus $\|\bm{\mathsf{A}}\|_2 \gtrsim \max\{G_1, 1\}$.
We now solve the linear system
\begin{align}
\bm{\mathsf{A}} \bm{\mathsf{p_b}} = \bm{\mathsf{f}}
\end{align}
via a truncated singular value decomposition which has been truncated at
precision $\varepsilon \|\bm{\mathsf{A}}\|_2$, and obtain a solution
\begin{align}
\hat{\bm{\mathsf{p}}}
=
(\bm{\mathsf{A}} + \bm{\mathsf{E}})_\nu^\dagger
(\bm{\mathsf{f}} + \bm{\mathsf{e}})
,
\end{align}
where
$\bm{\mathsf{E}} \in \C^{(2k-1)^2 \times k^2}$
and
$\bm{\mathsf{e}} \in \C^{(2k-1)^2}$
with
$\|\bm{\mathsf{E}}\|_2 < \varepsilon \|\bm{\mathsf{A}}\|_2 / 2$
and
$\|\bm{\mathsf{e}}\|_2 \lesssim \varepsilon \|\bm{\mathsf{f}}\|_2$,
and where $(\bm{\mathsf{A}} + \bm{\mathsf{E}})_\nu^\dagger$ is the pseudoinverse
of the $\nu$-truncated singular value decomposition of
$\bm{\mathsf{A}} + \bm{\mathsf{E}}$, so that
\begin{align}
\hat{\sigma}_\nu
\geq
\varepsilon \|\bm{\mathsf{A}}\|_2
>
\hat{\sigma}_{\nu+1}
,
\end{align}
where 
$\hat{\sigma}_1 \ge \hat{\sigma}_2 \ge
\cdots \ge \hat{\sigma}_{k^2}$ 
denote the singular values of $\bm{\mathsf{A}} + \bm{\mathsf{E}}$,
defining $\hat{\sigma}_{k^2+1} = 0$.
By Theorem~\ref{tsvd:thm}, we have that
\begin{align}
\| \hat{\bm{\mathsf{p}}} \|_2
\leq
\frac{1}{\varepsilon \|\bm{\mathsf{A}}\|_2}
(
2 \varepsilon \|\bm{\mathsf{A}}\|_2 \|\bm{\mathsf{p_b}}\|_2
+
\| \bm{\mathsf{e}} - \bm{\mathsf{r}} \|_2
)
+
\|\bm{\mathsf{p_b}}\|_2
,
\end{align}
and
\begin{align}
\| \bm{\mathsf{A}} \hat{\bm{\mathsf{p}}} - \bm{\mathsf{f}} \|_2
\leq
5 \varepsilon \|\bm{\mathsf{A}}\|_2 \|\bm{\mathsf{p_b}}\|_2
+
\frac{3}{2} \| \bm{\mathsf{e}} - \bm{\mathsf{r}} \|_2
+
\|\bm{\mathsf{r}}\|_2
.
\end{align}
Combining the above inequalities with (\ref{levnum:pbnorm}),
(\ref{levnum:dfnorm}), (\ref{levnum:anormupper}), (\ref{levnum:anormlower}) and
the fact that
\begin{align}
\max\{G_1, 1\} \min\{W, 1\}
\leq
2 G_1
,
\end{align}
we conclude that
\begin{align}
\label{levnum:pnorm}
\| \hat{\bm{\mathsf{p}}} \|_2
\lesssim
C(W_0) \, \frac{\min\{G_1, 1\}}{G_0}
\left\| f \right\|_{L^\infty([-1,1]^2)}
,
\end{align}
and
\begin{align}
\label{levnum:errnorm}
\| \bm{\mathsf{A}} \hat{\bm{\mathsf{p}}} - \bm{\mathsf{f}} \|_2
\lesssim
\varepsilon \, C(W_0) \, \frac{G_1}{G_0}
\left\| f \right\|_{L^\infty([-1,1]^2)}
.
\end{align}

We now let $\hat{p}$ be the unique bivariate polynomial of degree less than $k$
whose values on the $k \times k$ tensor product Chebyshev grid are given by the
entries of $\hat{\bm{\mathsf{p}}}$,
and let $\hat{r}$ be the unique bivariate polynomial of degree less than
$2k-1$ whose values on the $(2k-1) \times (2k-1)$ tensor product Chebyshev grid
equal to the entries of the vector
$
\hat{\bm{\mathsf{r}}}
=
\bm{\mathsf{A}} \hat{\bm{\mathsf{p}}} - \bm{\mathsf{f}}
$.
It then follows that
\begin{align}
\label{levnum:phatfun}
\pp{\hat{p}}{x}(x,y)
+
i \, \mathsf{P}_k\biggl[ \pp{g}{x} \biggr](x,y) \, \hat{p}(x,y)
=
\mathsf{P}_{2k-1}[f](x,y) + \hat{r}(x,y)
\end{align}
for all $(x,y) \in [-1,1]^2$,
since both sides of (\ref{levnum:phatfun}) are bivariate polynomials of degree
less than $2k-1$ which agree on the $(2k-1) \times (2k-1)$ tensor product
Chebyshev grid.
Finally, we combine (\ref{levnum:gxcheb}), (\ref{levnum:fcheb}),
(\ref{levnum:pnorm}) and (\ref{levnum:errnorm}) in order to conclude that
\begin{align}
\label{levnum:pdehigherr}
\begin{split}
&
\left|
\pp{\hat{p}}{x}(x,y) + i \pp{g}{x}(x,y) \hat{p}(x,y) - f(x,y)
\right|\\
&\leq
\left\|
\hat{r}
+
\left( \mathsf{P}_{2k-1}[f] - f \right)
+
i \hat{p}
\left( \pp{g}{x} - \mathsf{P}_k\biggl[ \pp{g}{x} \biggr] \right)
\right\|_{L^\infty([-1,1]^2)}\\
&\lesssim
\varepsilon \, C(W_0) \, \frac{G_1}{G_0}
\left\| f \right\|_{L^\infty([-1,1]^2)}
\end{split}
\end{align}
for all $(x,y) \in [-1,1]^2$.

To complete our analysis for the case in which $G_0 > 0$, we let
\begin{align}
I
=
\int_{-1}^{1} \int_{-1}^{1}
f(x,y) \exp(i g(x,y))
\, \mathrm{d}x \, \mathrm{d}y
\end{align}
be the true value of the oscillatory integral we aim to evaluate using the Levin
method, and consider the estimate $\hat{I}$ of $I$ given by
\begin{align}
\label{levnum:estimate}
\hat{I}
=
\int_{-1}^{1} \hat{p}(1,y) \exp(i g(1,y)) \, \mathrm{d}y
-
\int_{-1}^{1} \hat{p}(-1,y) \exp(i g(-1,y)) \, \mathrm{d}y
,
\end{align}
where we assume that the univariate integrals are computed exactly.
We observe that
\begin{align}
\hat{I}
=
\int_{-1}^{1} \int_{-1}^{1}
\pp{}{x}( \hat{p}(x,y) \exp(i g(x,y)) )
\, \mathrm{d}x \, \mathrm{d}y
,
\end{align}
and it follows from this and (\ref{levnum:pdehigherr}) that
\begin{align}
\label{levnum:highinterr}
\begin{split}
| \hat{I} - I |
&\leq
\int_{-1}^{1} \int_{-1}^{1}
\left|
\pp{\hat{p}}{x}(x,y) + i \pp{g}{x}(x,y) \hat{p}(x,y) - f(x,y)
\right|
\mathrm{d}x \, \mathrm{d}y\\
&\lesssim
\varepsilon \, C(W_0) \, \frac{G_1}{G_0}
\left\| f \right\|_{L^\infty([-1,1]^2)}
.
\end{split}
\end{align}
We note that the ratio $G_1 / G_0$ is small whenever $g_x$ does not vary in
magnitude too much over the domain, which is a reasonable assumption for an
adaptive integration scheme.
We conclude that (\ref{levnum:highinterr}) shows that the absolute error in the
estimate is bounded independently of the magnitude of $g_x$, provided that $g_x$
does not vary too much over the domain.

\subsection[Case II: G1 <= 1/4]{Case II: $G_1 \leq 1/4$}

We now move on to the case in which $G_1 \leq 1/4$.
We invoke Theorem~\ref{leveq:thmlow} to see that there exists a bandlimied
function $p_b$ such that
\begin{align}
\begin{split}
&\left|
\pp{p_{b}}{x}(x,y) + i \pp{g}{x}(x,y) p_{b}(x,y) - f(x,y)
\right|\\
&\leq
3 \varepsilon \left( 1 + \frac{G_1}{1 - 2 G_1} \right)
\left\| f \right\|_{L^\infty([-1,1]^2)}
\end{split}
\end{align}
for all $(x,y) \in [-1,1]^2$, and both
\begin{align}
\left\| p_{b} \right\|_{L^\infty([-1,1]^2)}
&\leq
\frac{2}{1 - 2 G_1}
\left\| f \right\|_{L^\infty([-1,1]^2)}
\end{align}
and
\begin{align}
\left\| \pp{p_{b}}{x} \right\|_{L^\infty([-1,1]^2)}
&\leq
4 \left( 1 + \frac{G_1}{1 - 2 G_1} \right)
\left\| f \right\|_{L^\infty([-1,1]^2)}
.
\end{align}
From (\ref{leveq:defn}) and the assumption that $G_1 \leq 1/4$, we see that the
maximum bandlimit of $p_b$ is bounded, which implies that the Chebyshev
coefficients of $p_b$ are bounded by a rapidly decaying function which is
independent of $G_1$.
We can thus choose an integer $k$ independently of $G_1$ such that
(\ref{levnum:pbcheb}) through (\ref{levnum:fcheb}) hold.
Proceeding as in the previous subsection, we see that
\begin{align}
\begin{split}
&
\left|
\pp{}{x}\mathsf{P}_k[p_b](x,y)
+
i \, \mathsf{P}_k\biggl[ \pp{g}{x} \biggr](x,y) \, \mathsf{P}_k[p_b](x,y)
-
\mathsf{P}_{2k-1}[f](x,y)
\right|\\
&\lesssim
\varepsilon \left( 1 + \frac{G_1}{1 - 2 G_1} \right)
\left\| f \right\|_{L^\infty([-1,1]^2)}
,
\end{split}
\end{align}
for all $(x,y) \in [-1,1]^2$.
If we define $\bm{\mathsf{p_b}}$, $\bm{\mathsf{f}}$, $\bm{\mathsf{r}}$ and
$\bm{\mathsf{A}}$ as before, then we have that
\begin{align}
\| \bm{\mathsf{p_b}} \|_2
\lesssim
\frac{1}{1 - 2 G_1}
\left\| f \right\|_{L^\infty([-1,1]^2)}
,
\end{align}
and
\begin{align}
\| \bm{\mathsf{r}} \|_2
\lesssim
\varepsilon \left( 1 + \frac{G_1}{1 - 2 G_1} \right)
\left\| f \right\|_{L^\infty([-1,1]^2)}
.
\end{align}
We again solve the linear system
\begin{align}
\bm{\mathsf{A}} \bm{\mathsf{p_b}} = \bm{\mathsf{f}}
\end{align}
via a truncated singular value decomposition which has been truncated at
precision $\varepsilon \|\bm{\mathsf{A}}\|_2$, and, by Theorem~\ref{tsvd:thm},
we see that the obtained solution $\hat{\bm{\mathsf{p}}}$ satisfies
\begin{align}
\| \hat{\bm{\mathsf{p}}} \|_2
\lesssim
\frac{1}{1 - 2 G_1}
\left\| f \right\|_{L^\infty([-1,1]^2)}
,
\end{align}
and
\begin{align}
\| \bm{\mathsf{A}} \hat{\bm{\mathsf{p}}} - \bm{\mathsf{f}} \|_2
\lesssim
\varepsilon \, \frac{1}{1 - 2 G_1}
\left\| f \right\|_{L^\infty([-1,1]^2)}
.
\end{align}
Defining the polynomials $\hat{p}$ and $\hat{r}$ as before, we conclude that
\begin{align}
\left|
\pp{\hat{p}}{x}(x,y) + i \pp{g}{x}(x,y) \hat{p}(x,y) - f(x,y)
\right|
\lesssim
\varepsilon \, \frac{1}{1 - 2 G_1}
\left\| f \right\|_{L^\infty([-1,1]^2)}
\end{align}
for all $(x,y) \in [-1,1]^2$.
It then follows from this that the absolute error in the estimate $\hat{I}$ of
$I$ defined by (\ref{levnum:estimate}) satisfies
\begin{align}
\label{levnum:lowinterr}
| \hat{I} - I |
\lesssim
\varepsilon \, \frac{1}{1 - 2 G_1}
\left\| f \right\|_{L^\infty([-1,1]^2)}
.
\end{align}
Since $G_1 \leq 1/4$, we have that $1/2 \leq 1 - 2G_1 \leq 1$, so the constant
in (\ref{levnum:lowinterr}) is small.

\section{Algorithm description}

\label{section:algorithm}

In this section, we describe an adaptive delaminating Levin method
for the numerical evaluation of the integral
\begin{align}
\label{algorithm:oscint}
\int_{-1}^{1} \int_{-1}^{1} f(x,y) \exp(i g(x,y)) \, \mathrm{d}x \, \mathrm{d}y
.
\end{align}

We recall that the rigorous analysis presented in the preceding section
concerns a spectral collocation method in which the solution to the PDE
(\ref{levnum:pde}) is represented by a bivariate polynomial of degree less than
$k$, while the PDE is enforced on the $(2k-1) \times (2k-1)$ tensor product
Chebyshev grid.
This is necessary to establish (\ref{levnum:phatfun}) since the product
$\mathsf{P}_k[g_x] \hat{p}$ is of degree less than $2k-1$.
The algorithm we describe in this section only enforces the PDE on the
$k \times k$ tensor product Chebyshev grid.
This change enables the efficient solution of the PDE using the delaminating
Levin method, where the PDE is solved over the fibers
$\{ \bm{x}_{i,k}^{(j)} \}_{i=1}^{k}$, independently for each $1 \leq j \leq k$.
We justify this modification as follows.

When the delaminating Levin method is used adaptively, the integration domain is
typically subdivided until both $f$ and $g_x$ can be represented by polynomials
of degree much lower than $k$.
This is often the case for adaptive procedures since adaptive subdivision
typically leads to over-discretization of the inputs.
If $g_x$ is large and has been over-discretized, the obtained solution $\hat{p}$
will also be a polynomial of degree somewhat lower than $k$.
This is because the discretized PDE does not have a nontrivial null space when
$g_x$ is large, and so $\hat{p}$ will closely agree with the slowly-varying
solution $p_b$, which can be represented by a polynomial of low degree.
So in the case when $g_x$ is large, the product $g_x \hat{p}$ can often be
represented accurately by a polynomial of degree less than $k$.

On the other hand, when $g_x$ is small enough, the discretized null space is
nontrivial, but its elements can be represented by polynomials of low degree.
Accordingly, in this case, the product $g_x \hat{p}$ can also often be
represented by a polynomial of degree less than $k$.

Only when $g_x$ is of moderate size, not too large or too small, do we
expect $\hat{p}$ to be of degree close to $k-1$, which means that the product
$g_x \hat{p}$ is of degree larger than $k-1$ when $g_x$ is nontrivial.
Consequently, when the delaminating Levin method is employed, only those
subrectangles on which $g_x$ falls into a relatively narrow range of moderate
values will require additional subdivision.
Thus, we do not expect many additional subdivisions to occur.

We first describe the adaptive subdivision procedure before detailing the fixed
order delaminating Levin method.
The adaptive algorithm takes as input a subdivision tolerance parameter
$\varepsilon > 0$, an integer $k$ specifying the order of the Chebyshev spectral
collocation method used to discretize the PDE on each subrectangle considered,
and an external subroutine which evaluates the functions $f$ and $g$ at a
specified collection of points.
It maintains an estimated value $\textrm{val}$ of (\ref{algorithm:oscint}) and a
list of subrectangles.
Initially, the list of subrectangles contains only $[-1,1]^2$ and the value of
the estimate is set to zero.
While the list of subrectangles is nonempty, the adaptive algorithm repeats the
following steps:
\begin{enumerate}
\item Remove a subrectangle $[a,b] \times [c,d]$ from the list of subrectangles.

\item Compute an estimate $\textrm{val}_0$ of
\begin{align}
\textrm{val}_0
=
\int_{c}^{d} \int_{a}^{b} f(x,y) \exp(i g(x,y)) \, \mathrm{d}x \, \mathrm{d}y
\end{align}
using the delaminating Levin method.

\item Compute an estimate $\textrm{val}_i$ of
\begin{align}
\textrm{val}_i
=
\int_{R_i} f(x,y) \exp(i g(x,y)) \, \mathrm{d}x \, \mathrm{d}y
\end{align}
for $1 \leq i \leq 4$, where
\begin{align}
\label{algorithm:subrecs}
\begin{aligned}
&R_1
=
\biggl[ a, \frac{a + b}{2} \biggr] {\times} \biggl[ c, \frac{c + d}{2} \biggr]
,\,
&&R_2
=
\biggl[ \frac{a + b}{2}, b \biggr] {\times} \biggl[ c, \frac{c + d}{2} \biggr]
,\\
&R_3
=
\biggl[ a, \frac{a + b}{2} \biggr] {\times} \biggl[ \frac{c + d}{2}, d \biggr]
,\,
&&R_4
=
\biggl[ \frac{a + b}{2}, b \biggr] {\times} \biggl[ \frac{c + d}{2}, d \biggr]
,
\end{aligned}
\end{align}
using the delaminating Levin method.

\item If
$
\textrm{val}_0
-
(\textrm{val}_1 + \textrm{val}_2 + \textrm{val}_3 + \textrm{val}_4)
<
\varepsilon$,
then update the estimate $\textrm{val}$ by letting
$\textrm{val} = \textrm{val} + \textrm{val}_0$.
Otherwise, add the subrectangle $R_i$ to the list of subrectangles for all
$1 \leq i \leq 4$.
\end{enumerate}
When the list of subrectangles becomes empty, the adaptive algorithm terminates
and returns the estimate $\textrm{val}$ for (\ref{algorithm:oscint}).

We now set forth the delaminating Levin method, which estimates the value of
\begin{align}
\label{algorithm:subrect}
\int_{c}^{d} \int_{a}^{b} f(x,y) \exp(i g(x,y)) \, \mathrm{d}x \, \mathrm{d}y
,
\end{align}
where $[a,b] \times [c,d]$ is a bounded rectangle in $\R^2$.
It takes as input the domain $[a,b] \times [c,d]$, an integer $k$ which controls
the order of the Chebyshev spectral collocation method, a truncation tolerance
parameter $\varepsilon_0 > 0$, and an external subroutine for evaluating the
functions $f$ and $g$.
The delaminating Levin method proceeds as follows:
\begin{enumerate}
\item Use the external subroutine provided by the user to evaluate the functions
$f$ and $g$ at the nodes of the $k \times k$ tensor product Chebyshev grid
translated from $[-1,1]^2$ to $[a,b] \times [c,d]$, which we also denote by
\begin{align}
\left\{
\bm{x}_{i,k}^{(j)}
\right\}_{i,j=1}^{k}
,
\end{align}
where
\begin{align}
\bm{x}_{i,k}^{(j)}
=
\bm{x}_{i + k (j-1),k^2}
=
( t_{i,k}, t_{j,k} )
\end{align}
for $1 \leq i, j \leq k$,
so that $\{t_{i,k}\}_{i=1}^{k}$ and $\{t_{j,k}\}_{j=1}^{k}$ are the Chebyshev
nodes translated to $[a,b]$ and $[c,d]$, respectively.

\item For each $1 \leq j \leq k$,

\begin{enumerate}
\item Compute the approximate values
\begin{align}
\tilde{g_x\Bigl(\bm{x}_{1,k}^{(j)}\Bigr)},
\ldots,
\tilde{g_x\Bigl(\bm{x}_{k,k}^{(j)}\Bigr)}
\end{align}
of the partial derivative of $g$ with respect to $x$ using the formula
\begin{align}
\begin{pmatrix}
\tilde{g_x\Bigl(\bm{x}_{1,k}^{(j)}\Bigr)} \\
\tilde{g_x\Bigl(\bm{x}_{2,k}^{(j)}\Bigr)} \\
\vdots \\
\tilde{g_x\Bigl(\bm{x}_{k,k}^{(j)}\Bigr)}
\end{pmatrix}
=
\bm{\mathsf{D}}_k
\begin{pmatrix}
g\Bigl(\bm{x}_{1,k}^{(j)}\Bigr) \\
g\Bigl(\bm{x}_{2,k}^{(j)}\Bigr) \\
\vdots \\
g\Bigl(\bm{x}_{k,k}^{(j)}\Bigr)
\end{pmatrix}
.
\end{align}

\item Form the matrix
$
\bm{\mathsf{A}}_j
=
\bm{\mathsf{D}}_k + i \bm{\mathsf{G}}_j
$,
where $\bm{\mathsf{G}}_j$ is the $k \times k$ diagonal matrix with entries given
by
\begin{align}
(\bm{\mathsf{G}}_j)_{i,i}
=
\tilde{g_x\Bigl(\bm{x}_{i,k}^{(j)}\Bigr)}
\end{align}
for all $1 \leq i \leq k$.

\item Construct the singular value decomposition
\begin{align}
\bm{\mathsf{A}}_j
=
\bm{\mathsf{U}}_j
\begin{pmatrix}
\sigma_1^{(j)} &                &        &                \\
               & \sigma_2^{(j)} &        &                \\
               &                & \ddots &                \\
               &                &        & \sigma_k^{(j)} \\
\end{pmatrix}
\bm{\mathsf{V}}_j^{*}
\end{align}
of the matrix $\bm{\mathsf{A}}_j$.

\item Find the largest integer $1 \leq \nu_j \leq k$ such that
$\sigma_{\nu_j}^{(j)} \geq \varepsilon_0 \|\bm{\mathsf{A}}_j\|_2$.
If no such integer exists, set the values of
\begin{align}
\tilde{p\Bigl(\bm{x}_{1,k}^{(j)}\Bigr)},
\ldots,
\tilde{p\Bigl(\bm{x}_{k,k}^{(j)}\Bigr)}
\end{align}
to zero.

\item Otherwise, let
\begin{align}
\label{algorithm:linsolve}
\begin{pmatrix}
\tilde{p\Bigl(\bm{x}_{1,k}^{(j)}\Bigr)} \\
\tilde{p\Bigl(\bm{x}_{2,k}^{(j)}\Bigr)} \\
\vdots \\
\tilde{p\Bigl(\bm{x}_{k,k}^{(j)}\Bigr)}
\end{pmatrix}
=
\frac{b-a}{2}
\left( \bm{\mathsf{A}}_j \right)_{\nu_j}^{\dagger}
\begin{pmatrix}
f\Bigl(\bm{x}_{1,k}^{(j)}\Bigr) \\
f\Bigl(\bm{x}_{2,k}^{(j)}\Bigr) \\
\vdots \\
f\Bigl(\bm{x}_{k,k}^{(j)}\Bigr)
\end{pmatrix}
,
\end{align}
where $\left( \bm{\mathsf{A}}_j \right)_{\nu_j}^{\dagger}$ is the pseudoinverse
of the $\nu_j$-truncated singular value decomposition of $\bm{\mathsf{A}}_j$.
\end{enumerate}

\item 
\label{item:univar}
Let $\hat{p}$ be the bivariate polynomial of degree less than $k$ whose
values on the $k \times k$ tensor product Chebyshev grid are given by
\begin{align}
\hat{\bm{\mathsf{p}}}
=
\begin{pmatrix}
\tilde{p\bigl(\bm{x}_{1,k^2}\bigr)} &
\tilde{p\bigl(\bm{x}_{2,k^2}\bigr)} &
\cdots &
\tilde{p\bigl(\bm{x}_{k^2,k^2}\bigr)}
\end{pmatrix}
^T
.
\end{align}
Return the estimate
\begin{align}
\label{algorithm:bdryint}
\int_{c}^{d} \hat{p}(b,y) \exp(i g(b,y)) \, \mathrm{d}y
-
\int_{c}^{d} \hat{p}(a,y) \exp(i g(a,y)) \, \mathrm{d}y
\end{align}
computed using the univariate adaptive Levin method presented in \cite{Chen2024},
with subdivision tolerance parameter
$
\varepsilon^{\mbox{\tiny 1D}}
=
\beta
\max\{ \|\hat{\bm{\mathsf{p}}}\|_\infty / \|\bm{\mathsf{f}}\|_\infty, 1\}
\varepsilon
$
and truncation tolerance parameter
$\varepsilon_0^{\mbox{\tiny 1D}} = \varepsilon_0$, for some safety factor
$0 < \beta < 1$.

\end{enumerate}

This choice of subdivision tolerance parameter $\varepsilon^{\mbox{\tiny
1D}}$
for the univariate adaptive Levin
method in step~\ref{item:univar}
ensures that the univariate adaptive Levin converges. Since
the two-dimensional
adaptive Levin method 
is expected to converge for a relative tolerance of
$\varepsilon/\norm{\bm{\mathsf{f}}}_\infty$,
the univariate adaptive Levin method should also converge 
for the same relative tolerance. This means that, if the right hand side function
$\norm{\bm{\mathsf{\hat p}}}_\infty$ is large, which can happen if the Levin PDE
does not admit a slowly-varying solution for a particular subrectangle, 
then the absolute tolerance must be relaxed to
$\varepsilon \|\hat{\bm{\mathsf{p}}}\|_\infty / \|\bm{\mathsf{f}}\|_\infty$.
On the other hand, if $\norm{\bm{\mathsf{\hat p}}}_\infty$ is smaller than
$\norm{\bm{\mathsf{f}}}_\infty$, then the absolute tolerance parameter can 
be set to $\varepsilon$, since there is no benefit to evaluating the univariate
integrals to an absolute tolerance smaller than $\varepsilon$.

A distinctive feature of the delaminating Levin method described above is the
use of the univariate adaptive Levin method, as opposed to a univariate Levin
method with a fixed order.
When resonance points are present, but the PDE (\ref{levnum:pde}) admits a
smooth solution, a univariate Levin method of a fixed order can fail, which will
cause the algorithm to subdivide the subrectangle.
As shown in \cite{Chen2024}, the univariate adaptive Levin method can
efficiently evaluate univariate oscillatory integrals with stationary points,
requiring at worst a number of subdivisions that scales logarithmically with the
magnitude of the derivative of the phase function.
By using an adaptive method for the univariate integrals in
(\ref{algorithm:bdryint}), the subrectangles are only subdivided when the PDE
does not admit a smooth solution.
This is far more efficient than adaptively subdividing in two dimensions in
order to accurately compute the univariate integrals in
(\ref{algorithm:bdryint}) using a method of fixed order.

We note that the choice to delaminate along the $x$-direction is arbitrary, and
that the discussion of the Levin PDE and its numerical discretization presented
above can be trivially adapted if delamination is performed along the
$y$-direction instead.
In view of the bound (\ref{levnum:highinterr}), it may at first seem
advantageous to select the direction in which the ratio of the maximum to
minimum value of the corresponding partial derivative of $g$ is small.
Note however that evaluating the smaller of the two partial derivatives can be
an ill-conditioned procedure,
especially when the size of the other partial derivative is much larger.
As a result, we find that delaminating along the direction with the largest
partial derivative is most efficient, in that it requires the fewest
subdivisions.

\begin{remark}
In our implementation of the algorithm, we used a rank-revealing QR
decomposition in lieu of the truncated singular value decomposition to solve the
linear systems arising from the discretization of the PDE.
This was found to be about twice as fast overall and leads to no apparent loss
in accuracy.
\end{remark}

\begin{remark}
Both the delaminating Levin method and the univariate Levin method involve
solving linear systems of the form
\begin{align}
\label{algorithm:linsys}
( \bm{\mathsf{D}}_k + \bm{\mathsf{G}} )
\bm{\mathsf{p}}
=
\bm{\mathsf{f}}
,
\end{align}
where $\bm{\mathsf{D}}_k$ is the $k \times k$ Chebyshev differentiation matrix
and $\bm{\mathsf{G}}$ is a $k \times k$ diagonal matrix.
When $\bm{\mathsf{G}}$ is invertible, (\ref{algorithm:linsys}) can be solved
extremely rapidly via the iteration
\begin{align}
\bm{\mathsf{p}}_0 = \bm{0}
\quad
\text{and}
\quad
\bm{\mathsf{p}}_{n+1}
=
-\bm{\mathsf{G}}^{-1} \bm{\mathsf{D}}_k \bm{\mathsf{p}}_n
+
\bm{\mathsf{G}}^{-1} \bm{\mathsf{f}}
,
\end{align}
as long as there is some matrix norm $\| \cdot \|_p$ such that
$\| \bm{\mathsf{G}}^{-1} \bm{\mathsf{D}}_k \|_p \ll 1$.
Since
\begin{align}
\| \bm{\mathsf{G}}^{-1} \bm{\mathsf{D}}_k \|_\infty
\leq
\| \bm{\mathsf{G}}^{-1} \|_\infty \| \bm{\mathsf{D}}_k \|_\infty
,
\end{align}
this condition is straightforward to check in practice.
Furthermore, the iteration can be implemented in $\cO(k \log k)$ operations
using the discrete Chebyshev transform \cite{Mason2003}, which results in
substantial speedup when $k$ is large.
We note that this is often the case when Levin methods are used non-adaptively.
\end{remark}

\section{Numerical experiments}

\label{section:experiments}

In this section, we present the results of numerical experiments conducted to
illustrate the properties of the adaptive delaminating Levin method.
We implemented our algorithm in Fortran and compiled our code with
version~14.2.0 of the GNU Fortran compiler.
All experiments were performed on a desktop computer equipped with an Intel Core
i9-13900K processor and 64GB of memory.
No attempt was made to parallelize our code.

In all of our experiments,
the value of the parameter $k$, which determines the order of the bivariate
Chebyshev spectral collocation method, was taken to be $7$.
The subdivision tolerance parameter $\varepsilon$ was set to $10^{-12}$, and the
truncation tolerance parameter
$\varepsilon_0$ was set to $\beta_0 \varepsilon / \|\bm{\mathsf{f}}\|_\infty$
with $\beta_0 = 1/2$.
For computing the univariate boundary integrals, we used a $12$-point Chebyshev
spectral collocation method for the univariate adaptive Levin method, and set
the safety factor $\beta$ to $10^{-1}$.
To account for task scheduling and other vagaries of modern computing
environments, all reported times were obtained by averaging the cost of each
calculation over 100 runs.

In order to estimate the error in the results produced by the adaptive
delaminating Levin method, we compared the results of the adaptive delaminating
Levin method with the results produced by adaptive tensor product Gauss-Legendre
quadrature for evaluating integrals of the form
\begin{align}
\label{exp:integral}
\int_{-1}^{1} \int_{-1}^{1} f(x,y) \, \mathrm{d}x \, \mathrm{d}y,
\end{align}
except when explicit
formulas for the integrals were available.
Our implementation of the adaptive tensor product Gauss-Legendre quadrature is
written in Fortran and is quite standard.
It maintains a list of subrectangles, which is initialized with the single
rectangle $[-1,1]^2$, and a running tally of the value of the integral.
As long as the list of subrectangles is nonempty, the algorithm removes a
subrectangle $[a,b] \times [c,d]$ from the list, and compares the value of
\begin{align}
\label{exp:gauss0}
\int_{c}^{d} \int_{a}^{b} f(x,y) \, \mathrm{d}x \, \mathrm{d}y
\end{align}
as computed by a $10^2$-point tensor product Gauss-Legendre quadrature
rule to the value of the sum
\begin{align}
\sum_{i=1}^{4} \int_{R_i} f(x,y) \, \mathrm{d}x \, \mathrm{d}y
,
\end{align}
where the subrectangles $R_i$ are defined in (\ref{algorithm:subrecs}), and
where each integral is separately estimated with a $10^2$-point tensor
product Gauss-Legendre quadrature rule.
If the difference is larger than $\varepsilon$, where $\varepsilon$ is a
tolerance parameter specified by the user, then the subrectangles $R_i$ are
inserted into the list of subrectangles.
Otherwise, the value of (\ref{exp:gauss0}) is added to the running tally of the
integral (\ref{exp:integral}).
In all of our experiments, the tolerance parameter for the adaptive Gaussian
quadrature code was taken to be $\varepsilon = 10^{-14}$.

\subsection{Integrals involving elementary functions (with explicit formulas)}

\label{section:experiments:elem}

In the experiments described in this subsection, we consider the integrals
\begin{align}
\begin{split}
I_1(\lambda)
&=
\int_{0}^{1} \int_{0}^{1}
\cos(x+y) \exp(i \lambda (x + y + x^2 + y^2))
\, \mathrm{d}x \, \mathrm{d}y\\
&=
\frac{i \pi}{8 \lambda}
\sum_{\nu=0}^{1}
e^{i \left( 2 \nu - \frac{(1 + \lambda)^2}{2 \lambda} \right)}
\erf\left(
-\frac{(-1)^\nu + \lambda}{\sqrt{2 \lambda} (1 + i)}
,
-\frac{(-1)^\nu + 3\lambda}{\sqrt{2 \lambda} (1 + i)}
\right)^2
\end{split}
\end{align}
and
\begin{align}
\begin{split}
I_2(\lambda)
&=
\int_{0}^{2} \int_{0}^{2}
\frac{\exp(i \lambda (\arctan(x) + \arctan(y)))}{(1 + x^2) (1 + y^2)}
\, \mathrm{d}x \, \mathrm{d}y\\
&=
-\left( \frac{1 - \exp(i \lambda \arctan(2))}{\lambda} \right)^2
,
\end{split}
\end{align}
where $\erf(z_0,z_1) = \erf(z_1) - \erf(z_0)$ is the generalized error function.
The integral $I_1(\lambda)$ can be found in \cite{Levin1982,Li2009}, and
the integral $I_2(\lambda)$ can be computed via a change of variables.

\begin{figure}[t!]
\hfil
\includegraphics[width=.40\textwidth]{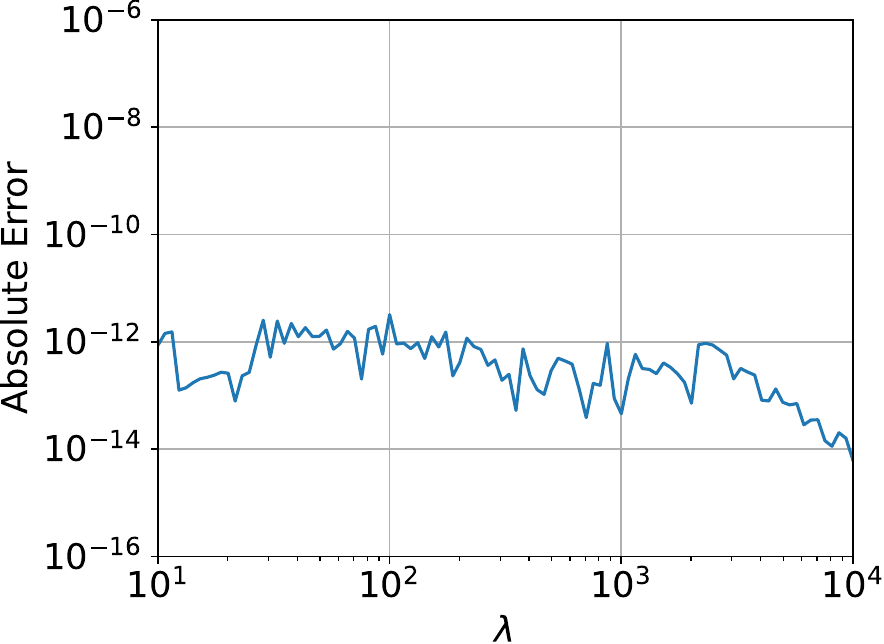}
\hfil
\includegraphics[width=.40\textwidth]{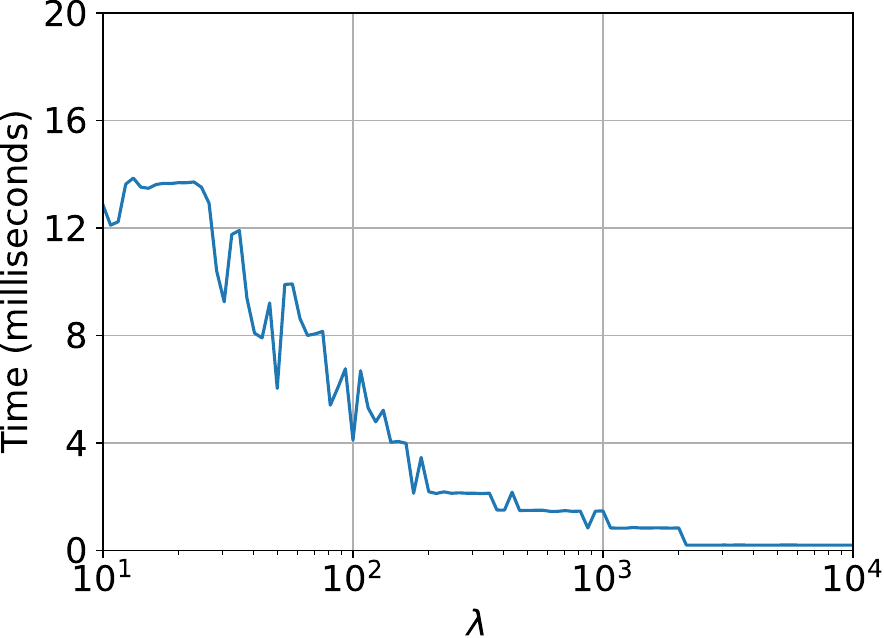}
\hfil

\hfil
\includegraphics[width=.40\textwidth]{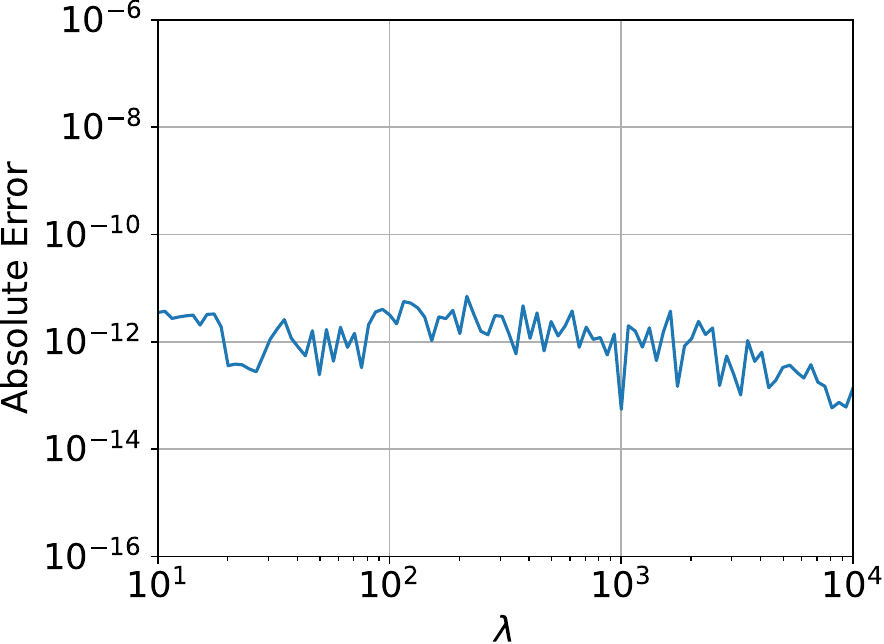}
\hfil
\includegraphics[width=.40\textwidth]{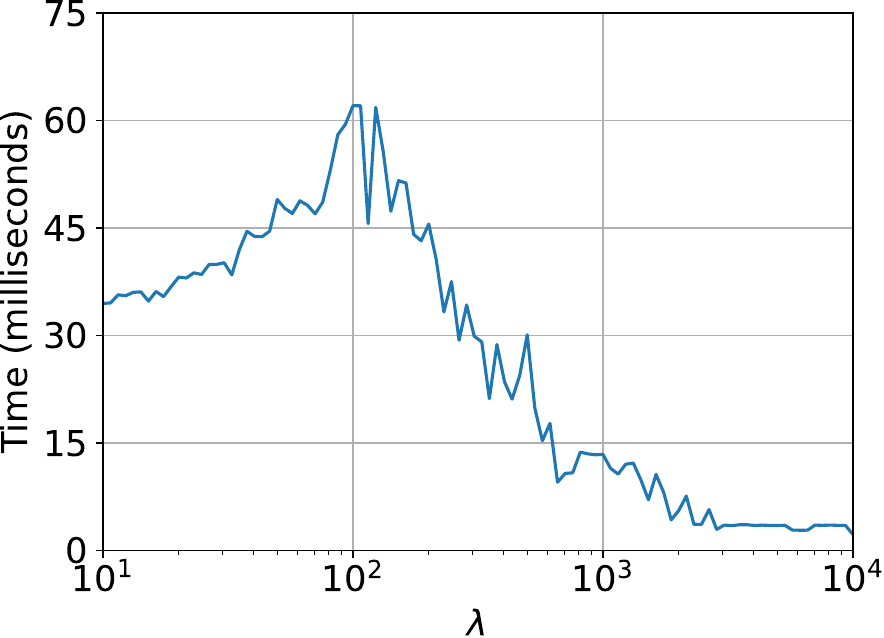}
\hfil

\caption{\small
The results of the experiments of Section~\ref{section:experiments:elem}.
The first row of plots pertains to the integral $I_1(\lambda)$, and the second
pertains to $I_2(\lambda)$.
In each row, the plot on the left gives the absolute error in the calculation of
the integral as a function of $\lambda$, and the plot on the right shows the
running time in milliseconds as a function of $\lambda$.
}
\label{figure:elemplots1}
\end{figure}

We sampled $\ell = 100$ equispaced points $x_1, \ldots, x_\ell$ in the interval
$[1,4]$, and used the adaptive delaminating Levin method to evaluate the
integrals $I_1(\lambda)$ and $I_2(\lambda)$ for each
$\lambda = 10^{x_1}, \ldots, 10^{x_\ell}$.
Figure~\ref{figure:elemplots1} gives the time taken to evaluate these
integrals and the absolute error in the obtained values.

\subsection{Integrals involving the Bessel functions (with explicit formulas)}

\label{section:experiments:bessel}

\begin{figure}[t!]
\hfil
\includegraphics[width=.40\textwidth]{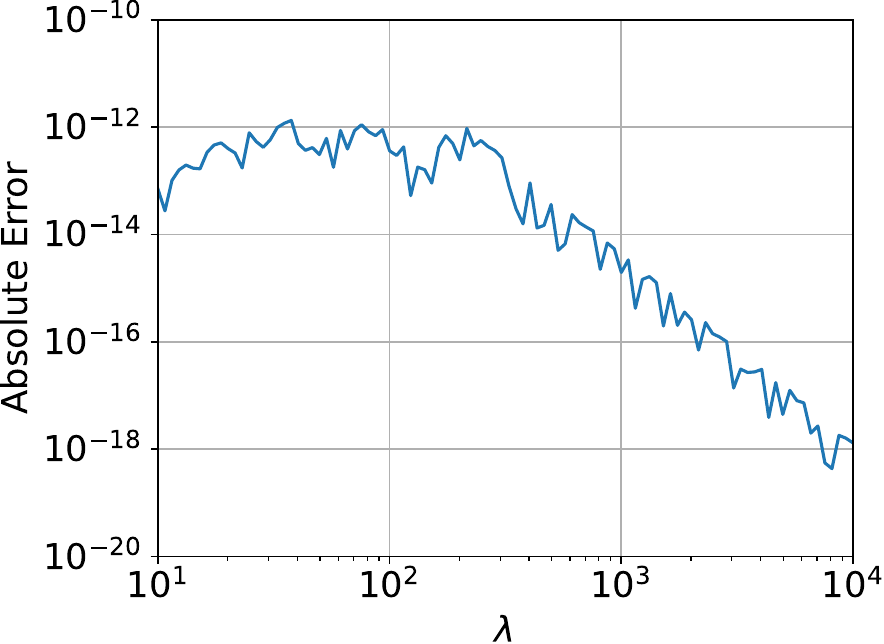}
\hfil
\includegraphics[width=.40\textwidth]{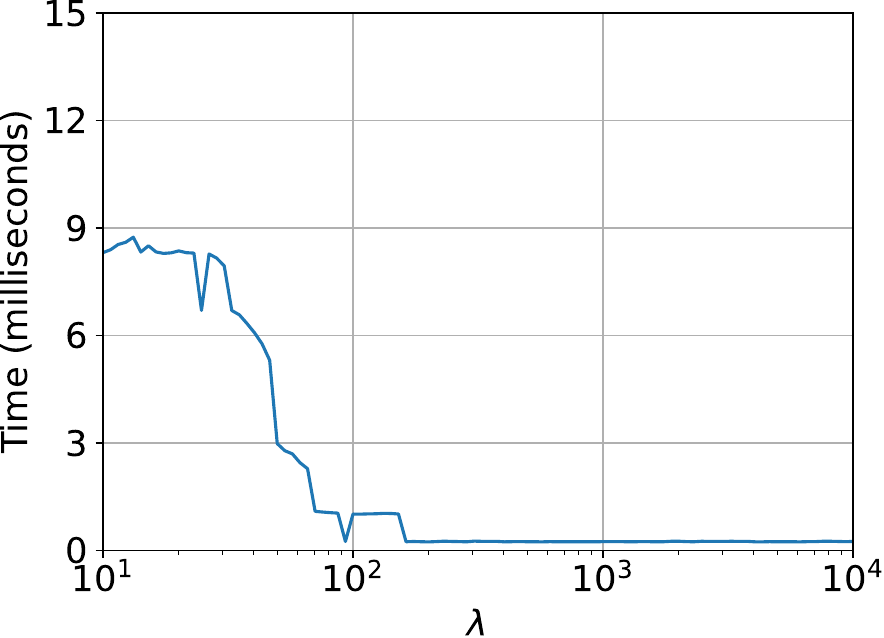}
\hfil

\hfil
\includegraphics[width=.40\textwidth]{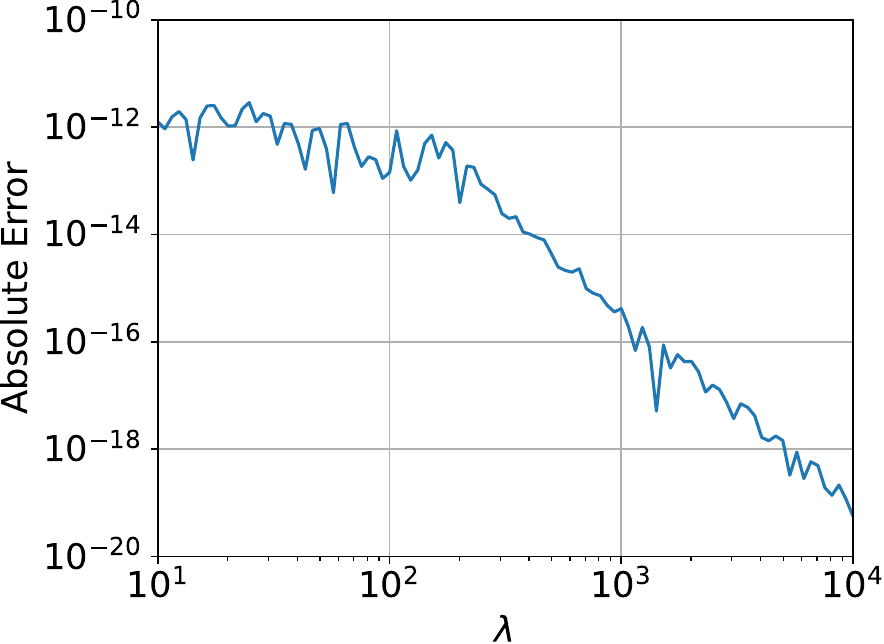}
\hfil
\includegraphics[width=.40\textwidth]{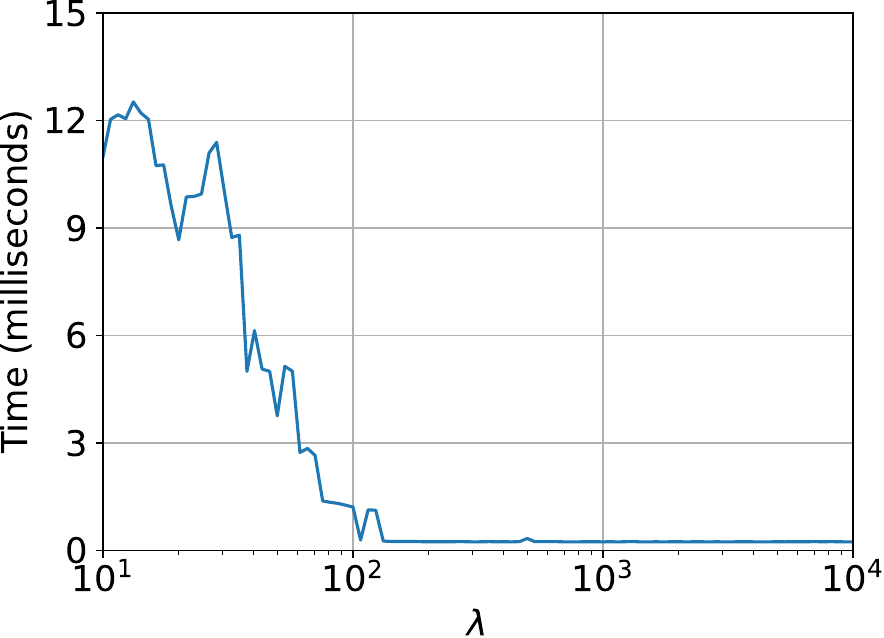}
\hfil

\caption{\small
The results of the experiments of Section~\ref{section:experiments:bessel}.
The first row of plots pertains to the integral $I_3(\lambda)$, and the second
pertains to $I_4(\lambda)$.
In each row, the plot on the left gives the absolute error in the obtained value
of the integral as a function of $\lambda$, and the plot on the right shows the
running time in milliseconds as a function of $\lambda$.
}
\label{figure:phaseplots1}
\end{figure}

In the experiments described in this subsection, we used the adaptive
delaminating Levin method to evaluate the integrals
\begin{align}
\begin{split}
I_3(\lambda)
&=
\int_{1}^{2} \int_{1}^{2}
x y H_0^{(1)}(\lambda x y)
\, \mathrm{d}x \, \mathrm{d}y\\
&=
\frac
{-H_0^{(1)}(\lambda) + 2 H_0^{(1)}(2 \lambda) - H_0^{(1)}(4 \lambda)}
{\lambda^2}
,
\end{split}
\end{align}
and
\begin{align}
\begin{split}
I_4(\lambda)
&=
\int_{1}^{2} \int_{1}^{2}
x^5 y^3 H_2^{(1)}(\lambda x y)
\, \mathrm{d}x \, \mathrm{d}y\\
&=
\frac
{H_4^{(1)}(\lambda) - 20 H_4^{(1)}(2 \lambda) + 64 H_4^{(1)}(4 \lambda)}
{\lambda^2}
.
\end{split}
\end{align}
Both of the above integrals can be found in \cite{Chen2009}.

The Hankel functions of the first kind $H_\nu^{(1)}$ are defined via the formula
\begin{align}
\label{experiments:hankel}
H_{\nu}^{(1)}(x)
=
J_\nu(x) + i Y_\nu(x)
,
\end{align}
where $J_\nu$ and $Y_\nu$ are the Bessel function of the first and second kinds
of order $\nu$.
For each $\nu \in \{0, 2\}$, we constructed a phase function
$\psi_\nu^{\mbox{\tiny bes}}$ for the normal form
\begin{align}
y''(x) + \left( 1 + \frac{\frac{1}{4} - \nu^2}{x^2} \right) y(x) = 0
\quad \text{in $(0,\infty)$}
\end{align}
of Bessel's differential equation using the algorithm of \cite{Bremer2023},
allowing the Bessel functions of the first and second kinds of order $\nu$ to
be computed as
\begin{align}
\label{experiments:besphase}
J_\nu(x)
=
\sqrt{\frac{2}{\pi x}}
\frac
{\sin\left( \psi_\nu^{\mbox{\tiny bes}}(x) \right)}
{\sqrt{\frac{\mathrm{d}}{\mathrm{d}x}\psi_\nu^{\mbox{\tiny bes}}(x)}}
\quad
\text{and}
\quad
Y_\nu(x)
=
-
\sqrt{\frac{2}{\pi x}}
\frac
{\cos\left( \psi_\nu^{\mbox{\tiny bes}}(x) \right)}
{\sqrt{\frac{\mathrm{d}}{\mathrm{d}x}\psi_\nu^{\mbox{\tiny bes}}(x)}}
.
\end{align}
The representation
\begin{align}
H_\nu^{(1)}(x)
=
\sqrt
{\frac{2}{\pi x \frac{\mathrm{d}}{\mathrm{d}x}\psi_\nu^{\mbox{\tiny bes}}(x)}}
\exp
\left(
i \left( \psi_\nu^{\mbox{\tiny bes}}(x) - \frac{\pi}{2} \right)
\right)
\end{align}
of the Hankel functions of the first kind was then used in conjunction with the
adaptive delaminating Levin method to evaluate the above integrals.
The construction of the phase functions took $2.79$ ms and $1.29$ ms for
$H_0^{(1)}$ and $H_2^{(1)}$, respectively.

We sampled $\ell = 100$ equispaced points $x_1, \ldots, x_\ell$ in the interval
$[1,4]$, and used the adaptive delaminating Levin method to evaluate the
above integrals for each $\lambda = 10^{x_1}, \ldots, 10^{x_\ell}$.
The results are shown in Figure~\ref{figure:phaseplots1}.
The reported times account for both the evaluation of the phase functions and
the execution of the adaptive delaminating Levin method.
Note that the time taken to construct the phase function is not included,
as it represents a one-time cost that is independent of $\lambda$.

\subsection{Behavior in the presence of a stationary point}

\label{section:experiments:stationary}

\begin{figure}[t!]
\hfil
\includegraphics[width=.40\textwidth]{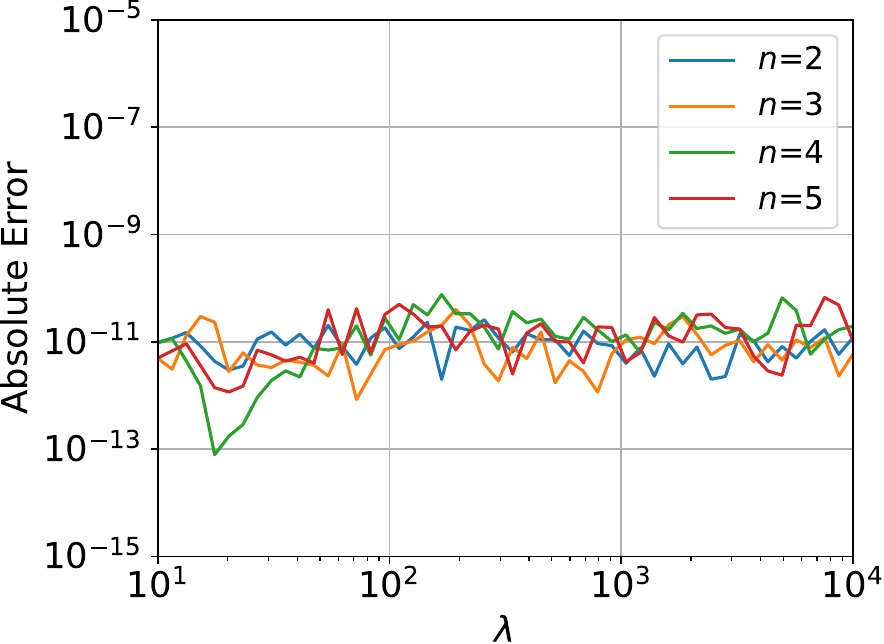}
\hfil
\includegraphics[width=.40\textwidth]{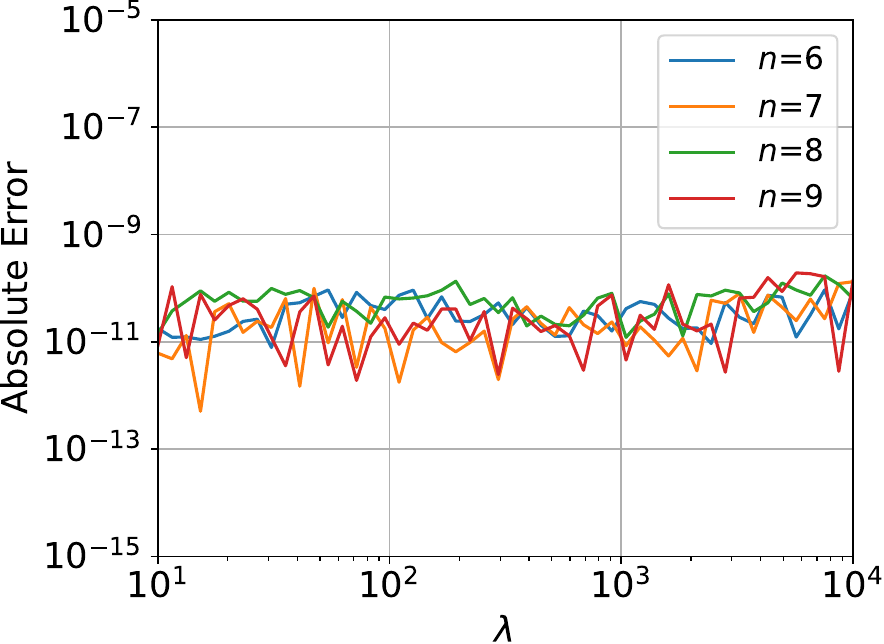}
\hfil

\hfil
\includegraphics[width=.40\textwidth]{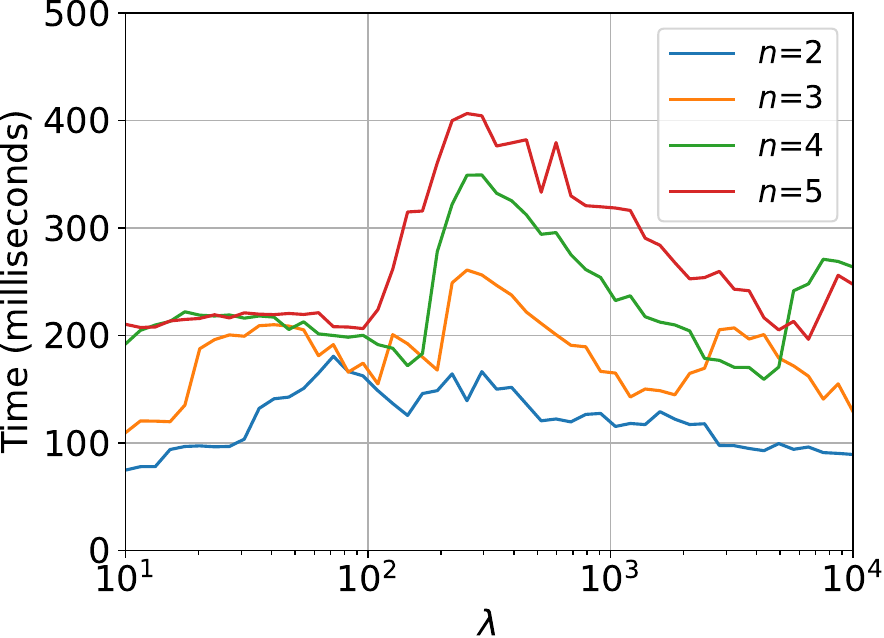}
\hfil
\includegraphics[width=.40\textwidth]{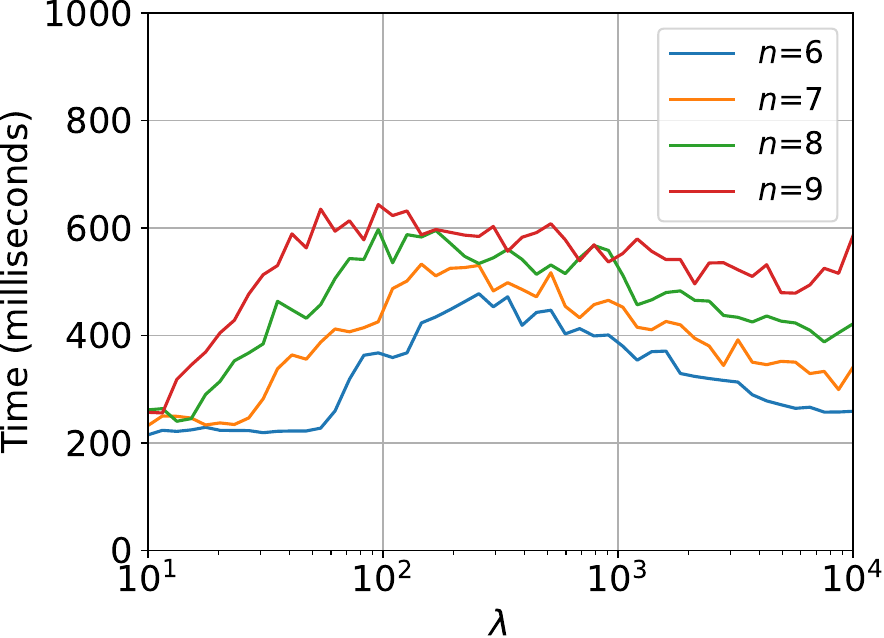}
\hfil

\caption{\small
The results of the experiments of Section~\ref{section:experiments:stationary}.
The first row of plots shows the absolute error in the calculated value of
$I_5(\lambda,n)$ as a function of $\lambda$ for the values of $n$ considered.
The second row of plots shows the running time of the adaptive delaminating
Levin method in milliseconds as a function of $\lambda$ for $n = 2$, $3$, $4$,
$5$, $6$, $7$, $8$, and $9$.
}
\label{figure:stationaryplots1}
\end{figure}

In the experiments described in this subsection, we consider the integral
\begin{align}
I_5(\lambda,n)
&=
\int_{-1}^{1} \int_{-1}^{1}
\frac{1}{1 + x^2 + y^2} \exp(i \lambda (x^n + y^n))
\, \mathrm{d}x \, \mathrm{d}y
\end{align}
in order to illustrate the behavior of the adaptive delaminating Levin method in
the presence of a stationary point.

We sampled $\ell = 50$ equispaced points $x_1, \ldots, x_\ell$ in the interval
$[1,4]$, and, for each $\lambda = 10^{x_1}, \ldots, 10^{x_\ell}$ and $n = 2$,
$3$, $4$, $5$, $6$, $7$, $8$, $9$, we evaluated $I_5(\lambda,n)$ using the
adaptive delaminating Levin method.
The results are shown in Figure~\ref{figure:stationaryplots1}.
We observe that the running time of our algorithm increases mildly with $n$ and,
for each fixed $n$, it is largely independent of $\lambda$ after peaking in the
moderate-frequency regime.

This behavior can be explained by the analysis presented in
Section~\ref{section:leveq}.
Since the phase function has a stationary point at the origin,
Theorem~\ref{leveq:thmlow} implies that the delaminating Levin method will yield
an accurate result on subrectangles of the form $[-\delta,\delta]^2$, provided
that the magnitude of $\nabla g$ is sufficiently small on the subrectangle.
Theorem~\ref{leveq:thmhigh} indicates that the delaminating Levin method will
also yield an accurate result on any subrectangle of $[-1,1]^2$ which is bounded
away from the stationary point, provided that the ratio of the maximum to
minimum value of the partial derivative of $g$ along the delamination direction
is small and that
\begin{align}
h(x',y')
=
f(\bm{u}^{-1}(x',y')) \det(D\bm{u}^{-1})(x',y')
\end{align}
can be approximated by a slowly-varying function with a small bandlimit over the
subrectangle.
The latter condition is more-or-less equivalent to the requirement that $h$ can
be represented by a Chebyshev expansion of small fixed order.
Thus, we expect the adaptive delaminating Levin method to subdivide the domain
until one or both of Theorems~\ref{leveq:thmhigh} and \ref{leveq:thmlow} apply
to the individual subrectangles.

Since $g_x(x,y) = \lambda n x^{n-1}$ and $g_y(x,y) = \lambda n y^{n-1}$, we see
that the condition
\begin{align}
\label{experiments:deltasize}
\delta < \left( \frac{C}{n \lambda} \right)^{\frac{1}{n-1}}
\end{align}
must be satisfied in order for $|g_x(x,y)|$ and $|g_y(x,y)|$ to be bounded by a
constant $C < 1$,
so we expect the adaptive delaminating Levin method to subdivide the domain into
at least
\begin{align}
\cO\left( \log\left( \frac{1}{\delta} \right) \right)
=
\cO\left( \frac{\log( n \lambda )}{n-1} \right)
\end{align}
subrectangles.
The algorithm will further subdivide the subrectangles outside this low
frequency region until, on each of the resulting subrectangles, the ratio of the
maximum to minimum value of the partial derivative of $g$ along the delamination
direction is small, and $h$ is accurately represented by a Chebyshev expansion
of small fixed order.

\begin{figure}
\hfil
\includegraphics[width=0.4\textwidth]{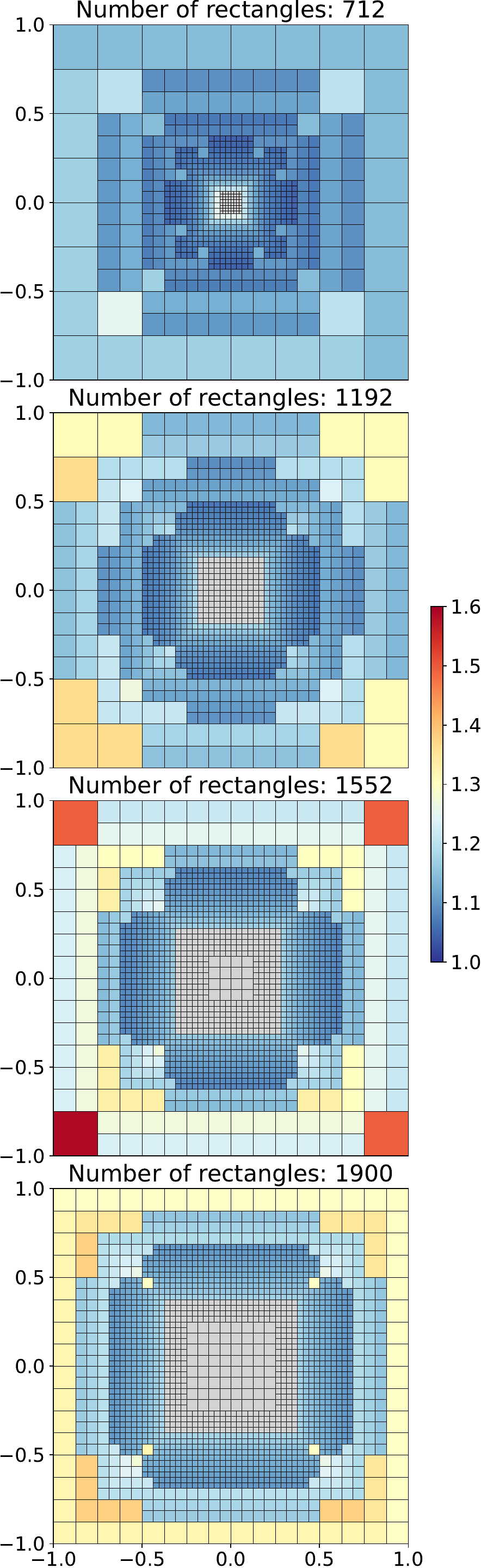}
\hfil

\caption{\small
The plots show the adaptive discretization of the domain produced by the
adaptive delaminating Levin method for computing the integral $I_5(300,n)$ with
$n=2$, $3$, $4$, and $5$.
The subrectangles are colored according to the ratio of the maximum to minimum
value of the partial derivative of $g$ along the delamination direction.
The subrectangles are grey if the relative maximum values of the partial
derivatives of $g$ are less than $1/4$.
}
\label{figure:stationaryrects1}
\end{figure}

Away from the stationary point, the low frequency region is enclosed by a
region in which the variation of the partial derivatives of $g$ is bounded,
so that, on each subrectangle $R$, the ratio $G_1 / G_0$ is small, where
$G_1 = \max_{(x,y) \in R} |g_v(x,y)|$ and $G_0 = \min_{(x,y) \in R} |g_v(x,y)|$,
and where $v$ is the delamination direction.
It is not hard to see that the subrectangles must be smallest in the
moderate frequency region bordering the low frequency region, and that they
become larger away from the stationary point.
As a result, the majority of the subrectangles lie on an annulus encircling the
low frequency region.
From (\ref{experiments:deltasize}), we see that the low frequency region
containing the stationary point expands with increasing $n$, so the area of the
annulus---and hence the number of subrectangles contained within---increases
with $n$.
It then follows that, when $\lambda$ is fixed, the total number of subrectangles
in the adaptive discretization of the domain must increase with $n$, which
explains why the running time of our algorithm increases with $n$.
Our analysis is consistent with the results shown in
Figure~\ref{figure:stationaryrects1}, which displays the adaptive discretization
of the domain produced by the adaptive delaminating Levin method for
$\lambda = 300$ and $n=2$, $3$, $4$, and $5$.
The subrectangles are colored according to the ratio of the maximum to minimum
value of the partial derivative of $g$ along the delamination direction, except
for subrectangles where the relative maximum values of the partial derivatives
of $g$ are less than $1/4$, which are shown in grey.

\subsection{Behavior in the presence of stationary and resonance points}

\label{section:experiments:difficult}

\begin{figure}[t!]
\hfil
\includegraphics[width=.40\textwidth]{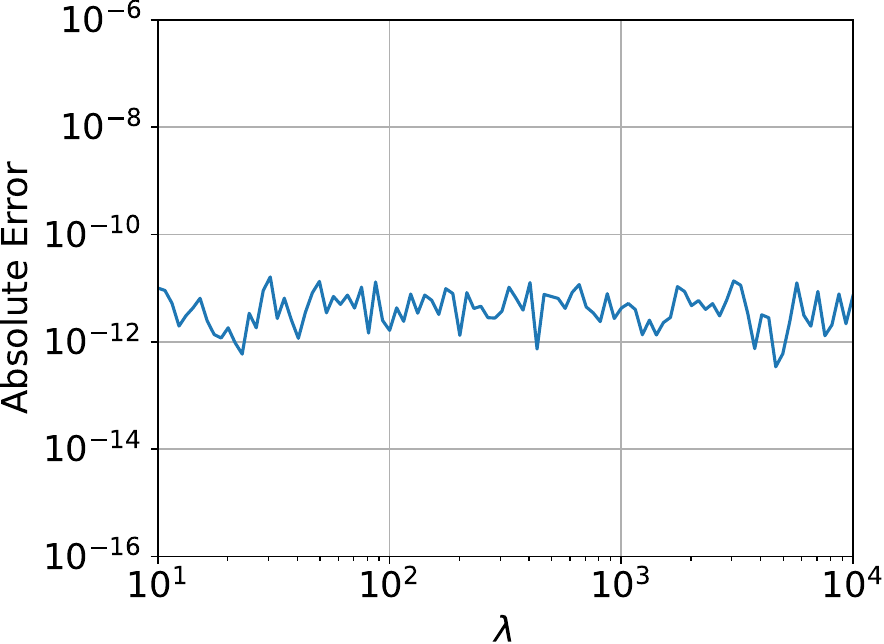}
\hfil
\includegraphics[width=.40\textwidth]{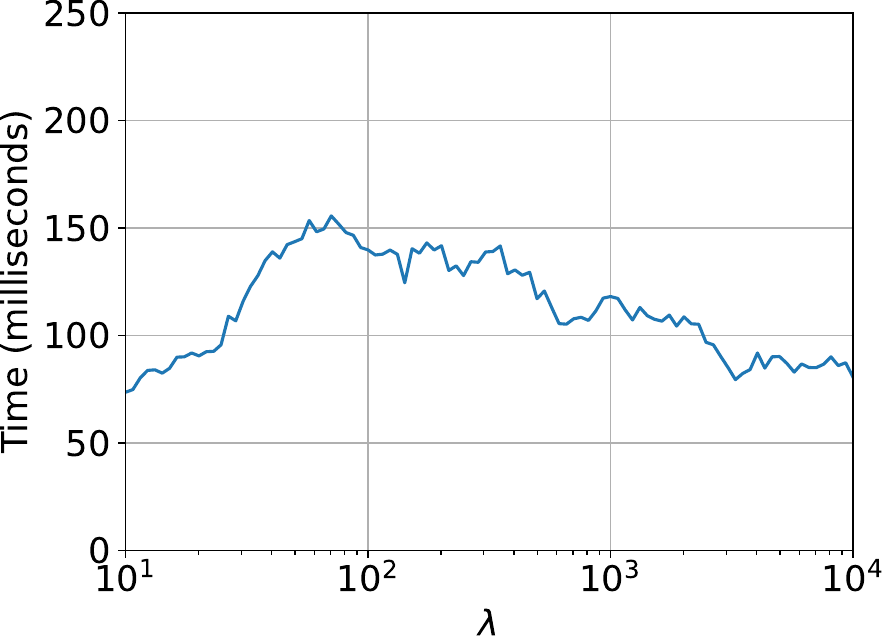}
\hfil

\hfil
\includegraphics[width=.40\textwidth]{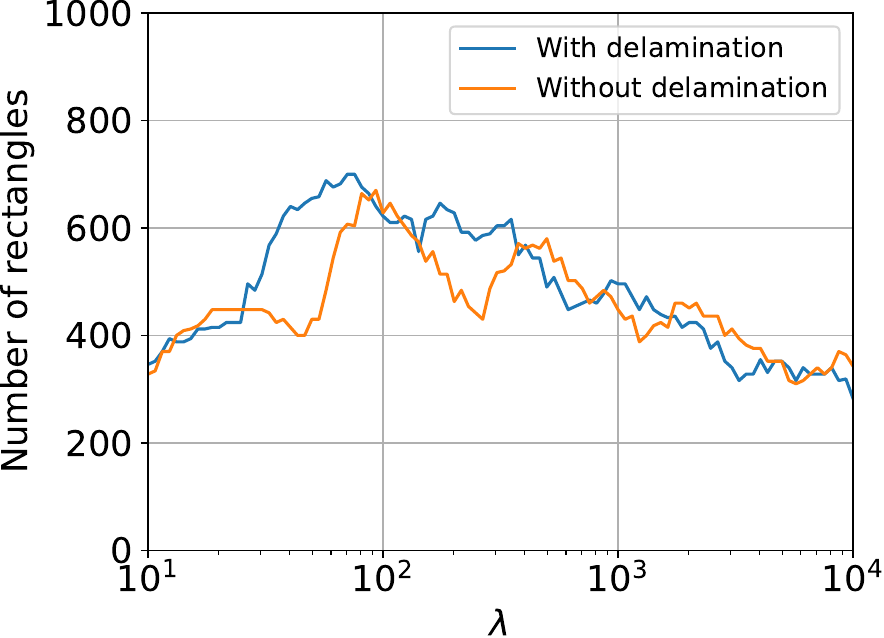}
\hfil
\includegraphics[width=.40\textwidth]{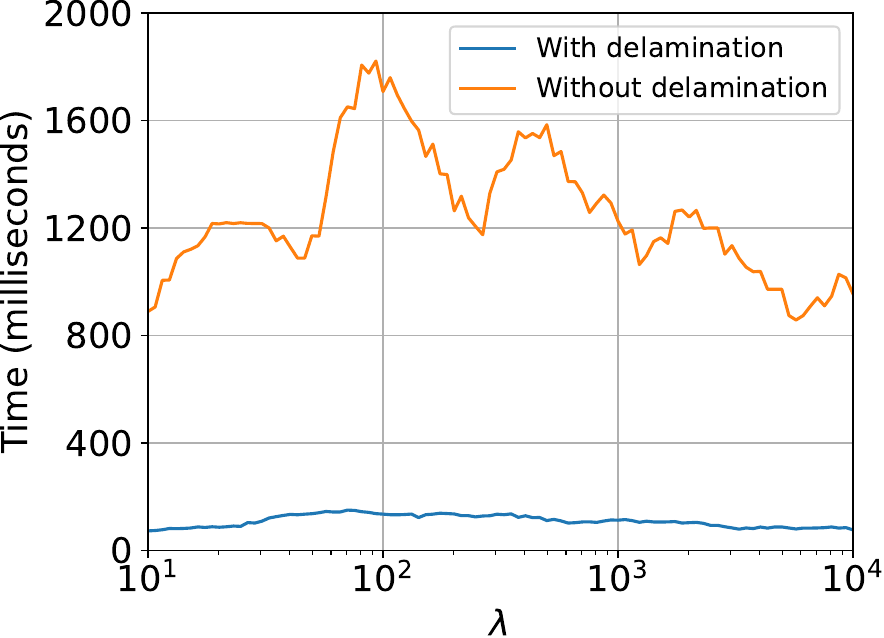}
\hfil

\caption{\small
The results of the experiments of Section~\ref{section:experiments:difficult}.
The first row of plots shows the absolute error in the calculated value of the
integral $I_6(\lambda)$ and the running time of the adaptive delaminating Levin
method in milliseconds as a function of $\lambda$.
The second row of plots compares both the number of rectangles in the adaptive
discretization of the domain and the running time, between the adaptive
delaminating Levin method and an adaptive method without delamination.
}
\label{figure:difficultplots1}
\end{figure}

In the experiments described in this subsection, we consider the integral
\begin{align}
I_6(\lambda)
&=
\int_{-1}^{1} \int_{-1}^{1}
(1 + xy) \exp(i \lambda (x^2 - xy - y^2))
\, \mathrm{d}x \, \mathrm{d}y
,
\end{align}
which is considered in~\cite{Huybrechs2007},
in order to understand the behavior of the adaptive delaminating Levin method
in the presence of both stationary and resonance points.
The phase function of $I_6(\lambda)$ is noteworthy
due to the stationary point at the origin, and the lines $y = 2x$ and
$y = -x/2$ where $g_x$ and $g_y$ vanish, respectively.

\begin{figure}
\centering
\mbox{\textsf{\small\hspace*{-1em} nondelaminating \hspace*{7em}
delaminating}}\medskip\\
\hfil
\hspace*{-3em}
\includegraphics[width=0.725\textwidth]{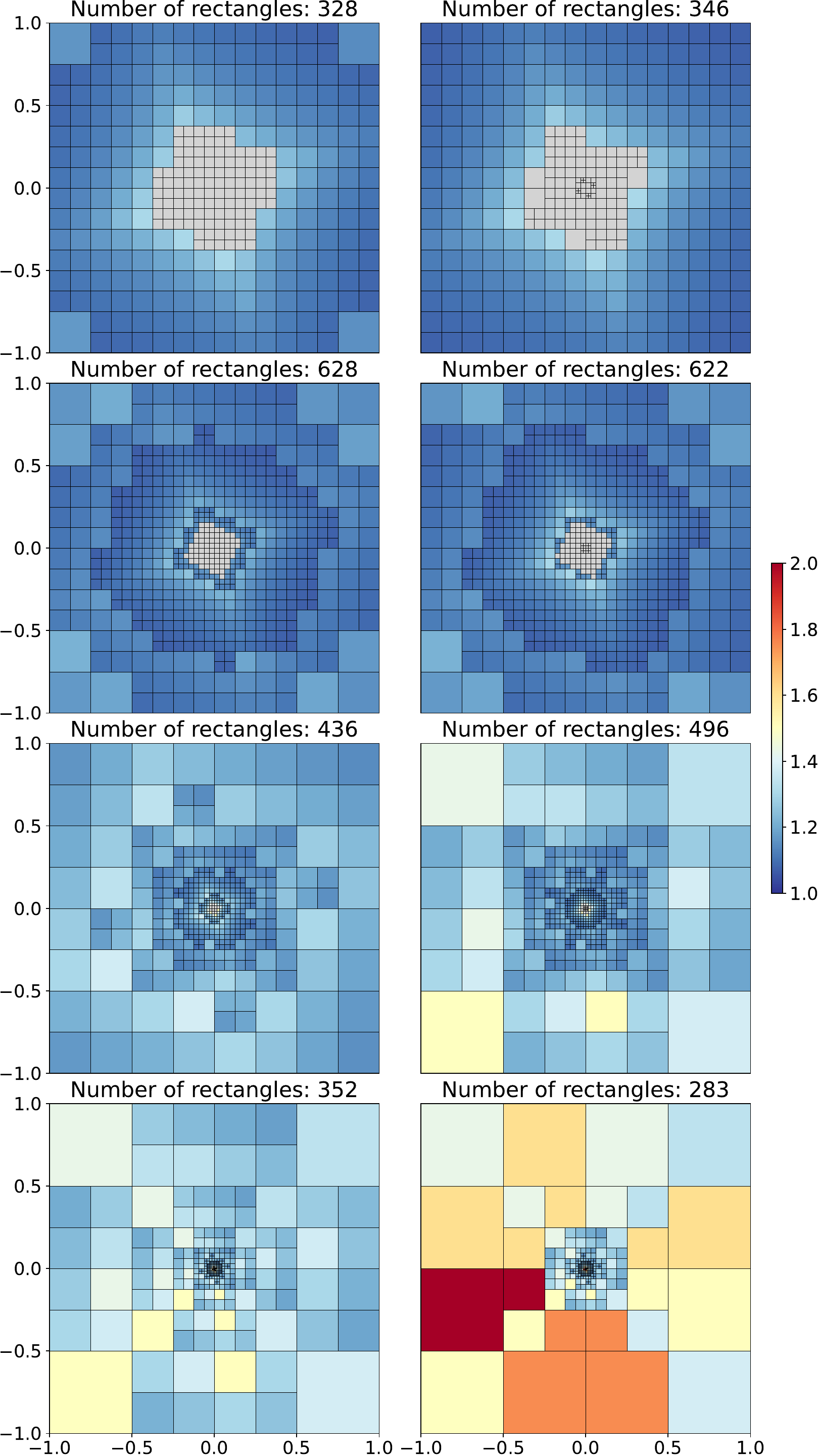}
\hfil

\caption{\small
The plots in the column on the left and right show the adaptive discretization
of the domain produced by the adaptive method that solves the rectangular
systems without delamination, and by the adaptive delaminating Levin method,
respectively, for computing the integral $I_6(\lambda)$ with $\lambda = 10^1$,
$10^2$, $10^3$, $10^4$.
The subrectangles are colored according to the ratio of the maximum to minimum
value of the partial derivative of $g$ along the delamination direction.
The subrectangles are grey if the relative maximum values of the partial
derivatives of $g$ are less than $1/4$.
}
\label{figure:difficultrects1}
\end{figure}

In the first set of experiments, we sampled $\ell = 100$ equispaced points
$x_1, \ldots, x_\ell$ in the interval $[1,4]$, and used the adaptive
delaminating Levin method to evaluate the integral $I_6(\lambda)$ for each
$\lambda = 10^{x_1}, \ldots, 10^{x_\ell}$.
The first row of Figure~\ref{figure:difficultplots1} presents the results.
We note that resonance points arise whenever the lines $y = 2x$ and $y = -x/2$
intersect the boundaries of a subrectangle.
Nevertheless, since the univariate Levin method is used to compute the boundary
integrals in (\ref{algorithm:bdryint}), the resonance points do not trigger any
additional subdivisions of subrectangles, as evidenced by
Figure~\ref{figure:difficultrects1}.

In a second set of experiments, we evaluated the above integral for each
$\lambda = 10^{x_1}, \ldots, 10^{x_\ell}$ using an adaptive method which solves
a rectangular linear system that enforces the collocation condition on the
$(2k - 1) \times (2k - 1)$ tensor product Chebyshev grid.
The results are shown in the second row of Figure~\ref{figure:difficultplots1}.
This approach follows the analysis presented in Section~\ref{section:levnum}
more closely, but the resulting rectangular linear system lacks a block
structure and must be solved monolithically.
Consequently, it incurs a significantly higher computational cost, especially as
the size of the collocation system increases with the value of $k$ and the
collocation grid size.
We observe that, while the adaptive delaminating Levin method produces slightly
more subdivisions in the moderate frequency regime, its overall running time is
substantially lower due to the smaller linear systems it solves.
Figure~\ref{figure:difficultrects1} shows the adaptive discretization of the
domain produced by the adaptive delaminating Levin method and the adaptive
method described above, for $\lambda = 10^1, 10^2, 10^3, 10^4$.
The subrectangles are colored according to the ratio of the maximum to minimum
value of the partial derivative of $g$ along the delamination direction.
The subrectangles where the relative maximum values of the partial derivatives
of $g$ are less than $1/4$ are colored grey.

\subsection{Behavior in the presence of many stationary points}

\label{section:experiments:stationary2}

\begin{figure}[t!]
\hfil
\includegraphics[width=.40\textwidth]{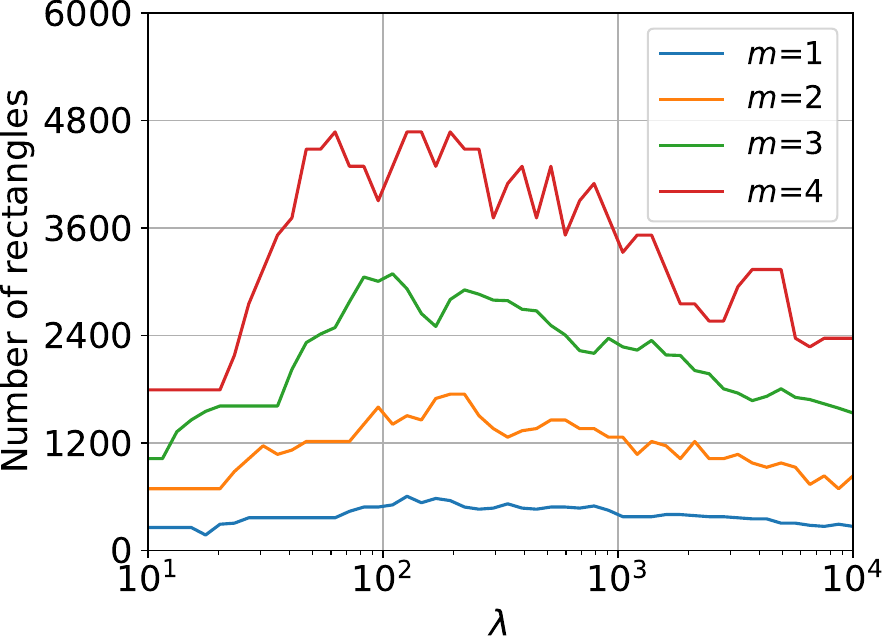}
\hfil
\includegraphics[width=.40\textwidth]{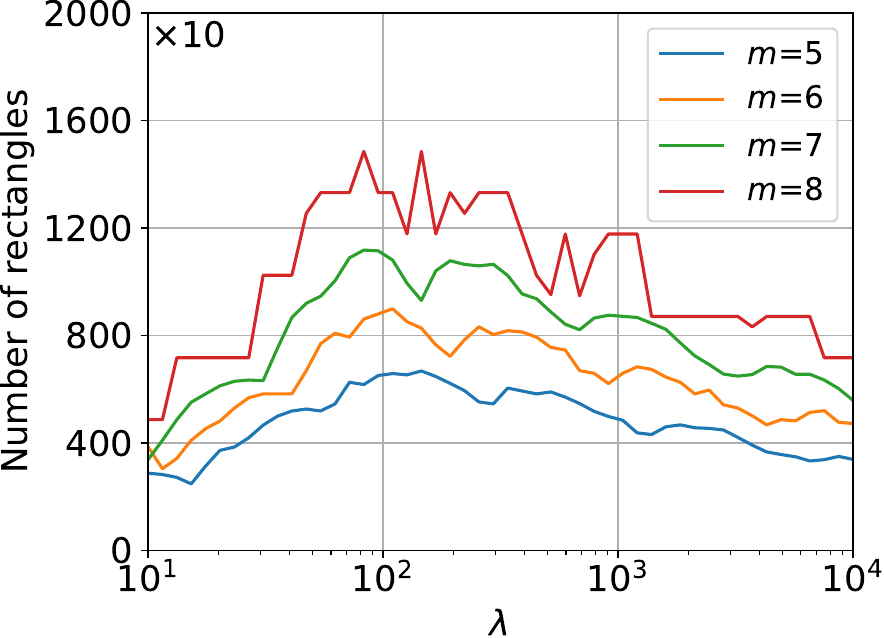}
\hfil

\hfil
\includegraphics[width=.40\textwidth]{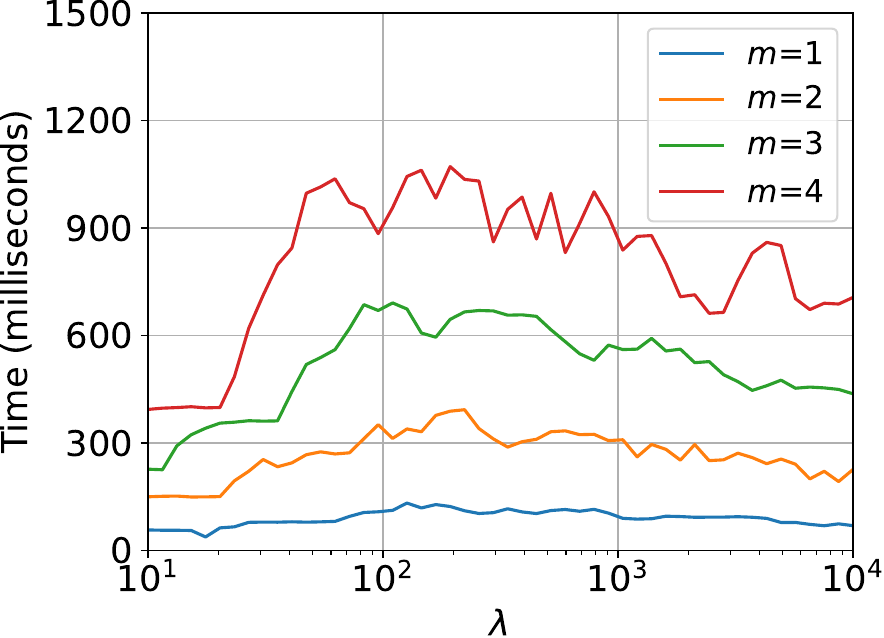}
\hfil
\includegraphics[width=.40\textwidth]{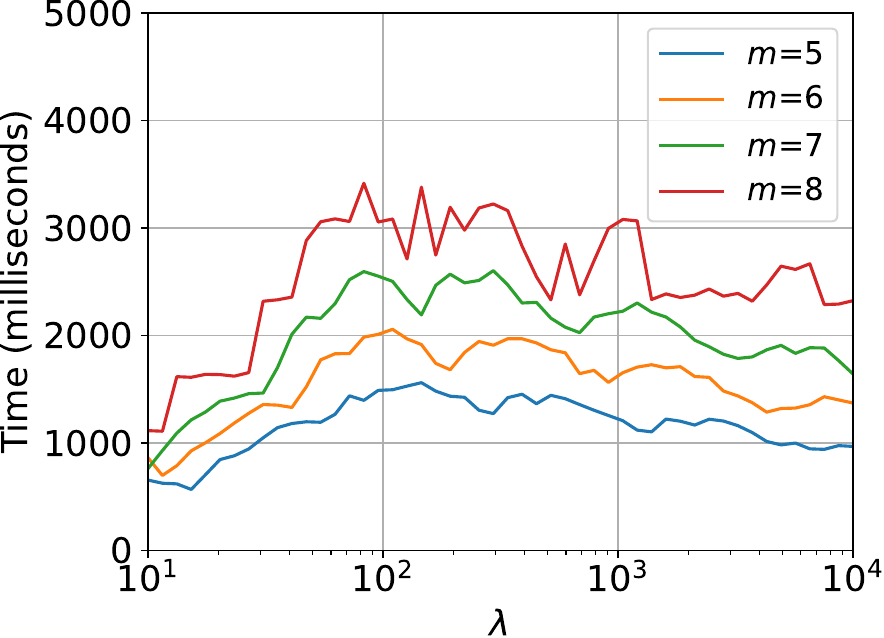}
\hfil

\hfil
\includegraphics[width=.40\textwidth]{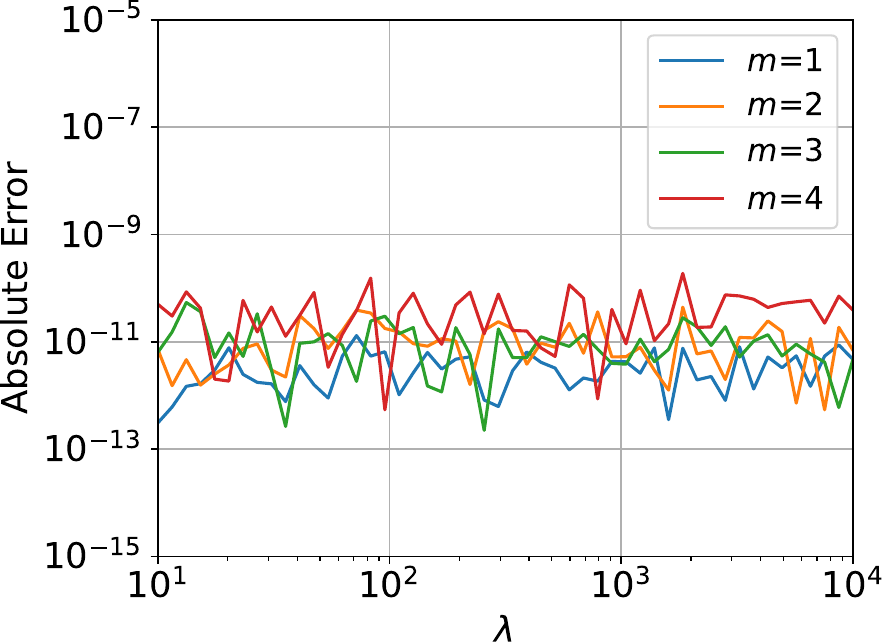}
\hfil
\includegraphics[width=.40\textwidth]{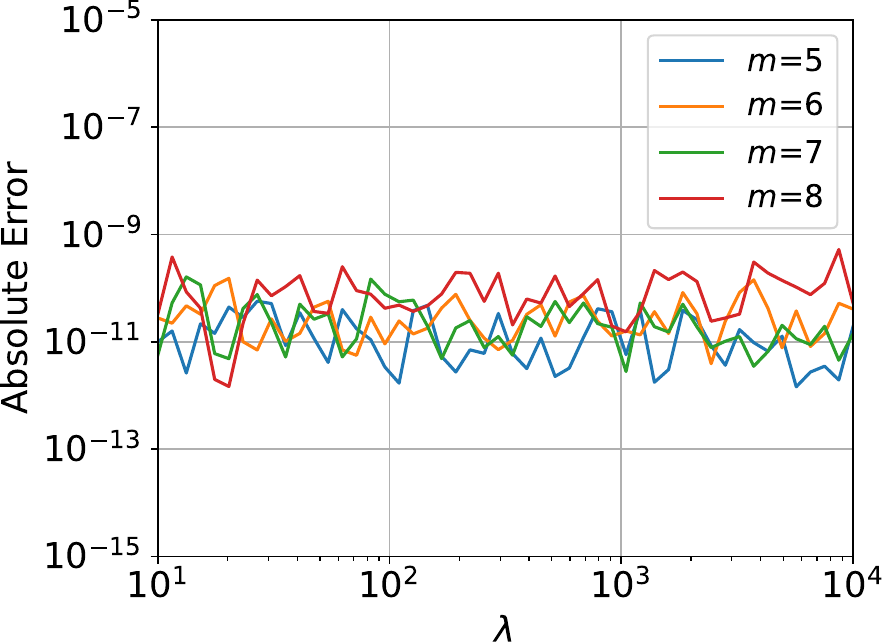}
\hfil

\caption{\small
The results of the experiments of Section~\ref{section:experiments:stationary2}.
The plots in the first row give the number of rectangles in the adaptively
determined discretization of $[0,1]^2$ used to compute the integral
$I_7(\lambda,m)$ as a function of $\lambda$ for the values of $m$ considered.
The plots in the second row show the runtime of the adaptive delaminating Levin
method as a function of $\lambda$ for the values of $m$ considered.
The third row of plots gives the absolute error in the calculated value of
$I_7(\lambda,m)$ as a function of $\lambda$ for $m = 1$, $2$, $3$, $4$, $5$, $6$,
$7$, and $8$.
}
\label{figure:stationaryplots2}
\end{figure}

In the experiments described in this subsection, we considered the integral
\begin{align}
I_7(\lambda,m)
&=
\int_{0}^{1} \int_{0}^{1}
\exp
\left(
i \lambda
\left(
\sin^2\left( \frac{\pi}{2} mx \right) + \sin^2\left( \frac{\pi}{2} my \right)
\right)
\right)
\mathrm{d}x \, \mathrm{d}y
,
\end{align}
which has $(m+1)^2$ stationary points in the unit square $[0,1]^2$.

We sampled $\ell = 50$ equispaced points $x_1, \ldots, x_\ell$ in the interval
$[1,4]$, and for each $\lambda = 10^{x_1}, \ldots, 10^{x_\ell}$ and $m = 1$, $2$,
$3$, $4$, $5$, $6$, $7$, $8$, we evaluated $I_7(\lambda,m)$ using the adaptive
delaminating Levin method.
The results are shown in Figure~\ref{figure:stationaryplots2}.
We observe that the running time of our algorithm grows linearly with the number
of stationary points, and is essentially independent of $\lambda$ for all values
of $m$ considered.

\section{Conclusions}

\label{section:conclusion}

We have shown that the Levin PDE admits a slowly-varying, approximate solution
across all frequency regimes, regardless of whether stationary and resonance
points are present, and that the adaptive delaminating Levin method,
when combined with a univariate adaptive Levin method,
rapidly and accurately evaluates bivariate oscillatory integrals over
rectangular domains.
We have also presented numerical experiments demonstrating the effectiveness of
our adaptive delaminating Levin method on a large class of bivariate oscillatory
integrals, including many with stationary and resonance points.

Our numerical method can be readily extended to evaluate bivariate oscillatory
integrals over general non-polytopal domains. When the integrand is
defined on a bounding box containing the domain, and
there exists a piecewise smooth parametrization of the boundary of the
domain, 
the Levin PDE can first be solved on the bounding box, and then
the univariate
adaptive Levin method can be used to evaluate the boundary integrals.
When the integrand cannot be evaluated outside the domain, but a smooth mapping
exists from a rectangle to the domain, the adaptive delaminating Levin method
can be applied on the pullback, with the Jacobian of the mapping appearing in
the integrand.
As a result, the method can be applied to any domain meshed with curved
quadrilateral elements, provided that the mapping functions for each element are
sufficiently smooth.

Finally, we note that both the delaminating Levin method and the theoretical
analysis presented in this paper are only weakly dependent on the dimensionality
of the ambient space, as tensor products are the primary theoretical tool.
As a result, our method can be easily generalized to evaluate volumetric
oscillatory integrals in dimensions three and higher.

\newpage

\bibliographystyle{hunsrt}
\bibliography{levin2d}

\newpage

\end{document}